\documentclass[amscd,amssymb,verbatim,12pt,letterpaper]{amsart}
\usepackage{latexsym}
\usepackage{euscript}

    \newcommand{\ie}{{\em i.e.}}
    \newcommand{\eg}{{\em e.g.}}
    \newcommand{\etc}{{\em etc.}}
    \newcommand{\etal}{{\em et al.}}

    \newtheorem{thm}{Theorem}[subsection]
    \newtheorem{prop}[thm]{Proposition}
    \newtheorem{lem}[thm]{Lemma}
    
    \newtheorem{cor}[thm]{Corollary}

    \newtheorem{hyp}[thm]{Hypothesis}   

    \theoremstyle{definition}

    \theoremstyle{remark}

    \newtheorem{rem}[thm]{Remark}

    \newcommand{\numbtext}[4]{\begin{equation}  \label{#1}
				  \left#2
				      \begin{array}{cc}
					  {} \parbox{4in}{#4} & {}
				      \end{array}
				  \right#3
			      \end{equation}}

    \newcommand{\fn}[3]{#1 \colon #2 \rightarrow #3}

    \newcommand{\Spec}{\operatorname{Spec}}

     \newcommand{\Hilb}{\operatorname{Hilb}}

    \newcommand{\pr}[3]{#1 \times_{#2} #3}

    \newcommand{\sheaf}[1]{\EuScript{#1}}
    \newcommand{\ideal}[1]{\mathfrak{#1}}

    \newcommand{\gf}{\mathsf{k}}   
    
    \newcommand{\Hn}{\textbf{\textup{H}}}
    
    \newcommand{\xlist}{\mathbf{x}}
    
    \newcommand{\jlist}{\mathbf{j}}
    \newcommand{\blist}{\mathbf{b}}
    \newcommand{\dlist}{\mathbf{d}}
    \newcommand{\elist}{\mathbf{e}}

    \newcommand{\lis}{\mathcal{S}} 

    \newcommand{\lcm}{\operatorname{lcm}}
    
    \newcommand{\nequiv}{\not \equiv}
    \newcommand{\nsim}{\not \sim}
    \newcommand{\ndiv}{\not \vert}

    \newcommand{\sa}{standard}
    \newcommand{\Sa}{Standard}
    \newcommand{\msa}{minimal \sa\ arrow}

    \newcommand{\rbox}{\EuScript{B}} 

\begin{document}

\title[Cotangent space at a monomial ideal]{The cotangent space at a monomial ideal of the Hilbert scheme of points of an affine space}

\author[M. Huibregtse]{Mark E. Huibregtse}

\address{Department of Mathematics and Computer Science\\
         Skidmore College\\
         Saratoga Springs, New York 12866}

\email{mhuibreg@skidmore.edu}

\date{\today} 


\begin{abstract}
	Let $\gf$ be an algebraically closed field.  We study the cotangent space of a point $t$ corresponding to a monomial ideal $I$ $\subseteq$ $\gf[x_1,\dots, x_r]$ in the Hilbert scheme of $n$ points of affine $r$-space (so $\dim_{\gf}(\gf[x_1,\dots, x_r]/I)$ $=$ colength of $I$ $=$ $n$).  Since $t$ lies in the closure of the locus corresponding to subschemes supported at $n$ distinct points of $\mathbb{A}^r_{\gf}$, one knows that the $\gf$-dimension of the cotangent space is always $\geq$ $r n$, and that $t$ is nonsingular if and only if the dimension equals $r n$.  We construct an explicit linearly independent set $\lis$ of cotangent vectors of size $r n$, and then explore conditions on $I$ under which $\lis$ either is or is not a basis of the cotangent space.  In particular, we give a condition on $I$ sufficient for $\lis$ to be a basis (equivalently, for $t$ to be nonsingular) that holds for every  monomial ideal in the case of $r$ $=$ $2$ variables, and that characterizes such ideals when $r$ $=$ $3$.  We also give an easily-checked condition on $I$ sufficient for $\lis$ not to be a basis.
\end{abstract}

\maketitle

\section{Introduction} \label{S:intro}

\subsection{Summary of results} \label{SS:overview}

Let $\gf$ be an algebraically closed field of any characteristic, and 
\[
    \mathbb{A}^r_{\gf} = \Spec(\gf[x_1, \dots, x_r]) = \Spec(\gf[\xlist])
\]
 the affine space of dimension $r$ over $\gf$.  The Hilbert scheme
\[
    \Hilb^{n}_{\mathbb{A}^r_{\gf}} = \Hn^n
\]
parameterizes the 0-dimensional closed subschemes
\[
    \Spec(\gf[\xlist]/I) \subseteq \mathbb{A}^r_{\gf}
\]
having length $n$, that is, 
\[
    \dim_{\gf}(\gf[\xlist]/I) = \text{\textbf{colength of }} I = n.
\]

In this paper we study the cotangent space of a point of $\Hn^n$ that corresponds to a monomial ideal $I$ (that is, $I$ is generated by monomials).  
If $I$ is a monomial ideal of colength $n$, and we let 
\[
    \beta = \{ \text{monomials } m \mid m \notin I \},
\]
then it is clear that $\beta$ is a $\gf$-basis of the quotient $\gf[\xlist]/I$; furthermore, $\beta$ has the property that for monomials $m_1$, $m_2$, 
\[
    m_1 \in \beta \text{ and } m_2\, |\, m_1 \Rightarrow m_2 \in \beta;
\]
we shall call any set of monomials $\beta$ with this property a \textbf{basis set} of monomials.  Let $U_{\beta}$ denote the (affine) open subscheme of $\Hn^n$ whose $\gf$-points $t$ are associated to ideals $I_t$ such that $\beta$ is a $\gf$-basis of the quotient $\gf[\xlist]/I_t$.  In particular, the monomial ideal $I$ $=$ $I_{\beta}$ that we started with has this property, so
\[
    t_{\beta} \in U_{\beta}, \text{ where } t_{\beta} \text{ is the point corresponding to } I_{\beta} = I_{t_{\beta}};
\]
we can therefore identify the cotangent space of $t_{\beta}$ with the $\gf$-vector space $M/M^2$, where $M$ is the maximal ideal of $t_{\beta}$ in the coordinate ring $R$ of $U_{\beta}$.  Since $t_{\beta}$ lies in the closure of the locus corresponding to closed subschemes supported at $n$ distinct points of $\mathbb{A}^r_{\gf}$, which is an $(r n)$-dimensional component of $\Hn^n$, we have that (Proposition \ref{prop:freedim})
\numbtext{eqn:proptext}{.}{.}
{
  \begin{itemize}
    \item $\dim_{\gf}(M/M^2)$ $\geq$ $r n$, and
    \item $t_{\beta}$ is nonsingular $\Leftrightarrow$
    $\dim_{\gf}(M/M^2)$ $=$ $r n$. 
  \end{itemize}
}  

Fortunately, we have a concrete description of $R$, based on the observation that for every point $t$ $\in$ $U_{\beta}$ and every monomial
\[
     x_1^{d_1}x_2^{d_2}\dots x_r^{d_r} = \xlist^{\dlist}\in I_{\beta},
\]
there is a unique polynomial
\[
       F_{\dlist}(t)\ =\ \xlist^{\dlist} - \sum_{\xlist^{\jlist} \in \beta} 
c^{\dlist}_{\jlist}(t) \cdot \xlist^{\jlist}\  \in\  I_t,\ \ c^{\dlist}_{\jlist}(t) \in \gf,
\]
since the quotient $\gf[\xlist]/I_t$ has $\gf$-basis $\beta$.  As $t$ ranges over $U_{\beta}$, the coefficients $c^{\dlist}_{\jlist}(t)$ define functions $c^{\dlist}_{\jlist}$ $\in$ $R$ that generate $R$ as a $\gf$-algebra; moreover, the relations among these functions can be completely described.  Note in particular that the point $t_{\beta}$ is the ``origin'' of $U_{\beta}$, the point at which all the functions $c^{\dlist}_{\jlist}$ vanish.

\begin{figure} 
\begin{picture}(374,96)
        \put(130,0){\line(1,0){24}}
        \put(130,0){\line(0,1){48}}
        \put(154,0){\line(0,1){48}}
        \put(130,48){\line(1,0){24}}

        \put(130,48){\line(2,1){24}}
        \put(154,48){\line(2,1){24}}
        \put(154,0){\line(2,1){24}}
        \put(178,12){\line(0,1){48}}

        \put(154,60){\line(1,0){24}}
        \put(154,60){\line(0,1){12}}
        \put(154,72){\line(1,0){48}}
        \put(154,72){\line(2,1){24}}
        \put(178,84){\line(1,0){48}}
        \put(202,72){\line(2,1){24}}
        \put(202,72){\line(0,-1){60}}
        \put(178,12){\line(1,0){24}}

        \put(202,12){\line(2,1){12}}
        \put(214,18){\line(1,0){36}}
        \put(214,18){\line(0,1){24}}
        \put(214,42){\line(2,1){12}}
        \put(214,42){\line(1,0){36}}

        \put(226,48){\line(0,1){36}}
        \put(226,48){\line(1,0){36}}

        \put(250,18){\line(2,1){12}}
        \put(250,18){\line(0,1){24}}

        \put(250,42){\line(0,-1){24}}
        \put(250,42){\line(2,1){12}}
        \put(262,48){\line(0,-1){24}}

        \put(130,0){\vector(-2,-1){12}}
        \put(112,6){$x_1$}
        \put(262,24){\vector(1,0){12}}
        \put(268,30){$x_2$}
        \put(178,84){\vector(0,1){12}}
        \put(160,90){$x_3$}

        \put(154,63){\vector(4,-1){94}}
        \put(136,60){$\xlist^{\blist}$}
        \put(142,62){\vector(1,0){10}}
        \put(266,40){$\xlist^{\jlist}$}
        \put(265,40){\vector(-1,0){14}}  
\end{picture}
\caption{A basis set $\beta$ (the monomials within the boxes) in $r$ = 3 variables with a non-\sa\ rigid arrow $c^{\blist}_{\jlist}$. The point $t_{\beta}$ $\in$ $\Hn^n$ is accordingly singular.  The arrow  is non-\sa\ because the vector $\jlist$ $-$ $\blist$ has $2$ negative components: the first and the third.} \label{fig:nonstrigid}
\end{figure}
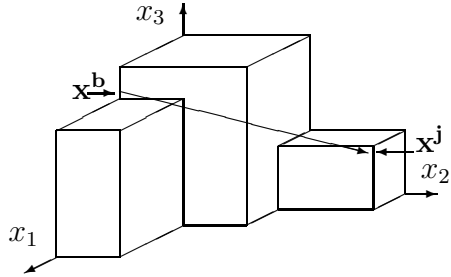

The maximal ideal $M$ $\subseteq$ $R$ that cuts out the point $t_{\beta}$ is therefore generated by the functions $c^{\dlist}_{\jlist}$.  Each of these functions can be viewed as an arrow pointing from $\xlist^{\dlist}$ to $\xlist^{\jlist}$ in the lattice of monomials obtained by identifying each monomial $\xlist^{\blist}$ with its $r$-tuple $\blist$ of exponents (for example, see Figures \ref{fig:nonstrigid} and \ref{fig:arroweg}).  It then turns out that translation of arrows (keeping the head inside $\beta$ and the tail outside $\beta$) corresponds to equivalence modulo $M^2$; more precisely (see Section \ref{sec:cotansp}),
\begin{itemize}
  \item If $c^{\dlist_1}_{\jlist_1}$ can be translated to $c^{\dlist_2}_{\jlist_2}$, then $c^{\dlist_1}_{\jlist_1}$ $\equiv$ $c^{\dlist_2}_{\jlist_2} \pmod{M^2}$.
  \item If $c^{\dlist_1}_{\jlist_1}$ $\equiv$ $0 \pmod{M^2}$, then $c^{\dlist_1}_{\jlist_1}$ can be translated so that its head exits $\beta$ across a hyperplane ($x_i$-degree $=$ $0$), and conversely.
  \item If $c^{\dlist_1}_{\jlist_1}$ $\equiv$ $c^{\dlist_2}_{\jlist_2}$ $\nequiv$ $0 \pmod{M^2}$, then $c^{\dlist_1}_{\jlist_1}$ can be translated to $c^{\dlist_2}_{\jlist_2}$.
\end{itemize}

Our first main result, Theorem \ref{thm:lindepconds}, states that if $S$ is a finite set of $d$ arrows such that no member of $S$ can be translated so that its head exits $\beta$ across a hyperplane ($x_i$-degree $=$ $0$), and no two distinct members of $S$ can be translated one to the other, then the $\gf$-span of $S$ in $M/M^2$ has dimension $d$, and conversely; in this case, we say that $S$ has \textbf{maximal rank} (mod ${M^2}$).  Using this criterion, we construct in Theorem \ref{thm:S} a  maximal rank (mod ${M^2}$) set $\lis$ of $r n$ arrows $c^{\dlist}_{\jlist}$ of a special kind that we call \textbf{minimal \sa} arrows: \emph{minimal} means that the tail $\xlist^{\dlist}$ is a minimal generator of the monomial ideal $I_{\beta}$, and \emph{standard} means that the vector $\jlist$ $-$ $\dlist$ has only one negative component.  (Since $\xlist^{\jlist}$ $\in$ $\beta$ and $\xlist^{\dlist}$ $\notin$ $\beta$, the vector of any arrow must always have at least one negative component, to avoid the contradiction $\xlist^{\dlist}$ $|$ $\xlist^{\jlist}$.)  Therefore, in view of (\ref{eqn:proptext}), 
\[
    t_{\beta} \text{ is nonsingular } \Leftrightarrow \lis \text{ spans } M/M^2 \Leftrightarrow \lis \text{ is a } \gf \text{-basis of } M/M^2.
\]
When any (and hence all) of these equivalent conditions holds, we say that $\beta$ is a \textbf{smooth} basis set.  (The set $\lis$ is in general not unique.)

\begin{rem} \label{rem:HaimanAck}
     The idea of arrow translation and its connection to equivalence modulo ${M^2}$  is due to Haiman \cite{Haiman:CN-HS} (in the case of two variables); it is a pleasure to acknowledge the inspiration of that beautiful article.
\end{rem}

The remainder of the paper explores conditions on $\beta$ sufficient to imply that $\beta$ either is or is not smooth.  The basic result is Theorem \ref{thm:necsufcond}, which in essence says that 
\begin{equation} \label{eqn:basres}
   \beta \text{ smooth } \Leftrightarrow 
      \left\{
        \begin{array}{l}
           \bullet\ c^{\dlist}_{\jlist} \text{ non-\sa\ } \Rightarrow c^{\dlist}_{\jlist} \equiv 0 \pmod{M^2}, \text{ and}\\
           \bullet\ c^{\dlist^{\prime}}_{\jlist^{\prime}} \text{ \sa\ and } \nequiv 0 \pmod{M^2} \Rightarrow c^{\dlist^{\prime}}_{\jlist^{\prime}} \\
                    \ \ \, \text{can be translated to an arrow in } \lis.
        \end{array} 
      \right.
\end{equation}

In Section \ref{SS:nonsmthcond}, we present a simple condition on $\beta$ that ensures the existence of non-\sa\ arrows that are \emph{not} $\equiv$ 0 modulo $M^2$.  One first checks that any minimal arrow $c^{\dlist}_{\jlist}$ with head $\xlist^{\jlist}$ \emph{maximal} in $\beta$ (for the partial ordering defined by divisibility of monomials) cannot be translated except to itself (Lemma \ref{lem:rigid}); we call such an arrow \textbf{rigid}.  It is then clear that any non-\sa\ rigid arrow will violate the first bulleted condition in (\ref{eqn:basres}), forcing $\beta$ to be non-smooth (Corollary \ref{cor:rigid}).  Figure \ref{fig:nonstrigid} illustrates such an arrow; a specific example is discussed in Section \ref{SS:nonsmootheg}.
\medskip

We go on to present three ways to construct smooth basis sets:

\begin{description}
    \item[Thickening]
         In Section \ref{S:thicken}, we show that a basis set $\beta$ in $r$ variables that is the ``thickening'' (\ref{eqn:liftdef}) of a smooth basis set $\beta_0$ in $r-1$ variables will be smooth (Theorem \ref{thm:thicksmth}).  As a consequence, we obtain that every \textbf{box} (a basis set that has the shape of a rectangular parallelepiped) is smooth (Proposition \ref{prop:boxsmth}).
    \item[Box addition]
          In Section \ref{S:boxes}, we show that a basis set $\beta$ that is obtained by adding a box to a smooth basis set in an appropriate way (see Subsection \ref{SS:boxadddefn}) is smooth (Theorem \ref{thm:boxaddnmain}).  Therefore, if we start with a box, and add a finite sequence of boxes to obtain a \textbf{compound box} (see Figure \ref{fig:cpbox}), we obtain a smooth basis set (Corollary \ref{cor:cpboxsmth}).  In particular, \emph{every} basis set in two variables is a compound box (see Figure \ref{fig:2varfig}) and therefore smooth, so we have obtained a variant proof of Haiman's result \cite[Prop.\ 2.4, p.\ 209]{Haiman:CN-HS}, using a lemma valid in arbitrary dimension.  In addition, we prove in Section \ref{S:threevars} that every smooth basis set in three variables is a compound box (Theorem \ref{thm:3varmain}).
    \item[Union of two boxes]
          In Section \ref{S:twoboxes}, we show that a basis set that is the union of two boxes  is a smooth basis set (Theorem \ref{thm:U2boxsmth}).  We present this result to show that the characterization of smooth basis sets in three variables (as compound boxes) does not carry over to higher dimensions, because for $r$ $\geq$ $4$ there are two-box unions that are not compound boxes (Example \ref{SS:2boxeg}).
\end{description} 

\medskip
Sections \ref{S:trunc} and \ref{S:truncsmth} are devoted to technical preparations needed for Sections \ref{S:boxes} and \ref{S:threevars}.  Section \ref{S:trunc} introduces the operation of \emph{truncation} of a basis set (see Figure \ref{Fig:trunc}), and Section \ref{S:truncsmth} establishes conditions under which a truncation of a smooth basis set will be smooth.  

\subsection{Table of contents} \label{SS:toc}

\begin{enumerate}
  \item[1.] Introduction.  
  \item[] \begin{enumerate}
           \item[1.1] Summary of results.
           \item[1.2] Table of contents.
        \end{enumerate}
  \medskip
  \item[2.] The Hilbert scheme of points $\Hn^n$ and its open subschemes $U_{\beta}$.  
  \item[] \begin{enumerate}
           \item[2.1] The Hilbert scheme $\Hn^n$.
           \item[2.2] The open subschemes $U_{\beta}$.
           \item[2.3] The coordinate ring of $U_{\beta}$.
           \item[2.4] The locus of reduced subschemes meets $U_{\beta}$.
        \end{enumerate}
  \medskip
  \item[3.] The cotangent space of the point $t_{\beta}$ $\in$ $\Hn^n$.  
  \item[] \begin{enumerate}
           \item[3.1] Arrows, translation, and congruence modulo $M^2$.
           \item[3.2] $\gf$-{L}inearly independent sets of arrows in $M/M^2$.
        \end{enumerate}
  \medskip
  \item[4.] A linearly independent set in $M/M^2$ of cardinality $r n$.  
  \item[] \begin{enumerate}
           \item[4.1] \Sa\ and minimal arrows.
           \item[4.2] Advancement of minimal standard arrows.
           \item[4.3] Shadow promotion.
           \item[4.4] Iterated shadow promotion.
           \item[4.5] Proof of Theorem \ref{thm:S}
           \item[4.6] Summary and terminology.
        \end{enumerate}
  \medskip
  \item[5.] Consequences of Theorem \ref{thm:S}.  
  \item[] \begin{enumerate}
           \item[5.1] Necessary and sufficient conditions for $t_{\beta}$ $\in$ $\Hn^n$ to be nonsingular.
           \item[5.2] A sufficient condition for $t_{\beta}$ $\in$ $\Hn^n$ to be singular.
           \item[5.3] Example: $\beta$ $=$ $\{1, x_1, x_2, x_3 \}$.
        \end{enumerate}
  \medskip
  \item[6.] Thickening of basis sets.  
  \item[] \begin{enumerate}
           \item[6.1] Definition of thickening.
           \item[6.2] Minimal generators of $I_{\beta}$.
           \item[6.3] A standard bunch $\lis$ for $\beta$.
           \item[6.4] Thickenings of smooth basis sets are smooth.
           \item[6.5] Example: ``Boxes'' are smooth basis sets.
        \end{enumerate}
  \medskip
  \item[7.] Truncation of basis sets.  
  \item[] \begin{enumerate}
           \item[7.1] Definition of truncation.
           \item[7.2] Minimal generators of $I_{\beta_t}$.
           \item[7.3] Lifting arrows from $\beta_t$ to $\beta$.
           \item[7.4] Non-standard arrows on $\beta$ and $\beta_t$.
        \end{enumerate}
  \medskip
  \item[8.] Sufficient conditions for a truncation to be smooth.  
  \item[] \begin{enumerate}
           \item[8.1] The additional hypothesis.
           \item[8.2] $x_j$-sub-bunches of arrows for $\beta$ and $\beta_t$.
           \item[8.3] $x_k$-standard arrows of $x_j$-height $\geq$ $h$.
           \item[8.4] Linear independence of lifts of $x_k$-sub-bunches.
           \item[8.5] $x_k$-sub-bunches of arrows for $\beta$ and $\beta_t$, $x_k$ $\neq$ $x_j$.
           \item[8.6] Main theorem on truncations.
        \end{enumerate}
  \medskip
  \item[9.] Addition of boxes to basis sets.  
  \item[] \begin{enumerate}
           \item[9.1] Definition of box addition.
           \item[9.2] Minimal generators of $I_{\beta^{\prime}}$.
           \item[9.3] Main theorem on box additions.
           \item[9.4] Compound boxes.
           \item[9.5] Example: Basis sets in two variables.
           \item[9.6] Example: The lexicographic point.
           \item[9.7] Example: $\beta$ $=$ $\{1,\, x_1,\, x_2,\, x_1 x_2,\, x_3 \}$.
        \end{enumerate}
  \medskip
  \item[10.] Smooth basis sets in three variables are compound boxes.  
  \item[] \begin{enumerate}
           \item[10.1] The main lemma.
           \item[10.2] Proof of Lemma \ref{lem:3varmainlem}.
           \item[10.3] The main theorem.
        \end{enumerate}
  \medskip
  \item[11.] The union of two boxes.  
  \item[] \begin{enumerate}
           \item[11.1] Notation.
           \item[11.2] Minimal generators of $I_{\beta}$.
           \item[11.3] Two-box unions are smooth basis sets.
           \item[11.4] Example: $\beta$ $=$ $\{1,\, x_1,\, x_2,\, x_1 x_2,\, x_3,\, x_4,\, x_3 x_4 \}$.
        \end{enumerate}
\end{enumerate}


\section{The Hilbert scheme of points $\Hn^n$ and its open subschemes $U_{\beta}$} \label{S:backgrnd}

     In this section we briefly recall the definition and some properties of the Hilbert scheme of points $\Hn^n$, and the open subschemes $U_{\beta}$ that form an open covering of $\Hn^n$.  Much of our subsequent work is based on the explicit representation of the coordinate ring of $U_{\beta}$ that is obtained in \cite{Huib:UConstr}.

\subsection{The Hilbert scheme $\Hn^n$} \label{SS:HnDef}

     The scheme $\Hn^n$ can be defined functorially as the parameter scheme 
of a universal flat and proper family $Z_n$ of zero-dimensional closed subschemes of 
$\mathbb{A}^r_{\gf}$ having length $n$. In other words, $Z_n$ $\subseteq$ 
$\pr{\Hn^n}{}{\mathbb{A}^r_{\gf}}$ is a closed subscheme that is finite and flat of 
degree $n$ over $\Hn^n$ and satisfies the following universal property:
\begin{equation} \label{E:uprop}
      \parbox{4in}{\emph{Let $T$ be a scheme over 
$\gf$.   Then the set of maps $\fn{f}{T}{\Hn^n}$ is in natural bijective 
correspondence with the set of closed subschemes $Z_f \subseteq 
\pr{T}{}{\mathbb{A}^r_{\gf}}$ that are finite and flat of degree $n$ over $T$; the bijection $f 
\mapsto Z_f$ is defined by $Z_f$ = $\pr{T}{\Hn^n}{Z_n}$.}}
\end{equation}
In particular, the inclusion of the $\gf$-point $t \in \Hn^n$ corresponds   
to a unique closed subscheme 
\[
  Z_t \subseteq \pr{t}{}{\mathbb{A}^r_{\gf}} \approx \mathbb{A}^r_{\gf};
\]
the map $t \mapsto Z_t$ defines a bijection 
from the set of $\gf$-points of $\Hn^n$ to the set of 0-dimension\-al closed 
subschemes of length $n$, or, equivalently, to the set of ideals 
\[
    I \subseteq \gf[x_1, \dots, x_r] = \gf[\xlist]
\]
of colength $n$.   We write $I_t$ to denote the ideal corresponding to the $\gf$-point $t$ $\in$ $\Hn^n$.

The existence of $\Hn^n$ can be established as follows: It is the open 
subscheme of $\Hilb^n_{\mathbb{P}^r_{\gf}}$ (a projective scheme the existence of 
which is a consequence of Grothen\-dieck's general construction given in 
\cite{GrothHilbSch}) arising from the inclusion of $\mathbb{A}^r_{\gf}$ in 
$\mathbb{P}^r_{\gf}$ as a standard affine.  Alternatively, $\Hn^n$ is a special case of the multigraded Hilbert scheme constructed by Haiman and Sturmfels \cite{Haiman-Sturmfels:MGHS}.  Elementary constructions of $\Hn^n$ are given in \cite{LakSkjGus:ElHilbSchConst} and \cite{Huib:UConstr}.

\subsection{The open subschemes $U_{\beta}$} \label{SS:Udef}

Let $\beta$ denote a nonempty set of $n$ monomials in the indeterminates 
$x_1, \dots, x_r$ that satisfies the following property: if $m_1$ $\in$ $\beta$ 
and $m_2$ is a monomial dividing $m_1$, then $m_2$ $\in$ $\beta$.   We will call 
such a set $\beta$ a \textbf{basis set}, and its members \textbf{basis 
monomials}.   As a set of $\gf$-points, we define $U_{\beta}$ as follows:
\begin{equation} \label{E:Ubeta}
   U_{\beta} = \{ t \in \Hn^n   \mid   \gf[\xlist]/I_t \text{ has $\gf$-basis 
} \beta \} \subseteq \Hn^n.
\end{equation}
Every point $t$ $\in$ $\Hn^n$ belongs to at least one set of the form 
$U_{\beta}$, a fact that is, in M.\ Haiman's words, ``often regarded as a part of 
Gr\"obner basis theory but actually goes back to Gordan \cite{Gordan1}'' \cite[p.\ 
207]{Haiman:CN-HS}.  Therefore, the sets $U_{\beta}$ form an open covering of $\Hn^n$.  In fact, the $U_{\beta}$ form an open \emph{affine} covering of $\Hn^n$; this is proved in \cite{Huib:UConstr}, where the coordinate ring of $U_{\beta}$ is explicitly obtained as a quotient of a polynomial ring.   Note that Haiman introduced the open subschemes $U_{\beta}$ in \cite{Haiman:CN-HS} for the case of the affine plane $(r = 2)$; he showed that they are affine using other means.

\begin{rem} \label{rem:free}
     Let $Z_{\beta}$ denote the restriction of the universal closed subscheme $Z_n$ to $U_{\beta}$ $\subseteq$ $\Hn^n$.  Then the direct image of the structure sheaf of $Z_{\beta}$ on $U_{\beta}$ is free with basis $\beta$, because $\beta$ is everywhere a local basis of the (locally free) direct image.  This also follows directly from the construction of $(U_{\beta}, Z_{\beta})$ given in \cite[Sec.\ 7.3]{Huib:UConstr}.
\end{rem}

Note that the monomials that do \emph{not} belong to a basis set $\beta$ generate a monomial ideal 
\[
    I\ =\ I_{\beta}\ \subseteq\ \gf[\xlist];
\]
it is immediate that the quotient $\gf[\xlist]/I_{\beta}$ has $\beta$ as $\gf$-basis, that is, $I_{\beta}$ corresponds to a point 
\begin{equation} \label{E:tbeta}
 t\ =\ t_{\beta}\ \in\ U_{\beta}.
\end{equation}
Conversely, if $I$ is a monomial ideal such that $\gf[\xlist]/I$ has $\gf$-dimension $n$, then the quotient has as $\gf$-basis the set of monomials not in $I$, and this is clearly a basis set $\beta$ of $n$ members.  

\subsection{The coordinate ring of $U_{\beta}$} \label{SS:Uring}

Let $R$ denote the coordinate ring of the affine scheme $U_{\beta}$, and let $I_{U_{\beta}}$ $\subseteq$ $R[\xlist]$ be the ideal that cuts out the universal closed subscheme $Z_{\beta}$ $\subseteq$ $\pr{U_{\beta}}{}{\mathbb{A}^r_{\gf}}$.  By Remark \ref{rem:free}, we have that the module $R[\xlist]/I_{U_{\beta}}$ is $R$-free with basis $\beta$.  Therefore, every monomial 
\[
     x_1^{d_1}x_2^{d_2}\dots x_r^{d_r} = \xlist^{\dlist}\in R[\xlist]
\]
is congruent modulo $I_{U_{\beta}}$ to a unique $R$-linear combination of the monomials $\xlist^{\jlist}$ $\in$ $\beta$, or, in other words, for every $\xlist^{\dlist}$ there is a unique polynomial of the form
\begin{equation} \label{E:Fmdef}
       F_{\dlist}\ =\ \xlist^{\dlist} - \sum_{\jlist \in \beta} 
c^{\dlist}_{\jlist} \cdot \xlist^{\jlist}\  \in\  I_{U_{\beta}},\ \ c^{\dlist}_{\jlist} \in R,
\end{equation}
where we abuse the notation (here and elsewhere) by writing $\jlist \in \beta$ for $\xlist^{\jlist} \in \beta$.  When writing monomials in this way, we will reserve $\jlist$ for basis monomials.  Note that \begin{equation} \label{E:speccs}
    \dlist \in \beta\ \Rightarrow 
        \left\{
        \begin{array}{l}
            c^{\dlist}_{\jlist} = 1 \text{ if } \dlist = \jlist, \text{ and}\vspace{.05in}\\
            c^{\dlist}_{\jlist} = 0 \text{ if } \dlist \neq \jlist.
        \end{array} \right. 
\end{equation}

Following Haiman \cite[p.\ 210]{Haiman:CN-HS}, we multiply the polynomial (\ref{E:Fmdef}) by the variable $x_i$, and then expand each monomial $x_i \cdot \xlist^{\jlist}$ using (\ref{E:Fmdef}) to obtain a polynomial of the form 
\[
    x_i \cdot \xlist^{\dlist} + (\text{$R$-linear combination of basis monomials})\ \in\ I_{U_{\beta}}, 
\]
which must therefore be equal to $F_{\dlist^{\prime}}$, where $\xlist^{\dlist^{\prime}}$ = $x_i \cdot \xlist^{\dlist}$.  Equating coefficients, we obtain the relations
 \begin{equation} \label{E:crelns}
     c^{\dlist^{\prime}}_{\jlist_0} -  \sum_{\jlist \in \beta} c^{\dlist}_{\jlist} \cdot c^{\jlist^{\prime}}_{\jlist_0}\ = \ 0,    
\end{equation}
where $\xlist^{\jlist^{\prime}}$ $=$ $x_i \cdot \xlist^{\jlist}$.  

For each coefficient $c^{\dlist}_{\jlist}$ such that $\dlist$ $\notin$ $\beta$ we introduce an indeterminate $C^{\dlist}_{\jlist}$, and let 
\begin{equation} \label{E:deldef}
    \fn{\delta}{\gf[( C^{\dlist}_{\jlist} )]}{R},\  C^{\dlist}_{\jlist} \mapsto c^{\dlist}_{\jlist},
\end{equation}
be the natural map.  We have the following
\begin{prop} \label{prop:genprop}
    The coordinate ring $R$ of $U_{\beta}$ is generated as a $\gf$-algebra by the coefficients $c^{\dlist}_{\jlist}$ such that $\dlist$ $\notin$ $\beta$; that is, the map $\delta$ is surjective.  Furthermore, the kernel of $\delta$ is generated by the polynomials 
\[\tau^{(\dlist,x_i)}_{\jlist_0} =  C^{\dlist^{\prime}}_{\jlist_0} -  \sum_{\jlist \in \beta} C^{\dlist}_{\jlist} \cdot C^{\jlist^{\prime}}_{\jlist_0},\ \ \dlist \notin \beta,\ \jlist_0 \in \beta,\ 1 \leq i \leq r,
\]
that are obtained from {\rm (\ref{E:crelns})} by replacing each coefficient $c^{\blist}_{\jlist}$ with the corresponding indeterminate $C^{\blist}_{\jlist}$, if $\blist$ $\notin$ $\beta$, or the appropriate constant --- $0$ or $1$, according to \emph{(\ref{E:speccs})}, and denoted $C^{\blist}_{\jlist}$ --- if $\blist$ $\in$ $\beta$.
\end{prop}

\emph{Proof:}
     The construction of $R$ given in \cite[Sec.\ 7]{Huib:UConstr} shows that $R$ is generated by a certain finite subset of the coefficients $c^{\dlist}_{\jlist}$, and the kernel of the map $\delta$ is generated by a finite set of polynomials (there denoted $\rho^{(\blist_1, \blist_2)}_{\jlist}$), each of which is in fact either equal to one of the polynomials $\tau$ or to a $\gf$-linear combination of two of them.   
\qedsymbol  

\begin{rem} \label{rem:tcoords}
    It is clear that the point $t_{\beta}$  (\ref{E:tbeta}) is the ``origin'' of $U_{\beta}$; that is, $t_{\beta}$ is the point at which all the coefficient functions $c^{\dlist}_{\jlist}$ with $\dlist$ $\notin$ $\beta$ vanish.
\end{rem}

\subsection{The locus of reduced subschemes meets $U_{\beta}$} \label{SS:lrsub}

Let $\Hn_{\circ}$ $\subseteq$ $\Hn^n$ denote the open subscheme parameterizing the closed subschemes of $\mathbb{A}^r_{\gf}$ that are supported at $n$ distinct points (and are therefore reduced).  Since $\Hn_{\circ}$ is irreducible and of dimension $r n$, its closure $\bar{\Hn}_{\circ}$ is a component of $\Hn^n$ of dimension $r n$.  It is well-known that the point $t_{\beta}$   corresponding to the monomial ideal $I_{\beta}$ $\subseteq$ $\gf[\xlist]$ lies in $\bar{\Hn}_{\circ}$; in fact, one can exhibit a one-parameter family of length-$n$ closed subschemes of $\mathbb{A}^r_{\gf}$ such that the fiber over 0 is the subscheme defined by $I_{\beta}$ and the general fiber is a reduced subscheme.  This family is constructed using Hartshorne's concept of \textbf{distraction} \cite{Hart:ConnHS} as adapted by Geramita \etal\ in \cite{Geramita}.  The construction works whenever the ground field is infinite (which is true for us, since $\gf$ is assumed algebraically closed) or sufficiently large.

Briefly, the construction goes like this.  First, choose an infinite (or sufficiently long) sequence of distinct elements $a_0, a_1, a_2, \dots$ of $\gf$.  Then, for every minimal generator 
\[
    x_1^{d_1} x_2^{d_2} \dots x_r^{d_r} = \xlist^{\dlist}
\]
of the monomial ideal $I_{\beta}$, form the homogeneous polynomial 
\[
    f_{\dlist}\ =\ \prod_{i=1}^r \left( \prod_{j=0}^{d_i - 1} (x_i - a_j \cdot w) \right)\ \in\ \gf[w][\xlist]. 
\]
One then proves that the homogeneous ideal $\mathcal{F}$ generated by the polynomials $f_{\dlist}$ is the ideal of all functions vanishing on the $n$ (distinct projective) points
\[
    [1 : a_{e_1} : a_{e_2} : \dots : a_{e_r} ] \in \mathbb{P}^r_{\gf},\ \ x_1^{e_1} x_2^{e_2} \cdots x_r^{e_r} = \xlist^{\elist} \in \beta
\]
(\cite[Theorem 2.2]{Geramita}).  In particular, if we view $w$ as a parameter, then for any nonzero value $w_1$ of $w$, the ideal $\mathcal{F}(w_1)$ $\subset$ $\gf[\xlist]$ obtained by replacing the variable $w$ by the constant $w_1$ is supported at the $n$ distinct points 
\[
    (w_1 \cdot a_{e_1}, w_1 \cdot a_{e_2}, \dots, w_1 \cdot a_{e_r}),\ \xlist^{\elist} \in \beta
\]
(when $w_1$ $=$ $0$, clearly $\mathcal{F}(0)$ $=$ $I_{\beta}$, the original monomial ideal of colength $n$).
One then checks that the quotient $\gf[w][\xlist]/\mathcal{F}$ is $\gf[w]$-free with basis $\beta$ as a $\gf[w]$-module.  (Indeed, from the form of the polynomials $f_{\dlist}$, it is clear that every minimal generator of the monomial ideal $I_{\beta}$ is congruent (mod ${\mathcal{F}}$) to a $k[w]$-linear combination of basis monomials \emph{that divide the minimal generator}.  Induction on the degree then shows that every monomial $m$ $\in$ $I_{\beta}$ is congruent (mod ${\mathcal{F}}$) to a $k[w]$-linear combination of basis monomials that divide $m$; therefore, the basis monomials span the quotient $\gf[w][\xlist]/\mathcal{F}$.  Since the $\gf$-dimension of each quotient $\gf[\xlist]/\mathcal{F}(w_1)$ is $\geq$ $n$, we find that $\beta$ is locally a basis everywhere, and hence a global basis, as desired.) In this way we obtain a finite and flat family of subschemes which defines a map $\fn{\varphi}{\Spec(\gf[t])}{U_{\beta}$ $\subseteq$ $\Hn^n}$; one sees easily that $\varphi(0)$ $=$ $t_{\beta}$ and $\varphi(w_1)$ $\in$ $\Hn_{\circ}$ $\cap$ $U_{\beta}$ for all $w_1$ $\neq$ $0$.  It follows at once that $t_{\beta}$ $\in$ $\bar{\Hn}_{\circ}$, which yields

\begin{prop} \label{prop:freedim}

    The $\gf$-dimension of the (Zariski) tangent and cotangent spaces of the point $t_{\beta}$ is $\geq$ $\dim(\bar{\Hn}_{\circ})$ $=$ $r n$, and $t_{\beta}$ is a smooth point of $\Hn^n$ if and only if 
\[
    \dim_{\gf}(\text{\rm cotangent space of } t_{\beta}) = r n. \ \ \text{\qedsymbol}
\]

\end{prop}

\begin{rem} \label{rem:anotherpf}
    We will give another proof that the $\gf$-dimension of the cotangent space of $t_{\beta}$ is at least $r n$ in Section \ref{S:lindset}.
\end{rem}

\begin{rem} \label{rem:redloc}
    The foregoing demonstrates that each open subscheme $U_{\beta}$ meets the locus $\Hn_{\circ}$ of reduced subschemes nontrivially.  This is shown in a different way in \cite[Sec.\ 2.4]{Huib:UConstr}.  
\end{rem}


\section{The cotangent space of the point $t_{\beta}$ $\in$ $\Hn^n$} \label{sec:cotansp}

     Recall that the cotangent space of the point $t_{\beta}$ $\in$ $\Hn^n$ is the $\gf$-vector space $\ideal{m}_{t_{\beta}}/\ideal{m}^2_{t_{\beta}}$, where $\ideal{m}_{t_{\beta}}$ is the maximal ideal of the local ring $\sheaf{O}_{t_{\beta}}$.  In light of Remark \ref{rem:tcoords}, we see that we can identify the cotangent space with the $\gf$-vector space $M/M^2$, where $M$ is the maximal ideal generated by the functions $c^{\dlist}_{\jlist}$ ($\dlist$ $\notin$ $\beta$) in the coordinate ring $R$ of $U_{\beta}$.  In \cite{Haiman:CN-HS}, Haiman observed, in the two-variable case, that congruence (mod ${M^2}$) can be visualized as translation-equivalence of arrows defined by the $c^{\dlist}_{\jlist}$.  We now extend Haiman's idea to any number $r$ of variables.

\subsection{Arrows, translation, and congruence modulo $M^2$} \label{SS:arrtrans}

     As stated by Haiman, with symbols adjusted, we have that \cite[p.\ 210]{Haiman:CN-HS} 

\begin{quotation} 
    [m]odulo $M^2$, the terms 
-$c^{\dlist}_{\jlist} \cdot c^{\jlist^{\prime}}_{\jlist_0}$ [in] equation (\ref{E:crelns}) reduce to zero for $\jlist^{\prime}$ $\notin$ $\beta$ and for $\jlist^{\prime}$ $\in$ $\beta$, $\jlist^{\prime}$ $\neq$ $\jlist_0$ [recall (\ref{E:speccs}) --- here we are assuming that $\dlist$ $\notin$ $\beta$, so that $c^{\dlist}_{\jlist}$ $\in$ $M$].  The remaining term is -$c^{\dlist}_{\jlist_1}$, where 
\[
    \jlist_1^{\prime} = \jlist_0, \text{ that is, } \xlist^{\jlist_1} \cdot x_i = \xlist^{\jlist_0},
\]
or zero if $\xlist^{\jlist_0}$ is not divisible by $x_i$.  Thus in $M/M^2$ we have
\begin{equation} \label{E:transeqn}
    c^{\dlist^{\prime}}_{\jlist_0}\ =\
    \left\{
      \begin{array}{l}
          c^{\dlist}_{\jlist_1}, \text { if } x_i \text{ divides } \xlist^{\jlist_0}, \text { and}\\
          0, \text{ otherwise}.
      \end{array}
    \right.
\end{equation}

It is convenient to depict each $c^{\dlist}_{\jlist}$ by an arrow from $\dlist$ to $\jlist$ (see Figure \ref{fig:arroweg}).  Equation (\ref{E:transeqn}) says that we may move these arrows [in the $x_i$-direction, $1$ $\leq$ $i$ $\leq$ $r$,] without changing their values modulo $M^2$, provided we keep the head inside $\beta$ and the tail outside.  More generally, as long as we keep the tail in the first [orthant] and outside $\beta$, we may even move the head across the [hyperplane ($x_i$-degree $=$ $0$)].  When this is possible, the value of the arrow (mod ${M^2}$) is zero.
\end{quotation}
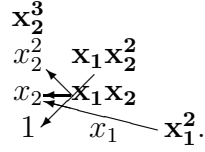
\begin{figure} 
\begin{picture}(374,72)
        \put(130,32){
               \mbox{$\begin{array}{ccc}
                \mathbf{x_2^3} & {} & {} \\
                x_2^2 & \mathbf{x_1 x_2^2} & {}     \\
                x_2   & \mathbf{x_1 x_2} & {}   \\
                1 & x_1 & \mathbf{x_1^2\text{.}} 
                     \end{array}
                    $}
                   }
        \put(170,35){\vector(-1,-1){20}}
        \put(161,27){\vector(-1,0){10}}
        \put(162,27){\vector(-1,1){10}}
        \put(194,14){\vector(-4,1){43}}
\end{picture}
\caption{The ``arrows'' $c^{(1,2)}_{(0,0)}$, $c^{(1,1)}_{(0,1)}$, $c^{(1,1)}_{(0,2)}$, and $c^{(2,0)}_{(0,1)}$.  Monomials not in $\beta$ are shown in boldface.} \label{fig:arroweg}
\end{figure}

Henceforth when we speak of an ``arrow'' $c^{\dlist}_{\jlist}$, we will assume that $\dlist$ $\notin$ $\beta$, so that there is a corresponding indeterminate $C^{\dlist}_{\jlist}$ (as defined in Section \ref{SS:Uring}).  We will say that two arrows $c^{\blist_1}_{\jlist_1}$ and $c^{\blist_2}_{\jlist_2}$ are \textbf{translation-equivalent} if the first can be moved to the second by a sequence of discrete steps in the various variable directions such that the head of the arrow remains inside $\beta$ and the tail remains outside $\beta$.  This clearly defines an equivalence relation on the set of arrows, which we denote $c^{\blist_1}_{\jlist_1}$ $\sim$ $c^{\blist_2}_{\jlist_2}$.  For example, in Figure \ref{fig:arroweg}, we see that $c^{(1,1)}_{(0,1)} \sim c^{(1,2)}_{(0,2)}$.  Abusing the language and notation, we say $c^{\blist}_{\jlist}$ is \textbf{translation-equivalent to} 0, and write $c^{\blist}_{\jlist} \sim 0$, to indicate that the arrow can be translated to a position such that one more step in some direction of decreasing degree would cause the head of the arrow to exit the first orthant, with the tail remaining a monomial outside of $\beta$.  The arrow $c^{(1,2)}_{(0,0)}$ in Figure \ref{fig:arroweg} provides an example: one further step in the decreasing $x_2$-direction would cause the head of this arrow to exit the first quadrant.

From the foregoing, it is clear that
\begin{equation} \label{E:simcong}
   \begin{array}{rcl}
        c^{\blist_1}_{\jlist_1} \sim c^{\blist_2}_{\jlist_2} & \Rightarrow & c^{\blist_1}_{\jlist_1} \equiv c^{\blist_2}_{\jlist_2} \pmod{M^2}, \text { and}\\
        c^{\blist}_{\jlist} \sim 0 & \Rightarrow & c^{\blist}_{\jlist} \equiv 0 \pmod{M^2}.
   \end{array}
\end{equation}
Furthermore, the reasoning in the quoted passage can be adapted to prove

\begin{prop} \label{prop:taufacts}
Let $\tau^{(\dlist,x_i)}_{\jlist_0}$ be any one of the polynomial generators of the kernel of the map $\delta$ \emph{(\ref{E:deldef})}.  Then 
    \begin{enumerate}
        \item[(a)] Each  term in $\tau^{(\dlist,x_i)}_{\jlist_0}$ is, up to sign, either a single indeterminate $C^{\blist}_{\jlist}$ or a product of two such indeterminates.
        \item[(b)] The number of linear terms in $\tau^{(\dlist,x_i)}_{\jlist_0}$ is equal to $1$ or $2$.
        \item[(c)] If there are two linear terms in $\tau^{(\dlist,x_i)}_{\jlist_0}$, then these terms have the form $C^{\dlist^{\prime}}_{\jlist_0}$ and -$C^{\dlist}_{\jlist_1}$, where $\xlist^{\dlist^{\prime}}$ $=$ $x_i \cdot \xlist^{\dlist}$ and $\xlist^{\jlist_0}$ $=$ $x_i \cdot \xlist^{\jlist_1}$; in particular, the signs differ, and the corresponding arrows $c^{\dlist^{\prime}}_{\jlist_0}$ and $c^{\dlist}_{\jlist_1}$ are translation-equivalent (\ie, $c^{\dlist^{\prime}}_{\jlist_0}$ $\sim$ $c^{\dlist}_{\jlist_1}$). 
        \item[(d)] If there is only one linear term  in $\tau^{(\dlist,x_i)}_{\jlist_0}$, then it is the term $C^{\dlist^{\prime}}_{\jlist_0}$, and the corresponding arrow $c^{\dlist^{\prime}}_{\jlist_0}$ $\sim$ $0$.\ \  \qedsymbol
    \end{enumerate}
\end{prop}

\subsection{$\gf$-{L}inearly independent sets of arrows in $M/M^2$} \label{SS:ldsets} 

Our main goal in this section is to prove the following 

\begin{thm} \label{thm:lindepconds}

     Let $S$ be a finite set of arrows $c^{\dlist}_{\jlist}$ having $d$ members.  Then $S$ has maximal rank $(\text{mod } M^2)$ (that is, the $\gf$-span of $S$ in $M/M^2$ has dimension $d$) if and only if the following conditions hold:   
\begin{enumerate}
    \item[(a)] $c^{\dlist}_{\jlist}$ $\nsim$ $0$ for all $c^{\dlist}_{\jlist}$ $\in$ $S$.
    \item[(b)] $c^{\dlist_{1}}_{\jlist_{1}}$ $\nsim$ $c^{\dlist_{2}}_{\jlist_{2}}$ for all $c^{\dlist_{1}}_{\jlist_{1}}$ $\neq$ $c^{\dlist_{2}}_{\jlist_{2}}$ in $S$.
\end{enumerate} 
\end{thm}

Before proceeding with the proof, we make a few preparations.
We first extend the notion of translation-equivalence to the set of indeterminates $C^{\blist}_{\jlist}$ in the obvious way:
\[
  \begin{array}{rcl}
    C^{\blist_1}_{\jlist_1} \sim C^{\blist_2}_{\jlist_2} & \Leftrightarrow &  c^{\blist_1}_{\jlist_1} \sim c^{\blist_2}_{\jlist_2}\vspace{.05in}\\
    C^{\blist}_{\jlist} \sim 0 & \Leftrightarrow & c^{\blist}_{\jlist} \sim 0.
  \end{array}
\]
Then Proposition \ref{prop:taufacts} immediately yields
\begin{lem} \label{lem:taufacts}
    Let $\tau^{(\dlist,x_i)}_{\jlist_0}$ be any one of the polynomial generators of the kernel of the map $\delta$ \emph{(\ref{E:deldef})}.  Then 
    \begin{enumerate}
        \item[(a)] If $\tau^{(\dlist,x_i)}_{\jlist_0}$ has two linear terms $C^{\dlist^{\prime}}_{\jlist_0}$ and -$C^{\dlist}_{\jlist_1}$, then $C^{\dlist^{\prime}}_{\jlist_0}$ $\sim$ $C^{\dlist}_{\jlist_1}$.
        \item[(b)] If $\tau^{(\dlist,x_i)}_{\jlist_0}$ has just one linear term $C^{\dlist^{\prime}}_{\jlist_0}$, then $C^{\dlist^{\prime}}_{\jlist_0}$ $\sim$ $0$.\ \  \qedsymbol
    \end{enumerate}
 
\end{lem}

Given a polynomial $P$ $\in$ $\gf[(C^{\blist}_{\jlist})]$ of degree $q$, we write
\[
    P\ = \ P^{(0)} + P^{(1)} + \dots + P^{(q)},
\]  
where each $P^{(j)}$ is homogeneous of degree $j$.  We also define the map
\begin{equation} \label{E:Edefn}
    \fn{E}{\gf[(C^{\blist}_{\jlist})]}{\gf[(\langle C^{\blist}_{\jlist}\rangle)]},\ C^{\blist}_{\jlist} \mapsto \langle C^{\blist}_{\jlist}\rangle,
\end{equation}
where the target is the ring of polynomials in the translation-equivalence classes $\langle C^{\blist}_{\jlist} \rangle$ of the indeterminates $C^{\blist}_{\jlist}$.  

\begin{lem} \label{lem:sublem1}
    Let $a_1, \dots, a_s$ $\in$ $\gf$, and $\tau_1, \dots, \tau_s$ be $s$ of the polynomials $\tau^{(\dlist,x_i)}_{\jlist_0}$.  Let $N$ denote the set of indices $n$ such that $\tau_n$ has just one linear term $C^{\blist_n}_{\jlist_n}$.  Then
\[
    E( \sum_{j=1}^s a_j \cdot \tau^{(1)}_j )\ = \ \sum_{n \in N} a_n \cdot \langle C^{\blist_n}_{\jlist_n}\rangle.
\]
\end{lem}

\emph{Proof:}
    It suffices to observe that for the $\tau_j$ containing two linear terms $C^{\blist^{\prime}_j}_{\jlist_0}$ and $-C^{\blist_j}_{\jlist_1}$, we have that
\[
    E (a_j \cdot \tau^{(1)}_j) \ =\ a_j \cdot \left( \langle C^{\blist^{\prime}_j}_{\jlist_0}\rangle - \langle C^{\blist_j}_{\jlist_1}\rangle \right) \ =\ 0,
\]
where the last equality follows from Lemma \ref{lem:taufacts}.  \qedsymbol

\bigskip

\noindent
\emph{Proof of Theorem} \ref{thm:lindepconds}:
From the implications (\ref{E:simcong}), it is clear that if either of the conditions (a), (b) in the statement of the theorem fails, then the dimension of the $\gf$-span of $S$ in $M/M^2$ is $<$ $d$; therefore, if the dimension is $d$, then (a) and (b) must hold.  To prove the converse, we argue by contradiction: Suppose that (a) and (b) hold, but that the $\gf$-span of $S$ in $M/M^2$ has dimension $<$ $d$.  Then there exists a nontrivial $\gf$-linear combination of (distinct) elements of $S$ that is congruent to 0 modulo $M^2$, say
\[
    \sum_{i = 1}^m \alpha_i \cdot c^{\dlist_i}_{\jlist_i}\ \equiv\ 0 \pmod{M^2}, 
\]
which implies that 
\[
    \sum_{i = 1}^m \alpha_i \cdot c^{\dlist_i}_{\jlist_i} - (\text{elt.\ of } M^2)\ =\ 0\ \in\ R\ \ \ (R = \text{coord.\ ring of } U_{\beta}).
\]
From the description of $R$ as a quotient of the polynomial ring $\gf[(C^{\dlist}_{\jlist})]$ provided by Proposition \ref{prop:genprop}, we see that we have an equation
\[
    \sum_{i = 1}^m \alpha_i \cdot C^{\dlist_i}_{\jlist_i} - (\text{terms  in $C^{\dlist}_{\jlist}$ of degree } \geq 2)\ =\ \sum_{j = 1}^s g_j \cdot \tau_j,
\]
where the coefficients $g_j$ $\in$ $\gf[(C^{\dlist}_{\jlist})]$ and each $\tau_j$ is one of the polynomials $\tau^{(\dlist,x_i)}_{\jlist_0}$ that generate the kernel of the map $\delta$ (\ref{E:deldef}).  Since the $\tau_j$ have only linear and quadratic terms (by Proposition \ref{prop:taufacts}), we see that the terms of the form $g^{(i)}_j \cdot \tau_j$ for $i \geq 1$ can be transposed to the left to yield an equation
\[
    \sum_{i = 1}^m \alpha_i \cdot C^{\dlist_i}_{\jlist_i} - (\text{terms  in $C^{\dlist}_{\jlist}$ of degree } \geq 2)\ =\ \sum_{j = 1}^s g^{(0)}_j \cdot \tau_j.
\]
Equating the degree-1 terms on both sides, we obtain
\[
    \sum_{i = 1}^m \alpha_i \cdot C^{\dlist_i}_{\jlist_i}\ =\ \sum_{j = 1}^s g^{(0)}_j \cdot \tau^{(1)}_j.
\]
As in Lemma \ref{lem:sublem1}, we let $N$ denote the set of indices $n$ such that $\tau^{(1)}_n$ consists of a single term $C^{\blist_n}_{\jlist_n}$.  Applying the map $E$ (\ref{E:Edefn}) to both sides and rewriting the RHS using the lemma, we find that
\[
    \sum_{i = 1}^m \alpha_i \cdot \langle C^{\dlist_i}_{\jlist_i}\rangle \ =\ \sum_{n \in N} g^{(0)}_n \cdot \langle C^{\blist_n}_{\jlist_n}\rangle.
\]
Recall that conditions (a) and (b) hold, by hypothesis, and that we began with a nontrivial linear combination of the $c^{\dlist}_{\jlist}$ $\in$ $S$, so that at least one of the coefficients $\alpha_i$ $\neq$ $0$.  If the corresponding equivalence class $\langle C^{\dlist_i}_{\jlist_i} \rangle$ does not appear on the RHS of the last equation, then the term $\alpha_i \cdot \langle C^{\dlist_i}_{\jlist_i} \rangle$ must cancel with one or more other terms $\alpha_j \cdot \langle C^{\dlist_j}_{\jlist_j} \rangle$ on the LHS; whence, $c^{\dlist_i}_{\jlist_i}$ $\sim$ $c^{\dlist_j}_{\jlist_j}$, which contradicts condition (b).  Therefore, we must have that  $\langle C^{\dlist_i}_{\jlist_i} \rangle$ $=$ $\langle C^{\dlist_n}_{\jlist_n} \rangle$ for some $n$ $\in$ $N$.  That is, we have $C^{\dlist_i}_{\jlist_i}$ $\sim$ $C^{\blist_n}_{\jlist_n}$, but since 
$C^{\blist_n}_{\jlist_n}$ $\sim$ $0$ by Lemma \ref{lem:taufacts}, it follows immediately that $C^{\dlist_i}_{\jlist_i}$  $\sim$ $0$ $\Leftrightarrow$ $c^{\dlist_i}_{\jlist_i}$ $\sim$ $0$, which contradicts condition (a), and the proof is complete.
\qedsymbol

\bigskip

As a corollary, we obtain the following converses to the implications (\ref{E:simcong}):

\begin{cor} \label{cor:lindepconds}

For $c^{\blist}_{\jlist},\ c^{\blist_1}_{\jlist_1},\ c^{\blist_2}_{\jlist_2}$ any arrows, we have that

  \begin{enumerate}
    \item[(a)] $c^{\blist}_{\jlist} \equiv 0 \pmod{M^2} \Rightarrow c^{\blist}_{\jlist} \sim 0$.
    \item[(b)] $c^{\blist_1}_{\jlist_1} \equiv c^{\blist_2}_{\jlist_2} \nequiv 0 \pmod{M^2} \Rightarrow c^{\blist_1}_{\jlist_1} \sim c^{\blist_2}_{\jlist_2}$.
  \end{enumerate}

\end{cor}

\emph{Proof:} 
    Apply the theorem to the sets $\{ c^{\blist}_{\jlist} \}$ and $\{ c^{\blist_1}_{\jlist_1},  c^{\blist_2}_{\jlist_2} \}$, neither of which has maximal rank (mod ${M^2}$), provided that $c^{\blist_1}_{\jlist_1}$ $\neq$ $c^{\blist_2}_{\jlist_2}$.
\qedsymbol


\section{A linearly independent set in $M/M^2$ of cardinality $r n$} \label{S:lindset}

     Given a basis set of monomials $\beta$ of size $n$ in $r$ variables (or equivalently the associated monomial ideal $I_{\beta}$), we exhibit a set $\lis$ of arrows $c^{\dlist}_{\jlist}$ of cardinality $r n$ whose $\gf$-span in $M/M^2$ has dimension $r n$.  This gives a second proof that the $\gf$-dimension of the cotangent space of the point $t_{\beta}$ $\in$ $\Hn^n$ is at least $r n$, as promised in Remark \ref{rem:anotherpf}.  More importantly, if $t_{\beta}$ $\in$ $\Hn^n$ is nonsingular, then $\lis$ must be a basis of the cotangent space.  Because of the particular form of the arrows in $\lis$, it is often easy to show in particular cases that $\lis$ does \emph{not} span the cotangent space; we conclude that $t_{\beta}$ is a singular point in such cases.  On the other hand, there are several families of basis sets $\beta$ for which we can prove that $\lis$ \emph{is} a basis of the cotangent space.

\subsection{\Sa\ and minimal arrows} \label{SS:stdarrows}

We begin by identifying the type of arrow that will belong to our set $\lis$.  First we define the \textbf{vector} of the arrow $c^{\dlist}_{\jlist}$ to be the tuple $\jlist - \dlist$ (recall that we only speak of an ``arrow'' when $\dlist$ $\notin$ $\beta$ $\Leftrightarrow$ $\xlist^{\dlist}$ $\in$ $I_{\beta}$).  Since $\xlist^{\dlist}$ $|$ $\xlist^{\jlist}$ is then impossible, the following is immediate:

\begin{lem}  \label{lem:vec1}
    For every arrow $c^{\dlist}_{\jlist}$, the vector of the arrow has at least one negative component.\ \qedsymbol
\end{lem}

We say that $c^{\dlist}_{\jlist}$ is a \textbf{\sa\ arrow} provided that the vector of the arrow has exactly one negative component; if this component is the $i$-th, corresponding to the variable $x_i$, then we say that $c^{\dlist}_{\jlist}$ is \textbf{\sa\ for} $x_i$, or $x_i$\textbf{-standard}.  Concretely, this means that the monomial $\xlist^{\dlist}$ at the tail of the arrow has a strictly larger $x_i$-degree than the monomial $\xlist^{\jlist}$ at the head, but the $x_k$-degree of $\xlist^{\dlist}$ is $\leq$ the $x_k$-degree of $\xlist^{\jlist}$ for all other variables $x_k$ $\neq$ $x_i$.  For example, in Figure \ref{fig:arroweg}, the arrows $c^{(1,1)}_{(0,1)}$, $c^{(1,1)}_{(0,2)}$, and $c^{(2,0)}_{(0,1)}$ are \sa\ for $x_1$, and $c^{(1,2)}_{(0,0)}$ is not a \sa\ arrow.

Recall that the monomials in $x_1, \dots, x_r$ are partially ordered by divisibility, and $I_{\beta}$ is generated by the minimal monomials in $I_{\beta}$ under this ordering (the \textbf{minimal generators} of $I_{\beta}$).  We call an arrow $c^{\dlist}_{\jlist}$ whose tail $\xlist^{\dlist}$ is a minimal generator a \textbf{minimal arrow}; if the arrow is \sa\ (for $x_i$), we call it a \textbf{minimal \sa\ arrow} (\textbf{for} $x_i$).  In light of the results of Section \ref{sec:cotansp}, the next lemma shows that in seeking a set of arrows to span the $\gf$-vector space $M/M^2$, it suffices to restrict one's attention to minimal arrows.  

\begin{lem} \label{lem:msa}
    Let $c^{\dlist}_{\jlist}$ be an arbitrary arrow, and $\xlist^{\blist}$ a minimal generator of $I_{\beta}$ such that $\xlist^{\blist}$ $|$ $\xlist^{\dlist}$ (at least one such minimal generator must clearly exist).  Then either $c^{\dlist}_{\jlist}$ $\sim$ $0$ or $c^{\dlist}_{\jlist}$ $\sim$ $c^{\blist}_{\jlist_1}$ for some $\jlist_1$ $\in$ $\beta$.
\end{lem}

\emph{Proof:}
    We can clearly translate $\xlist^{\dlist}$ to $\xlist^{\blist}$ by a series of degree-reducing steps in the various variable directions (each step involves dividing the head and tail of the arrow by one of the $x_k$).  Since dividing a basis monomial (at the head of the arrow) by $x_k$ either keeps us inside $\beta$ or causes us to exit the first orthant across the hyperplane ($x_k$-degree $=$ $0$), the result follows immediately. 
\qedsymbol
\bigskip

In fact, the set $\lis$ that we are out to construct will consist of minimal \emph{\sa} arrows, as our main theorem asserts:

\begin{thm} \label{thm:S}
     Let $\beta$ be a basis set of $n$ monomials in the variables $x_1$, $\dots$, $x_r$, and $I_{\beta}$ $\subseteq$ $\gf[\xlist]$ the associated monomial ideal.  Then there exists a set $\lis$ of  \msa s that has cardinality $r n$ and maximal rank $(\text{mod } M/M^2)$; that is, the $\gf$-span of $S$ in $M/M^2$ has dimension $r n$.
\end{thm}

The proof will be given in Section \ref{SS:Sproof}, after the necessary preparations have been made. 

\subsection{Advancement of minimal standard arrows} \label{SS:msadvance}

We begin with a host of definitions.  Let $\beta$, $I_{\beta}$, \etc, be as above.  Note first of all that, since $I_{\beta}$ has finite colength, there is for each variable $x_i$ a minimal exponent $w_i$ $>$ $0$ such that $x_i^{w_i}$ $\in$ $I_{\beta}$; we will call $w_i$ the $x_i$-\textbf{width} of $\beta$.  It is then clear that 
\begin{equation} \label{eqn:width}
    \text{$x_i$-degree}(\xlist^{\jlist}) < w_i \text{ for all }\xlist^{\jlist} \in \beta,\ 1 \leq i \leq r.
\end{equation} 
For example, the basis set shown in Figure \ref{fig:arroweg} has $x_1$-width $w_1$ $=$ $2$ and $x_2$-width $w_2$ $=$ 3.

Given $\xlist^{\jlist}$ $\in$ $\beta$, we define the $x_i$-\textbf{column} of $\xlist^{\jlist}$ to be the set of monomials $\xlist^{\jlist}$, $\xlist^{\jlist}/x_i$, $\xlist^{\jlist}/x_i^2$, \dots, $\xlist^{\jlist}/x_i^q$, where $q$ $=$ $x_i$-degree$(\xlist^{\jlist})$.  Clearly these monomials all belong to $\beta$.  We define the $x_i$-\textbf{shadow} of an arrow $c^{\dlist}_{\jlist}$ to be the set of arrows $c^{\dlist}_{\jlist^{\prime}}$ such that $\xlist^{\jlist^{\prime}}$ is in the $x_i$-column of $\xlist^{\jlist}$, and 
\[
    x_i\text{\textbf{-height}}(c^{\dlist}_{\jlist}) = x_i\text{-degree}(\xlist^{\dlist}).
\]
If $c^{\dlist}_{\jlist}$ is a standard arrow for $x_i$, we define its $x_i$-\textbf{offset} to be the $(r-1)$-tuple $v$ of non-negative integers obtained by deleting the $i$-th (negative) component of the vector $\jlist - \dlist$.  

We say that an $x_i$-standard arrow $c^{\dlist}_{\jlist}$ can be \textbf{advanced} if either $c^{\dlist}_{\jlist}$ $\sim$ $0$ or $c^{\dlist}_{\jlist}$ $\sim$ $c^{\dlist_1}_{\jlist_1}$ with the latter arrow having strictly smaller $x_i$-height.  In each case, the idea is that $c^{\dlist}_{\jlist}$ can be translated so that its tail moves closer to the hyperplane ($x_i$-degree $=$ $0$), provided that we allow the head to exit the first orthant across this hyperplane in the first case.  Note that the head of an $x_i$-\sa\ arrow can only be translated out of the first orthant across the hyperplane ($x_i$-degree $=$ $0$), since the tail must remain within the first orthant (and outside of $\beta$) during translation.

\begin{lem} \label{lem:shadadv}
    Suppose that $c^{\dlist}_{\jlist}$ is a standard arrow for $x_i$.  Then:
\begin{enumerate}
    \item[(a)] All of the arrows in the $x_i$-shadow of $c^{\dlist}_{\jlist}$ are $x_i$-standard and have equal offsets.
    \item[(b)]  If $c^{\dlist}_{\jlist}$ can be advanced, then so can all the arrows in its $x_i$-shadow. 
    \item[(c)] If $c^{\dlist}_{\jlist}$ $\nsim 0$ can be advanced, then we can advance this arrow until we reach a minimal \sa\ (for $x_i$) arrow $c^{\dlist^{\prime}}_{\jlist^{\prime}}$ $\sim$  $c^{\dlist}_{\jlist}$  that  cannot be advanced.
\end{enumerate}
\end{lem}

\emph{Proof:}
     Statement (a) is immediate.  Statement (b) follows from the observation that if $c^{\dlist}_{\jlist}$ can be translated to  $c^{\blist}_{\jlist_1}$, then  every arrow in the shadow of $c^{\dlist}_{\jlist}$ is simultaneously translated ``in parallel'' either to an arrow in the shadow of $c^{\blist}_{\jlist_1}$ or out of the first orthant (across the hyperplane ($x_i$-degree $=$ $0$)). Statement (c) follows easily from Lemma \ref{lem:msa}: after advancing $c^{\dlist}_{\jlist}$ to $c^{\dlist_1}_{\jlist_1}$ of smaller $x_i$-height, we translate by degree-reducing steps to a \msa\ $c^{\dlist^{\prime}}_{\jlist^{\prime}}$, and repeat as often as necessary.  
\qedsymbol

\subsection{Shadow promotion} \label{SS:shadprom}
The simultaneous translation and advancement of the arrows in a shadow leads to a process of arrow replacement that is the key to the proof of Theorem \ref{thm:S}; we call this process, which we proceed to describe, \emph{shadow promotion}.
Suppose that $c^{\dlist}_{\jlist}$ is an $x_i$-standard arrow that can be advanced, which implies by statement (b) of Lemma \ref{lem:shadadv} that every arrow in the $x_i$-shadow of $c^{\dlist}_{\jlist}$ can be advanced as well.  Consider a sequence of translation steps that advances $c^{\dlist}_{\jlist}$, and let $c^{\dlist_1}_{\jlist_1}$ be the first position in the sequence of steps from which it is possible to move one step in the decreasing $x_i$-direction either to reach an arrow of $x_i$-height  $<$ $x_i$-height of $c^{\dlist}_{\jlist}$ or to move the head of the arrow out of the first orthant; put another way, $c^{\dlist_1}_{\jlist_1}$ is the first position such that
\[
   \xlist^{\dlist_1}/x_i \in I_{\beta}, \text{ and } x_i\text{-degree}( \xlist^{\dlist_1}) = x_i\text{-degree}( \xlist^{\dlist}).
\]
Note that
\begin{equation} \label{eqn:xidegs}
     x_i\text{-degree}(\xlist^{\jlist_1}) =  x_i\text{-degree}(\xlist^{\jlist}).
\end{equation}

Let $\xlist^{\blist}$  be a minimal generator of $I_{\beta}$ that divides $\xlist^{\dlist_1}/x_i$; it is clear that
\[
    x_i\text{-degree}(\xlist^{\blist})\ <\   x_i\text{-degree}(\xlist^{\dlist}).
\]
We must have that 
\begin{equation} \label{eqn:ineq1}
    x_i\text{-degree}(\xlist^{\blist})\ >\ x_i\text{-degree}(\xlist^{\jlist_1}),
\end{equation}
since otherwise $ x_i\text{-degree}(\xlist^{\blist})\ \leq\ x_i\text{-degree}(\xlist^{\jlist_1})$ holds in addition to 
\begin{equation} \label{eqn:ineq2}    
        x_k\text{-degree}(\xlist^{\blist}) \leq x_k\text{-degree}(\xlist^{\dlist_1}) \leq x_k\text{-degree}(\xlist^{\jlist_1})\ \text{ for all } k \neq i,
\end{equation}
where the last inequality holds because $c^{\dlist}_{\jlist}$ and its translate $c^{\dlist_1}_{\jlist_1}$ are \sa\ arrows for $x_i$.  In other words, $\xlist^{\blist}$ $|$ $\xlist^{\jlist_1}$, which is a contradiction since $\xlist^{\jlist_1}$ is a basis monomial and $\xlist^{\blist}$ is not.  The contradiction establishes (\ref{eqn:ineq1}); in particular, we have that $x_i$-degree$(\xlist^{\blist})$ $>$ $0$.

Furthermore, (\ref{eqn:ineq2}) shows that we can translate the arrow $c^{\dlist_1}_{\jlist_1}$ in degree-reducing directions, excluding the $x_i$-th, to reach an arrow $c^{\dlist_2}_{\jlist_2}$ such that $\xlist^{\dlist_2}$ and $\xlist^{\blist}$ differ only in $x_i$-degree. (Because the arrows involved are $x_i$-standard, the head of the arrow can never leave the first orthant during any of these steps.)  We then have that 
\[
    x_i\text{-degree}(\xlist^{\jlist_2})\ = \ x_i\text{-degree}(\xlist^{\jlist_1})\ < \ x_i\text{-degree}(\xlist^{\blist});
\]
therefore, the arrow $c^{\blist}_{\jlist_2}$ is a \msa\ for $x_i$ having the same offset as $c^{\dlist_2}_{\jlist_2}$, $c^{\dlist_1}_{\jlist_1}$, and $c^{\dlist}_{\jlist}$, and having the same number of arrows in its shadow as $c^{\dlist_1}_{\jlist_1}$ and $c^{\dlist}_{\jlist}$ do (recall (\ref{eqn:xidegs})).  We will call the shadow of $c^{\blist}_{\jlist_2}$ the \textbf{promotion image} of the shadow of the original \sa\ arrow $c^{\dlist}_{\jlist}$; the promotion image has the same number of arrows as the original shadow, and the arrows in the promotion image are \msa s for $x_i$ having strictly smaller $x_i$-height than the original arrows.  Indeed, we can view shadow promotion as the replacement of every arrow $c^{\dlist}_{\jlist^{\prime}}$ in the shadow of $c^{\dlist}_{\jlist}$ with the arrow $c^{\blist}_{\jlist^{\prime}_2}$ in the shadow of $c^{\blist}_{\jlist_2}$ such that
\begin{equation} \label{E:headDeg}
    x_i\text{-degree}(\xlist^{\jlist^{\prime}_2})\ =\ x_i\text{-degree}(\xlist^{\jlist^{\prime}}). 
\end{equation}
  Since the promotion image depends on the path chosen to advance $c^{\dlist}_{\jlist}$, it is not unique in general.

For example, in Figure \ref{fig:arroweg}, the arrow $c^{(2,0)}_{(0,1)}$ is a \sa\ arrow for $x_1$ that can be advanced; its shadow is the singleton set $\{ c^{(2,0)}_{(0,1)} \}$.  The promotion image of this shadow is easily seen to be $\{ c^{(1,1)}_{(0,2)} \}$.

\subsection{Iterated shadow promotion}  \label{SS:isp}

     Select one of the variables $x_i$, and recall that $x_i^{w_i}$ $=$ $\xlist^{\dlist}$ is the least power of $x_i$ that belongs to the monomial ideal $I_{\beta}$; it is clear that $\xlist^{\dlist}$ is a minimal generator of the ideal, which we will call the $i$-th \textbf{corner monomial} of $\beta$.  Note that \emph{every} arrow $c^{\dlist}_{\jlist}$ (with tail at the $i$-th corner monomial) is minimal and standard for $x_i$, since only the $x_i$-degree decreases when we move from the tail to the head of the arrow.  In this case, the offset of the arrow is the $(r-1)$-tuple $v$ obtained by deleting the $i$-th component of $\jlist$.  Let $\lis_0(i,v)$ denote the set of all (minimal $x_i$-standard) arrows having tail $\xlist^{\dlist}$ and offset $v$, and let $c^{\dlist}_{\jlist_1}$ be the arrow in $\lis_0(i,v)$ whose head $\xlist^{\jlist_1}$ has maximal $x_i$-degree $q(i,v)$; it follows that $\lis_0(i,v)$ is the $x_i$-shadow of the arrow $c^{\dlist}_{\jlist_1}$, and the number of arrows in $\lis_0(i,v)$ is $q(i,v)$ $+$ $1$.

     By iterating the process of shadow promotion, we can construct a set of $q(i,v)+1$ minimal $x_i$-standard arrows that have offset $v$ and cannot be advanced.  Indeed, if none of the arrows in $\lis_0(i,v)$ can be advanced, then $\lis_0(i,v)$ is itself the desired set, which we will denote $\lis(i,v)$.  Otherwise, one or more of the arrows in $\lis_0(i,v)$ can be advanced; in this case we let $c^{\dlist}_{\jlist_2}$ $\in$ $\lis_0(i,v)$ be the advanceable arrow whose head has maximal $x_i$-degree, so that the set of all advanceable arrows in $\lis_0(i,v)$ is the $x_i$-shadow of $c^{\dlist}_{\jlist_2}$, by Lemma \ref{lem:shadadv}.  We then promote the $x_i$-shadow of $c^{\dlist}_{\jlist_2}$, as described in the previous subsection, and denote the (not necessarily unique) promotion image by $P_0(i,v)$; this permits us to form a new set of minimal $x_i$-standard arrows of offset $v$:
\[
    \lis_1(i,v) = \left(\lis_0(i,v) \setminus (x_i\text{-shadow of } c^{\dlist}_{\jlist_2})\right)\ \cup\ P_0(i,v).
\]

Since the number of arrows in the promotion image is equal to the number of arrows in the shadow being promoted, it is clear that the number of arrows in $\lis_1(i,v)$ is equal to the number of arrows in $\lis_0(i,v)$.  Furthermore, the $x_i$-height of the arrows in $P_0(i,v)$ is strictly less than the height of the arrows in $\lis_0(i,v)$, which is $w_i$.  If none of the arrows in $\lis_1(i,v)$ can be advanced, then $\lis_1(i,v)$ is the desired set $\lis(i,v)$; otherwise, the set of advanceable arrows in $\lis_1(i,v)$ is equal to the shadow of an advanceable arrow $c^{\dlist^{\prime}}_{\jlist^{\prime}}$ $\in$ $P_0(i,v)$.  Replacing this shadow by its promotion image $P_1(i,v)$, we obtain yet another set of minimal $x_i$-standard arrows of offset $v$:
\[
    \lis_2(i,v) = \left(\lis_1(i,v) \setminus (x_i\text{-shadow of } c^{\dlist^{\prime}}_{\jlist^{\prime}})\right)\ \cup\ P_1(i,v).
\]
Continuing in this way, we eventually arrive at the desired set $\lis(i,v)$ of minimal $x_i$-standard arrows of offset $v$, none of which can be advanced.  The process must terminate because at each stage the $x_i$-height of the arrows that can be advanced is strictly less than the corresponding height at the previous stage, and this height cannot decrease to 0.

\subsection{Proof of Theorem \ref{thm:S}} \label{SS:Sproof}

For each of the variables $x_i$, $1$ $\leq$ $i$ $\leq$ $r$,  we will construct a set $\lis(i)$ consisting of $n$ minimal standard arrows for $x_i$ such that no two of the arrows in $\lis(i)$ are translation-equivalent to each other, and none is $\sim$ $0$.  But then the same is true of
\[
    \lis\ = \ \bigcup_{i=1}^r \lis(i),
\]
since $x_i$- and $x_j$-standard arrows cannot be translation-equivalent if $i$ $\neq$ $j$; in particular, the cardinality of $\lis$ is $r n$.  Theorem \ref{thm:lindepconds} now yields that the $\gf$-span of $\lis$ in $M/M^2$ has dimension $r n$, as desired.

It remains to construct the sets $\lis_i$, but this is not difficult.  Simply form the sets $\lis(i,v)$ for every possible offset $v$ of a minimal $x_i$-standard arrow having tail the $i$-th corner monomial $x_i^{w_i}$, as in Subsection \ref{SS:isp}, and let 
\[
    \lis(i)\ = \ \bigcup_{\{\text{offsets } v\}} \lis(i,v). 
\]
  It is clear that every monomial $\xlist^{\jlist}$ $\in$ $\beta$ is the head of an arrow in one of the initial sets $\lis_0(i,v)$.  Furthermore, the sets $\lis(i,v)$ are pairwise disjoint, since arrows of different offsets cannot be equal.  Therefore, counting elements in the various sets, we find that
\[
  \begin{array}{rcl}
       | \lis(i)| & = & \sum_{\{\text{offsets } v\}} |\lis(i,v)| \\
           {}     & = & \sum_{\{\text{offsets } v\}} |\lis_0(i,v)| \\
           {}     & = & |\beta| = n.
  \end{array}
\]
By construction, the arrows $c^{\blist}_{\jlist}$ $\in$ $\lis(i,v)$ are minimal $x_i$-standard arrows of offset $v$ that cannot be advanced; in particular, we have that $c^{\blist}_{\jlist}$ $\nsim$ $0$.  Therefore, the set $\lis(i)$ consists of $n$ minimal $x_i$-standard arrows, none of which are $\sim$ $0$.  We must now show that no two distinct arrows in $\lis(i)$ are translation-equivalent to one another.  To this end, let $c^{\blist_1}_{\jlist_1}$ and $c^{\blist_2}_{\jlist_2}$ be distinct arrows in $\lis(i)$.  If these arrows have different offsets, then they cannot possibly be translation-equivalent, so suppose that they each have offset $v$; that is, $c^{\blist_1}_{\jlist_1}$, $c^{\blist_2}_{\jlist_2}$ $\in$ $\lis(i,v)$.  If the two arrows have different heights (\ie, different tails), then, since neither arrow can be advanced, it is again clear that they cannot be translation-equivalent.  If the two arrows have the same height, then they have the same tail; therefore, if they were translation-equivalent, they would be equal, which we are assuming is not the case.  We conclude that distinct arrows in $\lis(i)$ cannot be translation-equivalent, and the proof is complete.
\qedsymbol

\begin{rem} \label{rem:noadvance}
    It is clear that none of the arrows in the set $\lis$ constructed in Theorem \ref{thm:S} can be advanced. 
\end{rem}

\begin{rem} \label{rem:notunique}
    Since the promotion image of a shadow is not necessarily unique, neither are the sets $\lis(i)$ and $\lis$.  
\end{rem}

\subsection{Summary and terminology} \label{SS:sumterm}

We will call a set $\lis(i)$ constructed as in Subsection \ref{SS:Sproof} a \textbf{standard $x_i$-sub-bunch}, and the union 
\[
    \lis\ = \ \bigcup_{i=1}^r \lis(i),
\]
a \textbf{standard bunch}, of arrows for $\beta$. 
These sets have the following properties:
\begin{description}
  \item[$\lis(i)$]    
     contains $n$ $=$ $|\beta|$ minimal $x_i$-standard arrows; has maximal rank (mod $M^2$); and no arrow in the set can be advanced.
  \item[$\lis$]  is the union of sets $\lis(i)$, $1$ $\leq$ $i$ $\leq$ $r$, and accordingly: contains $r n$ minimal standard arrows; has maximal rank (mod $M^2$); and no arrow in the set can be advanced.
\end{description} 

We will have occasion to consider more general sets of arrows $\lis^{\prime}(i)$ that we call \textbf{near-standard} $x_i$\textbf{-sub-bunches}: these are similar to standard $x_i$-sub-bunches in that they have cardinality $n$, have maximal rank (mod $M^2$), and consist of $x_i$-standard arrows that cannot be advanced; they differ in that the arrows they contain need not be \emph{minimal}.  For example, one way to obtain a near-standard $x_i$-sub-bunch is to replace some or all of the arrows $c^{\dlist}_{\jlist}$ $\in$ $\lis(i)$ (a standard $x_i$-sub-bunch) with arrows 
\[
    c^{\dlist^{\prime}}_{\jlist^{\prime}} \sim c^{\dlist}_{\jlist},\ \text{ such that }\ x_i\text{-height}({c^{\dlist^{\prime}}_{\jlist^{\prime}}}) =  x_i\text{-height}({c^{\dlist}_{\jlist}}).
\]
We will call the union 
\[
    \lis^{\prime}\ = \ \bigcup_{i=1}^r \lis^{\prime}(i),
\]
 of near-standard sub-bunches $\lis^{\prime}(i)$ a \textbf{near-standard bunch} of arrows for $\beta$; it is clear that $\lis^{\prime}$ consists of $r n$ standard arrows that cannot be advanced, and has maximal rank (mod $M^2$).  
\bigskip

We have the following useful corollaries of the proof of Theorem \ref{thm:S}:

\begin{cor} \label{cor:degofffact} 
    Let $c^{\dlist}_{\jlist}$ be an $x_i$-standard arrow for the basis set $\beta$.  Then every standard $x_i$-sub-bunch $\lis(i)$ of arrows for $\beta$ contains a unique arrow $c^{\dlist^{\prime}}_{\jlist^{\prime}}$ such that 
\[
    x_i\text{\emph{-degree}}(\xlist^{\jlist^{\prime}}) = x_i\text{\emph{-degree}}(\xlist^{\jlist})\ \text{ and\ \  \emph{offset}}(c^{\dlist^{\prime}}_{\jlist^{\prime}}) = \text{\emph{offset}}(c^{\dlist}_{\jlist}).
\]
\end{cor}

\emph{Proof:}
     Let 
\[
    p = x_i\text{-degree}(\xlist^{\dlist}),\ m = \xlist^{\dlist}/x_i^p,\ \text{and } \xlist^{\jlist_0} = \xlist^{\jlist} / m \in \beta.
\]
Then $\xlist^{\jlist_0}$ has offset $v$ from the corner monomial $x_i^{w_i}$ = $\xlist^{\blist}$, so $c^{\blist}_{\jlist_0}$ $\in$ $\lis_0(i,v)$, and is the only such arrow whose head has $x_i$-degree $=$ $x_i$-degree$(\xlist^{\jlist})$.  Under iterated shadow promotion, arrows are replaced by other arrows having the same offset and $x_i$-degree of the head (\ref{E:headDeg}), so we see that any standard $x_i$-sub-bunch $\lis(i)$ contains a unique arrow  $c^{\dlist_1}_{\jlist_1}$ as stated in the lemma, the arrow that ultimately replaces $c^{\blist}_{\jlist_0}$.
\qedsymbol

\begin{cor} \label{cor:crnrArrLem} 
    If $c^{\dlist}_{\jlist}$ is an $x_i$-\sa\ arrow that cannot be advanced, and whose tail is the corner monomial $\xlist^{\dlist}$ $=$ $x_i^{w_i}$, then $c^{\dlist}_{\jlist}$ belongs to every standard $x_i$-sub-bunch $\lis(i)$.
\end{cor}
     
\emph{Proof:}
    Clear from the construction. 
\qedsymbol


\section{Consequences of Theorem \ref{thm:S}} \label{S:conseq}

\subsection{Necessary and sufficient conditions for $t_{\beta}$ $\in$ $\Hn^n$ to be nonsingular} \label{SS:necsufcond}
By Proposition \ref{prop:freedim}, we know that the point $t_{\beta}$ $\in$ $\Hn^n$ is nonsingular if and only if the cotangent space $M/M^2$ has $\gf$-dimension $r n$.  Therefore, if $\lis^{\prime}$ is any standard or near-standard bunch of arrows for $\beta$ (see Subsection \ref{SS:sumterm} for terminology), we have that
\[
    t_{\beta} \text{ is nonsingular } \Leftrightarrow \lis^{\prime} \text{ spans } M/M^2 \Leftrightarrow \lis^{\prime} \text{ is a } \gf \text{-basis of } M/M^2.
\]
We will call $\beta$ a \textbf{smooth} basis set whenever any of these equivalent conditions holds.  

\begin{thm}  \label{thm:necsufcond}
     In order that the basis set $\beta$ be smooth, it is necessary and sufficient that the following conditions hold:
\begin{enumerate}
    \item[(a)] Every non-standard arrow $c^{\blist}_{\jlist}$ is translation-equivalent to $0$.
    \item[(b)] There exists a near-standard bunch of arrows $\lis^{\prime}$ such that for any \sa\ arrow $c^{\dlist}_{\jlist}$ $\nsim$ $0$, there exists an arrow  $c^{\dlist^{\prime}}_{\jlist^{\prime}}$ $\in$ $\lis^{\prime}$ such that $c^{\dlist}_{\jlist}$ $\sim$ $c^{\dlist^{\prime}}_{\jlist^{\prime}}$.
\end{enumerate}
\end{thm} 

\emph{Proof:}
    We first prove the necessity.  Let $\lis^{\prime}$ $=$ $\lis$ be the standard bunch given by Theorem \ref{thm:S}; since $\beta$ is assumed smooth, $\lis^{\prime}$ is a $\gf$-basis of $M/M^2$.  Therefore, given a non-standard arrow $c^{\blist}_{\jlist}$, we have that the $\gf$-span in $M/M^2$ of the enlarged set $\lis^{\prime}$ $\cup$ $\{ c^{\blist}_{\jlist} \}$ of $r n$ $+$ $1$ elements has dimension $r n$.  Since no two distinct arrows in $\lis^{\prime}$ are translation-equivalent, none is translation-equivalent to $0$, and no \sa\ arrow can be translation-equivalent to a non-\sa\ arrow, Theorem \ref{thm:lindepconds} yields that $c^{\blist}_{\jlist}$ $\sim$ $0$.  Now let $c^{\dlist}_{\jlist}$ $\nsim$ $0$ be a standard arrow.  If $c^{\dlist}_{\jlist}$ $\notin$ $\lis^{\prime}$, then again we must have that the $\gf$-span in $M/M^2$ of the enlarged set $\lis^{\prime}$ $\cup$ $\{ c^{\dlist}_{\jlist}\}$ has dimension $<$ $r n$ $+$ $1$; Theorem \ref{thm:lindepconds} now yields that $c^{\dlist}_{\jlist}$ $\sim$ $c^{\dlist^{\prime}}_{\jlist^{\prime}}$ for some $c^{\dlist^{\prime}}_{\jlist^{\prime}}$ $\in$ $\lis^{\prime}$.  

    It remains to prove the sufficiency; that is, assuming that conditions (a) and (b) hold, we must show that $\beta$ is smooth, for which it is enough to prove that $\lis^{\prime}$ spans the $\gf$-vector space $M/M^2$.  However, condition (a) implies that $M/M^2$ is spanned by the standard arrows, and condition (b) ensures that every standard arrow $c^{\dlist}_{\jlist}$ $\nsim$ $0$ is in the $\gf$-linear span of $\lis^{\prime}$, so we are done. 
\qedsymbol

\subsection{A sufficient condition for $t_{\beta}$ $\in$ $\Hn^n$ to be singular} \label{SS:nonsmthcond}

     It is often easy to identify basis sets $\beta$ for which the conditions (a) and (b) of Theorem \ref{thm:necsufcond} do \emph{not} hold, so that the corresponding point $t_{\beta}$ $\in$ $\Hn^n$ is singular; we now present one simple way to do this.  (Subsequent sections of the paper will be devoted to identifying basis sets $\beta$ for which conditions (a) and (b) \emph{do} hold.)

     If $\xlist^{\jlist}$ $\in$ $\beta$ is maximal among monomials in $\beta$ for the divisibility ordering, we call $\xlist^{\jlist}$ a \textbf{maximal basis monomial}.  If $\xlist^{\dlist}$ is a minimal generator of $I_{\beta}$ and $\xlist^{\jlist}$ is a maximal basis monomial, we call the arrow $c^{\dlist}_{\jlist}$ \textbf{rigid}, because we have

\begin{lem} \label{lem:rigid}
    An arrow $c^{\dlist}_{\jlist}$ such that $\xlist^{\dlist}$ is a minimal generator of $I_{\beta}$ and $\xlist^{\jlist}$ is a maximal basis monomial cannot be translated except to itself.  Such an arrow must belong to any $\gf$-basis $\mathcal{B}$ of $M/M^2$ consisting of arrows.
\end{lem}

\emph{Proof:}
    Since the tail of the arrow is a minimal generator of $I_{\beta}$, it is clear that the first step in any nontrivial translation of $c^{\dlist}_{\jlist}$ must be in a degree-increasing direction.  However, such a step would cause the head of the arrow to exit $\beta$ and enter $I_{\beta}$, which is forbidden.  Therefore, $c^{\dlist}_{\jlist}$ cannot be translated except to itself.

Now suppose that $\mathcal{B}$ is a $\gf$-basis of $M/M^2$ consisting of arrows, and that $c^{\dlist}_{\jlist}$ $\notin$ $\mathcal{B}$.  Then the $\gf$-span in $M/M^2$ of the set of arrows $\mathcal{B}$ $\cup$ $\{ c^{\dlist}_{\jlist} \}$ has dimension $<$ $|\mathcal{B}|$ $+$ $1$; moreover, $c^{\dlist}_{\jlist}$ $\nsim$ $0$, since  it has no nontrivial translations.  It then follows from Theorem \ref{thm:lindepconds} that $c^{\dlist}_{\jlist}$ $\sim$ $c^{\dlist_1}_{\jlist_1}$ for some $c^{\dlist_1}_{\jlist_1}$ $\in$ $\mathcal{B}$, but then $c^{\dlist}_{\jlist}$ $=$ $c^{\dlist_1}_{\jlist_1}$ by rigidity, which contradicts $c^{\dlist}_{\jlist}$ $\notin$ $\mathcal{B}$.  We conclude that $c^{\dlist}_{\jlist}$ $\in$ $\mathcal{B}$.
\qedsymbol

\begin{cor} \label{cor:rigid}
    If the basis set $\beta$ has a non-\sa\ rigid arrow $c^{\dlist}_{\jlist}$, then condition \rm{(a)} of Theorem \rm{\ref{thm:necsufcond}} is false, so $\beta$ is not smooth.
\end{cor}

\emph{Proof:}
    It is clear that $c^{\dlist}_{\jlist}$ is a non-\sa\ arrow that is not translation-equivalent to $0$.
\qedsymbol

\subsection{Example: $\beta$ $=$ $\{1, x_1, x_2, x_3 \}$} \label{SS:nonsmootheg}

\begin{figure} 
\begin{picture}(374,84)
        \put(156,13){$x_1$}
        \put(122,0){$\mathbf{x_1^2}$}
        \put(190,26){$1$}
        \put(224,26){$x_2$}
        \put(258,26){$\mathbf{x_2^2}$}
        \put(188,12){$\mathbf{x_1 x_2}$}
        \put(190,50){$x_3$}
        \put(190,74){$\mathbf{x_3^2}$}
        \put(152,37){$\mathbf{x_1 x_3}$}
        \put(214,50){$\mathbf{x_2 x_3}$}
        \put(178,39){\vector(4,-1){40}}
        \put(210,48){\vector(-3,-2){40}}
        \put(202,19){\vector(-1,4){7}}
        \put(190,22){\vector(-2,-1){60}} 
        \put(190,22){\vector(1,0){90}}   
        \put(190,22){\vector(0,1){70}}   
\end{picture}
\caption{The basis set $\beta$ $=$ $\{1, x_1, x_2, x_3 \}$, with the minimal generators of $I_{\beta}$ shown in boldface.  The arrows $c^{(1,0,1)}_{(0,1,0)}$, $c^{(0,1,1)}_{(1,0,0)}$, and $c^{(1,1,0)}_{(0,0,1)}$ are non-\sa\ and rigid, demonstrating that $t_{\beta}$ $\in$  $\Hn^4$ is singular.} \label{fig:nsrigid}
\end{figure}
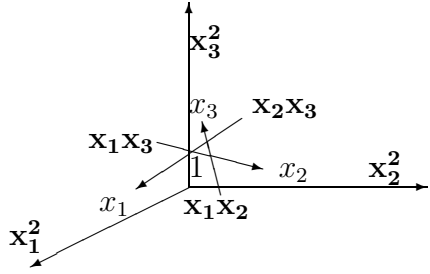

Figure \ref{fig:nsrigid} illustrates the basis set $\beta$ $=$ $\{1, x_1, x_2, x_3\}$, the smallest basis set in three variables that has non-\sa\ rigid arrows, and therefore fails to be smooth by Corollary \ref{cor:rigid}.

Going further, we recall that by Lemma \ref{lem:msa}, the $\gf$-vector space $M/M^2$ is spanned by minimal arrows.  However, in this case, there are $6 \cdot 4$ $=$ $24$ minimal arrows, 6 of which are $\sim 0$ (the ones with head $\xlist^{(0,0,0)} = 1$), and the other 18 of which are rigid, and so must belong to any basis of $M/M^2$ consisting of arrows, by Lemma \ref{lem:rigid}.  We conclude from this that the $\gf$-dimension of the cotangent space of $t_{\beta}$ is 18.  It is known (see, 
\eg, \cite[Sec.\ 1, p.\ 147]{Iarrob:Sitges}, or \cite[Sec.\ 5.1, p.\ 443]{Notari-Spreafico}) that $\Hn^4$ is irrreducible of dimension $3 \cdot 4$ $=$ $12$, so we reconfirm that $\beta$ is not smooth.

Finally, note that three of the 18 rigid minimal arrows are non-standard (see Figure \ref{fig:nsrigid}), leaving 15 rigid \msa s.  It is clear that any near-standard bunch $\lis^{\prime}$ (which contains $3 \cdot 4$ $=$ $12$ standard arrows) must exclude at least three of these rigid standard arrows, showing that $\beta$ also fails to satisfy condition (b) of Theorem \ref{thm:necsufcond}.  


\section{Thickening of basis sets} \label{S:thicken}

     In this section of the paper, we begin to identify certain families of smooth basis sets $\beta$.  Recall that a basis set is smooth if it is associated to a nonsingular point $t_{\beta}$ $\in$ $\Hn^n$, and therefore satisfies conditions (a) and (b) of Theorem \ref{thm:necsufcond}.  We first consider a natural way to obtain a smooth basis set in $r$ variables by ``thickening'' a smooth basis set in $r-1$ variables.
     
\subsection{Definition of thickening} \label{SS:thickdefn}

    Suppose that $\beta_0$ is a basis set of $n_0$ monomials in the variables $x_1$, $x_2$, \dots, $x_{r-1}$, and $w_r$ is a positive integer.  We define a basis set in $r$ variables as follows:
\begin{equation} \label{eqn:liftdef}
    \beta = \{ \xlist^{\jlist_0} \cdot x_r^{s} \mid \xlist^{\jlist_0} \in \beta_0,\ 0 \leq s \leq w_r-1\};
\end{equation}
it is indeed easy to check that $\beta$ is a basis set.  We say that $\beta$ is a \textbf{thickening} of $\beta_0$ from $r-1$ to $r$ variables.  Note that the number of monomials in $\beta$ is
\begin{equation} \label{eqn:liftcard}
    |\beta| = n = n_0 \cdot w_r.
\end{equation}

\subsection{Minimal generators of $I_{\beta}$} \label{SS:thickmg}

\begin{lem} \label{Lem:thick}
    Let $\beta_0$ and $\beta$ be as above.  Then a minimal generator of the ideal $I_{\beta}$ is either a minimal generator $\xlist^{\blist}$ of the ideal $I_{\beta_0}$ (and therefore not divisible by $x_r$) or else is the corner monomial $x_r^{w_r}$.  Furthermore, every arrow with tail $x_r^{w_r}$ is a \msa\ for $x_r$, and none of these arrows can be advanced.
\end{lem}

\emph{Proof:}
    Let $m$ $=$ $\xlist^{\dlist}$ be a minimal generator of $I_{\beta}$.  If $m$ is not divisible by $x_r$, then for any other variable $x_j$, $1$ $\leq$ $j$ $\leq$ $r-1$, we have that either $m/x_j$ is a monomial in $\beta_0$ or is undefined; it follows that $m$ is a minimal generator of $I_{\beta_0}$ $\subseteq$ $\gf[x_1, \dots, x_{r-1}]$.  If $m$ is divisible by $x_r$, then $m/x_r$ $\in$ $\beta$.  Therefore, $m/x_r$ $=$ $\xlist^{\jlist_0} \cdot x_r^{s}$, where $\xlist^{\jlist_0}$ $\in$ $\beta_0$ and $0$ $\leq$ $s$ $\leq$ $w_r-1$.  In fact, we must have $s$ $=$ $w_r-1$, since otherwise $m$ $\in$ $\beta$, which contradicts $m$ $\in$ $I_{\beta}$.  But then $x_r^{w_r}$ $|$ $m$, and since it is clear that $x_r^{w_r}$ is a minimal generator of $\beta$, we must have $m$ $=$ $x_r^{w_r}$.  As in Section \ref{SS:isp}, we have that every arrow having tail the corner monomial $x_r^{w_r}$ is a minimal $x_r$-standard arrow.  None of these arrows can be advanced, since $x_r^{w_r}$ is the only minimal generator that is divisible by $x_r$, and hence is the only minimal generator that can serve as the tail of an $x_r$-standard arrow; this completes the proof of the lemma.
\qedsymbol

\subsection{A standard bunch $\lis$ for $\beta$} \label{SS:thickS}

Let $\beta$ be the thickening (\ref{eqn:liftdef}) of $\beta_0$.  For the proof of our main result on thickenings, it is convenient to have available a standard bunch $\lis$ for $\beta$ that is closely related to a standard bunch $\lis_0$ for $\beta_0$.  We proceed to construct such a set of arrows.

Let $\pi$ denote the projection map taking monomials in $x_1$, \dots, $x_r$ to monomials in $x_1$, \dots, $x_{r-1}$ defined by 
\[
  m \stackrel{\pi}{\mapsto} m/x_r^{x_r\text{-deg}(m)}.
\]
  If $c^{\dlist}_{\jlist}$ is an arrow for $\beta$ such that $\pi(\xlist^{\dlist})$ $\notin$ $\beta_0$, then it is clear that
\[
    \pi(c^{\dlist}_{\jlist})\ = \ c^{\pi(\dlist)}_{\pi(\jlist)}
\]
is an arrow for $\beta_0$ that we will call the \textbf{projection} of $c^{\dlist}_{\jlist}$.  Furthermore, if we can translate $c^{\dlist}_{\jlist}$ to $c^{\dlist_1}_{\jlist_1}$ in such a way that the tails of the arrows in the translation path always project to monomials $\notin$ $\beta_0$, then the projections of the arrows define a translation path from $\pi(c^{\dlist}_{\jlist})$ to $\pi(c^{\dlist_1}_{\jlist_1})$ that only involves steps in the directions $x_1$, \dots, $x_{r-1}$.  In this case we will say that the translation path \textbf{projects} from $\beta$ to $\beta_0$.

For each $c^{\dlist_0}_{\jlist_0}$ such that 
\[
    x_r\text{-degree}(\xlist^{\dlist_0}) = 0\ \text{ and }\ x_r\text{-degree}(\xlist^{\jlist_0}) = 0 
\]
(that is, $c^{\dlist_0}_{\jlist_0}$ is an arrow for $\beta_0$), define
\[
    \text{tower}(c^{\dlist_0}_{\jlist_0}) = \{ c^{\dlist_0}_{\jlist^{\prime}}  \mid \xlist^{\jlist^{\prime}} = \xlist^{\jlist_0} \cdot x_r^s,\ 0 \leq s \leq w_r-1 \}.
\]
It is then clear that any translation of $c^{\dlist_0}_{\jlist_0}$ in directions other than the $x_r$-th can be applied to the entire tower ``in parallel''.  

\begin{lem} \label{lem:towerlem}
    For $x_i$ $\neq$ $x_r$, let $c^{\dlist}_{\jlist}$ be an $x_i$-\sa\ arrow for $\beta$ that can be advanced.  Then the projection
\[
    \pi(c^{\dlist}_{\jlist})\ =\ c^{\dlist_0}_{\jlist_0}
\]
is an $x_i$-\sa\ arrow for $\beta_0$ that can be advanced; consequently, all the arrows in  \emph{tower}$(c^{\dlist_0}_{\jlist_0})$ can be advanced in parallel.
\end{lem}

\emph{Proof:}
    We have that
\[
  x_r\text{-degree}(\xlist^{\dlist}) \leq x_r\text{-degree}(\xlist^{\jlist}) < r,
\]
so the tail $\xlist^{\dlist}$ must be divisible by a minimal generator of $I_{\beta_0}$, by Lemma \ref{Lem:thick}.  Therefore, the projection $c^{\dlist_0}_{\jlist_0}$ is defined.  Since the vector $\jlist_0$ $-$ $\dlist_0$ is the same as the original vector $\jlist$ $-$ $\dlist$ except in the $x_r$-component (which is $0$ for the former and $\geq$ $0$ for the latter), it follows easily that $c^{\dlist_0}_{\jlist_0}$ is $x_i$-\sa.  Furthermore, the translation path that advances $c^{\dlist}_{\jlist}$ consists entirely of $x_i$-standard arrows, so it projects from $\beta$ to $\beta_0$, implying that the projection $c^{\dlist_0}_{\jlist_0}$ can be advanced (as an arrow for $\beta_0$).  As observed earlier, all the arrows in tower$(c^{\dlist_0}_{\jlist_0})$ then advance in parallel.
\qedsymbol
\bigskip

We now select a particular standard bunch of arrows $\lis$ for $\beta$ (Theorem \ref{thm:S} ensures that at least one such set exists).  Recall that $\lis$ can be constructed as the union of standard $x_i$-sub-bunches $\lis(i)$, $1$ $\leq$ $i$ $\leq$ $r$; each $\lis(i)$ consists of $n$ \msa s for $x_i$ that cannot be advanced.  By the last statement of Lemma \ref{Lem:thick}, together with Corollary \ref{cor:crnrArrLem}, we have that $\lis(r)$ is the set of all $n$ arrows having tail $x_r^{w_r}$.  For any $x_i$ $\neq$ $x_r$, we can construct $\lis(i)$ as in Subsection \ref{SS:Sproof}, taking care to translate all towers in parallel, using Lemma \ref{lem:towerlem}.  That is, if $c^{\dlist_0}_{\jlist}$ is a \emph{minimal} $x_i$-\sa\ arrow that can be advanced, (so the tail $\xlist^{\dlist_0}$ is a minimal generator of $\beta_0$, by Lemma \ref{Lem:thick}), then the projection
\[
    \pi(c^{\dlist_0}_{\jlist})\ =\ c^{\dlist_0}_{\jlist_0}
\] 
can be advanced (as an arrow for $\beta_0$); whence, all the arrows in tower$(c^{\dlist_0}_{\jlist_0})$ (which contains the original arrow $c^{\dlist_0}_{\jlist}$) can be advanced in parallel.  With this understanding, we see that the construction (for $\beta$) of an $x_i$-standard sub-bunch $\lis(i)$ can be carried out as follows: First, restricting to the basis set $\beta_0$ and the variables $x_1$, \dots, $x_{r-1}$, choose a standard $x_i$-sub-bunch $\lis_0(i)$, as in the proof of Theorem \ref{thm:S}.  Then set
\[
    \lis(i) = \bigcup_{c^{\blist}_{\jlist_0} \in \lis_0(i)} \text{tower}(c^{\blist}_{\jlist_0}).
\]

Hence, by setting 
\[
    \lis_0 = \cup_{i=1}^{r-1}\lis_0(i),\ \lis = \cup_{i=1}^{r}\lis(i),
\]
we obtain a standard bunch $\lis$ for $\beta$ that ``lies over'' a standard bunch $\lis_0$ for $\beta_0$ in the sense that every arrow in $\lis(i)$ projects to an arrow in $\lis_0(i)$ for $1$ $\leq$ $i$ $\leq$ $r-1$.   

\subsection{Thickenings of smooth basis sets are smooth} \label{SS:thicksmth}

\begin{thm} \label{thm:thicksmth}
    If $\beta_0$ is a smooth basis set in $r-1$ variables, then for each integer $w_r$ $>$ $0$, the thickening $\beta$ \emph{(\ref{eqn:liftdef})} is a smooth basis set.
\end{thm} 

\emph{Proof:}
    It suffices to show that $\beta$ satisfies conditions (a) and (b) of Theorem \ref{thm:necsufcond}. 
    First suppose that $c^{\dlist}_{\jlist}$ is a non-\sa\ arrow for $\beta$.  If $x_r^{w_r}$ $|$ $\xlist^{\dlist}$, then we can translate the arrow $c^{\dlist}_{\jlist}$ in degree-decreasing steps so that its tail approaches $x_r^{w_r}$; since no non-\sa\ arrow can have tail $x_r^{w_r}$, we see that at some point the head of the arrow must exit the first orthant, demonstrating that $c^{\dlist}_{\jlist}$ $\sim$ $0$.  On the other hand, if $x_r^{w_r}$ $\ndiv$ $\xlist^{\dlist}$, then we can assume by Lemmas \ref{lem:msa} and \ref{Lem:thick} that we have translated our arrow to $c^{\blist}_{\jlist^{\prime}}$, with tail a minimal generator $\xlist^{\blist}$ of $I_{\beta_0}$ (we are done if the head of the arrow exits the first orthant during this translation).  Let $c^{\blist}_{\jlist_0}$ = $\pi(c^{\blist}_{\jlist^{\prime}})$, 
and observe that $c^{\blist}_{\jlist_0}$ is a non-\sa\ arrow for $\beta_0$.  Since $\beta_0$ is a smooth basis set, by hypothesis, we have that $c^{\blist}_{\jlist_0}$ $\sim$ $0$, using only translations in the first $r-1$ variable directions.  It is clear that the arrow $c^{\blist}_{\jlist^{\prime}}$ translates ``in parallel'' with $c^{\blist}_{\jlist_0}$, leading to the conclusion that $c^{\blist}_{\jlist^{\prime}}$ $\sim$ $0$.  Therefore $\beta$ satisfies condition (a) of Theorem \ref{thm:necsufcond}.

To prove that $\beta$ satisfies condition (b) of Theorem \ref{thm:necsufcond}, we use the standard bunches $\lis$ and $\lis_0$ constructed in Subsection \ref{SS:thickS}, and we let $c^{\dlist}_{\jlist}$ $\nsim$ $0$ be an $x_i$-\sa\ arrow for $\beta$.  If $i$ $=$ $r$, then we can translate in degree-reducing steps until we reach an arrow whose tail is the corner monomial $\xlist^{\dlist_1}$ $=$  $x_r^{w_r}$, so that $c^{\dlist}_{\jlist}$ $\sim$ $c^{\dlist_1}_{\jlist_1}$ $\in$ $\lis(r)$ $\subseteq$ $\lis$.  If $i$ $\neq$ $r$, then we can translate by degree-reducing steps until we reach an arrow $c^{\blist}_{\jlist_1}$ whose tail is one of the minimal generators $\xlist^{\blist}$ of $I_{\beta_0}$.
Let $c^{\blist}_{\jlist_0}$ $=$ $\pi(c^{\blist}_{\jlist_1})$, and note that $c^{\blist}_{\jlist_1}$ $\in$ tower$(c^{\blist}_{\jlist_0})$. 
Since $\beta_0$ is a smooth basis set, we know that there exists an arrow $c^{\blist^{\prime}}_{\jlist^{\prime}_0}$ $\in$ $\lis_0(i)$ such that $c^{\blist}_{\jlist_0}$ $\sim$ $c^{\blist^{\prime}}_{\jlist^{\prime}_0}$.  The translations only involve the variable directions $x_1$, \dots, $x_{r-1}$, and will carry $c^{\blist}_{\jlist_1}$ ``in parallel'' to an arrow $c^{\blist^{\prime}}_{\jlist^{\prime}_1}$ $\in$ tower$(c^{\blist^{\prime}}_{\jlist^{\prime}_0})$ $\subseteq$ $\lis(i)$.  We conclude that $\beta$ satisfies condition (b) of Theorem \ref{thm:necsufcond}.  This completes the proof that the thickening $\beta$ of $\beta_0$ is a smooth basis set.
\qedsymbol

\begin{rem} \label{rem:thickconv} 
     The converse of Theorem \ref{thm:thicksmth} also holds; that is, if $\beta$ is a smooth basis set that is a thickening (\ref{eqn:liftdef}) of $\beta_0$, then $\beta_0$ is smooth.  We leave this as an exercise for the reader. 
\end{rem}

\subsection{Example: ``Boxes'' are smooth basis sets} \label{SS:boxsmth}

Let $w_1$, \dots, $w_r$ be positive integers, and let 
\begin{equation} \label{E:rbox}
    \rbox(w_1,\ldots, w_r)\ = \ \{ x_1^{d_1} x_2^{d_2} \dots x_r^{d_r} \mid 0 \leq d_i < w_i,\ 1 \leq i \leq r\},
\end{equation}
which is clearly a basis set containing
\[
    |\rbox(w_1, \ldots, w_r)|\ =\ n\ =\ \prod_{i=1}^r w_i
\]
monomials; for obvious reasons, we call this type of basis set a \textbf{box}.  The following useful results are immediate:

\begin{lem} \label{lem:boxmingen}
    The minimal generators of $I_{\rbox(w_1, \ldots, w_r)}$ are the corner monomials $x_i^{w_i}$, $1$ $\leq$ $i$ $\leq$ $r$.
\qedsymbol
\end{lem}

\begin{lem} \label{lem:boxlcm}
    If $m_1$, $m_2$ $\in$ $\rbox(w_1, \dots, w_r)$, then 
\[
    \lcm(m_1, m_2) \in \rbox(w_1, \ldots, w_r).\ \ \text{\qedsymbol}
\]
\end{lem}

We could give an easy proof of the next proposition using Lemma \ref{lem:boxmingen}, but instead we offer a proof based on thickenings.

\begin{prop} \label{prop:boxsmth}
    The box $\rbox(w_1,\ldots, w_r)$ is a smooth basis set.
\end{prop}

\emph{Proof:}
    It is clear that $\rbox(w_1, \ldots, w_r)$ is a thickening of $\rbox(w_1, \ldots, w_{r-1})$, so if the latter is smooth, then so is the former, by Theorem \ref{thm:thicksmth}.  The desired result therefore follows by induction, provided that the result holds in the base case: boxes $\beta$ = $\rbox(w_1)$ in one variable.  But in this case, the monomial ideal $I_{\beta}$ has just one minimal generator $x_1^{w_1}$; therefore, there are exactly $w_1$ $=$ $|\beta|$ minimal arrows, all of which are standard and cannot be advanced; in fact, the set of minimal arrows is equal to the standard bunch $\lis$ $=$ $\lis(1)$, in the notation of the proof of Theorem \ref{thm:S}.  Since $M/M^2$ is spanned by minimal arrows (Lemma \ref{lem:msa}), we conclude that $\lis$ is a $\gf$-basis of $M/M^2$  $\Rightarrow$ $\beta$ is smooth, as desired.  
\qedsymbol


\section{Truncation of basis sets} \label{S:trunc} 

Truncation is another natural way to obtain a basis set from a given basis set.  In this section we will define the truncation operation and establish one of its key properties: If $\beta$ is such that every non-standard arrow $c^{\dlist}_{\jlist}$ is translation-equivalent to $0$ (that is, $\beta$ satisfies condition (a) of Theorem \ref{thm:necsufcond}), then any truncation of $\beta$ also has this property; in particular, this is so if $\beta$ is smooth.  Under an additional hypothesis, one can further prove that certain truncations of a smooth basis set are smooth (that is, also satisfy condition (b) of Theorem \ref{thm:necsufcond}); we discuss this in the next section. 

\subsection{Definition of truncation} \label{SS:truncdef}

Let $\beta$ be a basis set of $n$ monomials in the variables $x_1$, \dots, $x_r$.  Choose one of the variables $x_j$, and a positive integer $h$ such that $\beta$ contains at least one monomial that is divisible by $x_j^h$.  Then we define the $x_j$-\textbf{truncation of} $\beta$ \textbf{at height} $h$ to be the basis set
\begin{equation} \label{E:truncdef}
  \beta_t(x_j,h)\ =\ \beta_t\ =\ \{ m \mid x_j^h \cdot m \in \beta \}.
\end{equation}
That is, $\beta_t$ is obtained by discarding the monomials in $\beta$ with $x_j$-degree $<$ $h$, and dividing the remaining monomials by $x_j^h$.  Figure \ref{Fig:trunc} provides an example.  We write
\[
    n\ =\ |\beta|,\ n_t\ = \ |\beta_t|
\]
for the number of monomials in $\beta$, $\beta_t$, respectively.
\begin{figure} 
\begin{picture}(372,96)
        \put(102,42){
               \mbox{$\begin{array}{cccccc}
                x_2^6 & {} & {} & {} & {} & {}\\
                x_2^5 &{x_1 x_2^5} & {} &{} &{} & {}   \\
                x_2^4 & x_1 x_2^4 & x_1^2 x_2^4 & x_1^3 x_2^4 & {} & {}\\
                x_2^3 & x_1 x_2^3 & x_1^2 x_2^3 & x_1^3 x_2^3 & {} & {}\\
                x_2^2 & x_1 x_2^2 & x_1^2 x_2^2 & x_1^3 x_2^2 & x_1^4 x_2^2 & {}\\
                x_2 & x_1 x_2 & x_1^2 x_2 & x_1^3 x_2 & x_1^4 x_2 & {}\\
                1 & x_1 & x_1^2 & x_1^3 & x_1^4 & x_1^5 
                     \end{array}
                    $}
                   }
        \put(195,96){\line(0,-1){102}}
        \put(108,66){\line(1,0){21}}
        \put(108,66){\line(0,-1){70}}
        \put(108,-4){\line(1,0){84}}
        \put(192,-4){\line(0,1){14}}
        \put(192,10){\line(-1,0){30}}
        \put(162,10){\line(0,1){28}}
        \put(162,38){\line(-1,0){33}}
        \put(129,38){\line(0,1){28}}
\end{picture}
\caption{The small basis set enclosed by the polygon is the truncation at $x_1$-degree 3 of the basis set consisting of all of the displayed monomials.  It is obtained by dividing the monomials to the right of the vertical line by $x_1^3$.} \label{Fig:trunc}
\end{figure}

\subsection{Minimal generators of $I_{\beta_t}$} \label{SS:mgIbt}

We denote by $w_i$ (resp.\ $w^{(t)}_i$) the $x_i$-width of $\beta$ (resp.\ $\beta_t$); that is, $x_i^{w_i}$ (resp.\ $x_i^{w^{(t)}_i}$) is the $i$-th corner monomial of $I_{\beta}$ (resp.\ $I_{\beta_t}$).  
We have the following
\begin{lem} \label{lem:mgtrunc}
    Let $\beta$ be a basis set, and $\beta_t$ the $x_j$-truncation of $\beta$ at height $h$.  Then 
  \begin{enumerate}
    \item[(a)] If $m$ is a minimal generator of $I_{\beta}$ such that $x_j$-degree$(m)$ $\geq$ $h$, then $m/x_j^h$ $=$ $m_t$ is a minimal generator of $I_{\beta_t}$; consequently, for all $x_k$ $\neq$ $x_j$, we have that
\[
    x_k\text{-degree}(m)\ =\  x_k\text{-degree}(m_t)\ \leq\ w^{(t)}_k.
\]

    \item[(b)] The minimal generators $m$ of $I_{\beta}$ that have $x_j$-degree $>$ $h$ are in bijective correspondence with the minimal generators $m_t$ $=$ $m/x_j^h$ of $I_{\beta_t}$ that are divisible by $x_j$.
 
    \item[(c)] If $m_t^{\prime}$ is a minimal generator of $I_{\beta_t}$ that is not divisible by $x_j$, then $x_j^h \cdot m_t^{\prime}$ $=$ $m^{\prime}$ is an $x_j$-multiple of a minimal generator of $I_{\beta}$.

    \item[(d)] $w_j$ $=$ $w_j^{(t)} + h$ $>$ $h$, and for all $x_k$ $\neq$ $x_j$, $w_k$ $\geq$ $w^{(t)}_k$.
  \end{enumerate}
\end{lem}

\emph{Proof:}
    We begin by noting that a monomial $m$ $\in$ $I_{\beta}$ (resp.\ $I_{\beta_t}$) is a minimal generator if and only if for each variable $x_i$, either $x_i$ $\ndiv$ $m$ or $m/x_i$ $\in$ $\beta$ (resp.\ $\beta_t$).  

Let $m$ and $m_t$ be as in assertion (a).  Since $x_j^h \cdot m_t$ $=$ $m$ $\notin$ $\beta$, we have that $m_t$ $\notin$ $\beta_t$ $\Rightarrow$ $m_t$ $\in$ $I_{\beta_t}$.  For any variable $x_k$ $\neq$ $x_j$, note that
\[
    x_k\, |\, m_t \Rightarrow x_k\, |\, m \Rightarrow m/x_k \in \beta;
\]
 moreover,
\[
     x_j\text{-degree}(m/x_k) = x_j\text{-degree}(m) \geq h;
\]
therefore, $m_t/x_k$ $\in$ $\beta_t$.  If $x_j$ $|$ $m_t$, then 
\[
    x_j\text{-degree}(m) = \left(x_j\text{-degree}(m_t) + h\right)\ >\ h; 
\]
whence,
\[
    x_j\text{-degree}(m/x_j) \geq h, \text{ and as before } m/x_j \in \beta,
\]
which yields $m_t/x_j$ $\in$ $\beta_t$.  This completes the proof that $m_t$ is a minimal generator of $I_{\beta_t}$, and the stated consequence is immediate.

       Since the map $m$ $\mapsto$ $m_t$ is injective, to prove assertion (b), it remains to show that if $m^{\prime}_t$ is a minimal generator of $\beta_t$ that is divisible by $x_j$, then $x_j^h \cdot m^{\prime}_t$ $=$ $m^{\prime}$ is a minimal generator of $\beta$.  It is clear that $m^{\prime}$ $\in$ $I_{\beta}$, since $m^{\prime}_t$ $\notin$ $\beta_t$.  Furthermore, for any variable $x_i$, we have that
\begin{equation} \label{E:someimps}
    x_i\, |\, m^{\prime}\ \Rightarrow\ x_i\, |\, m^{\prime}_t\ \Rightarrow\ m^{\prime}_t/x_i \in \beta_t\ \Rightarrow\ m^{\prime}/x_i \in \beta
\end{equation}
(the first implication is trivial for $x_i$ $=$ $x_j$), so $m^{\prime}$ is indeed a minimal generator of $\beta$.

     Let now $m_t^{\prime}$, $m^{\prime}$ be as in assertion (c).  As in the preceding paragraph, $m^{\prime}$ $\in$ $I_{\beta}$.  We can translate $m^{\prime}$ to a minimal generator $m$ of $I_{\beta}$ by degree-reducing steps, but since the implications (\ref{E:someimps}) hold for all $x_i$ $\neq$ $x_j$, we see that we can only move in the direction of decreasing $x_j$-degree; whence, $m^{\prime}$ is an $x_j$-multiple of $m$, as desired. 
    
Turning to assertion (d), it is clear that
\[
  \begin{array}{rcl}
    x_{\ell}\, |\, x_j^{w^{(t)}_j+h}  & \Rightarrow  & x_{\ell} = x_j \\
     {}                            & \Rightarrow  &x_j^{w^{(t)}_j+h} / x_{\ell} = x_j^{w^{(t)}_j+h-1} = x_j^h \cdot x_j^{w^{(t)}_j-1}  \in x_j^{h} \cdot \beta_t \subseteq \beta,
  \end{array} 
\]
which implies that $x_j^{w^{(t)}_j+h}$ is a minimal generator of $I_{\beta}$; therefore, $w_j$ $=$ $w^{(t)}_j+h$.  Finally, for any $x_k$ $\neq$ $x_j$, we have that
\[
    \begin{array}{rcl}
        x_k^{w^{(t)}_k-1} \in \beta_t & \Rightarrow & x_j^h \cdot x_k^{w^{(t)}_k-1} \in \beta\\
                    {}            & \Rightarrow & x_k^{w_k} \ndiv\, (x_j^h \cdot x_k^{w^{(t)}_k-1})\\
                    {}            & \Rightarrow & w_k \geq w^{(t)}_k,
    \end{array}
\]
and the proof is complete.
\qedsymbol

\subsection{Lifting arrows from $\beta_t$ to $\beta$} \label{SS:truncS}

Let $\beta$ and $\beta_t$ be as above, and let $c^{\dlist}_{\jlist}$ be an arrow for $\beta_t$.  Since $\xlist^{\jlist}$ $\in$ $\beta_t$ (resp.\ $\xlist^{\dlist}$ $\notin$ $\beta_t$), we have that $\xlist^{\jlist^{\prime}}$ $=$ $\xlist^{\jlist}\cdot x_j^{h}$ $\in$ $\beta$ (resp.\ $\xlist^{\dlist^{\prime}}$ $=$ $\xlist^{\dlist}\cdot x_j^{h}$ $\notin$ $\beta$); therefore,  $c^{\dlist^{\prime}}_{\jlist^{\prime}}$ is an arrow for $\beta$. We will call $c^{\dlist^{\prime}}_{\jlist^{\prime}}$ the \textbf{lifting} of $c^{\dlist}_{\jlist}$ from $\beta_t$ to $\beta$, and say that $c^{\dlist^{\prime}}_{\jlist^{\prime}}$ \textbf{descends} to $c^{\dlist}_{\jlist}$.  We extend this terminology to sets of arrows in the obvious way.  It is clear that the arrows $c^{\dlist}_{\jlist}$ and $c^{\dlist^{\prime}}_{\jlist^{\prime}}$ have the same vector \[
    \jlist - \dlist\ =\ \jlist^{\prime} - \dlist^{\prime}.
\]

\subsection{Non-standard arrows on $\beta$ and $\beta_t$} \label{SS:zeroing}

\begin{thm} \label{thm:truncnsa} Let $\beta$ be a basis set with the property that every one of its non-standard arrows is translation-equivalent to $0$ (condition \emph{(a)} of Theorem \emph{\ref{thm:necsufcond}}), and let $\beta_t$ be the $x_j$-truncation of $\beta$ at height $h$ \emph{(\ref{E:truncdef})}.  Then $\beta_t$ also has this property; that is, every non-standard arrow for $\beta_t$ is translation-equivalent to $0$.
\end{thm}

\emph{Proof:} 
Let $c^{\dlist}_{\jlist}$ be a non-standard arrow for $\beta_t$.  We must show that $c^{\dlist}_{\jlist}$ $\sim$ $0$.  To do this, we consider the lifting $c^{\dlist^{\prime}}_{\jlist^{\prime}}$ of $c^{\dlist}_{\jlist}$ to $\beta$.  Suppose first that  
\[
     x_j\text{-degree}(\xlist^{\dlist}) \leq  x_j\text{-degree}(\xlist^{\jlist});
\]
that is, the $j$-th component of the vector $\jlist$ $-$ $\dlist$ is $\geq$ $0$.
By hypothesis, we can translate the arrow $c^{\dlist^{\prime}}_{\jlist^{\prime}}$ so that the head eventually leaves the first orthant.  In fact, we \emph{claim} that we can perform this translation without ever passing through a position $c^{\dlist_1}_{\jlist_1}$ for which 
\[
 x_j\text{-height}(c^{\dlist_1}_{\jlist_1}) < x_j\text{-height}(c^{\dlist^{\prime}}_{\jlist^{\prime}}) = x_j\text{-degree}(\xlist^{\dlist^{\prime}}).
\] 

\emph{Proof of claim:}  We argue by induction on the $x_j$-height of the original arrow $c^{\dlist^{\prime}}_{\jlist^{\prime}}$.  The case of height $0$ is trivial: no arrows can have $x_j$-height $<$ $0$.  So suppose that
\[
    x_j\text{-height}(c^{\dlist^{\prime}}_{\jlist^{\prime}}) = p > 0,
\]
and that the claim holds for all heights $<$ $p$.  Consider a sequence of steps that will translate the arrow's head out of the first orthant.  Let $c^{\dlist_0}_{\jlist_0}$ be the first point on this path (if any) from which translation to an $x_j$-height $<$ $p$ is possible; that is,
\[
     x_j\text{-degree}(\xlist^{\dlist_0}) = p, \text{ and } \xlist^{\dlist_0} / x_j \in I_{\beta}.
\]
Let $s$ $\geq$ $1$ be the largest integer such that 
\[
    \xlist^{\blist} = \xlist^{\dlist_0} / x_j^s \in I_{\beta}.
\]
Then the arrow $c^{\blist}_{\jlist_0}$ is non-standard with the same negative and non-negative components in its vector as the original arrow (the only change in the vector is an increased $j$-th component, which was $\geq$ $0$ to begin with).  In light of the hypothesis on $\beta$, and our induction hypothesis, we can translate the arrow $c^{\blist}_{\jlist_0}$ (of $x_j$-height $q$ $<$ $p$) so that the head exits the first orthant without ever reaching a position of height $<$ $q$.  It is now clear that the arrow $c^{\dlist_0}_{\jlist_0}$ can be translated ``in parallel'' with the arrow $c^{\blist}_{\jlist_0}$, and the height of the former never becomes $<$ $p$.  This shows that the claim holds for any non-standard arrow at $x_j$-height $p$, which completes the proof of the claim.
\bigskip

Returning to the proof of Theorem \ref{thm:truncnsa}, we now see that the lifted arrow $c^{\dlist^{\prime}}_{\jlist^{\prime}}$ can be translated so that the head leaves the first orthant and no position on the path has $x_j$-height less than the original height (which is $\geq$ $h$).  By descending the path to $\beta_t$, we obtain that the original arrow $c^{\dlist}_{\jlist}$ is translation-equivalent to $0$ for the basis set $\beta_t$.

It remains to show that the same conclusion holds when the $j$-th component of the vector of our non-standard arrow $c^{\dlist}_{\jlist}$ is $<$ $0$; that is,
\[
     x_j\text{-degree}(\xlist^{\dlist}) >  x_j\text{-degree}(\xlist^{\jlist}).
\]
   The lifted arrow $c^{\dlist^{\prime}}_{\jlist^{\prime}}$ is by hypothesis translation-equivalent to $0$ on $\beta$.  Note that any step $c^{\dlist_1}_{\jlist_1}$ on the corresponding translation path such that
\[
   x_j\text{-degree}(\xlist^{\jlist_1}) \geq h
\]
descends to $\beta_t$.  Therefore, by descending the initial segment of the path, up to the first position  $c^{\dlist_1}_{\jlist_1}$ for which the $x_j$-degree$(\xlist^{\jlist_1})$ $<$ $h$ (if any), we conclude that $c^{\dlist}_{\jlist}$ is translation-equivalent to $0$ for $\beta_t$, and we are done.
\qedsymbol


\section{Sufficient conditions for a truncation to be smooth} \label{S:truncsmth}

Let $\beta$ be a \emph{smooth} basis set in the variables $x_1, \dots, x_r$, and $\beta_t$ the $x_j$-truncation of $\beta$ at height $h$ (\ref{E:truncdef}).  In the last section (Theorem \ref{thm:truncnsa}) we saw that $\beta_t$ necessarily satisfies condition (a) of Theorem \ref{thm:necsufcond}.  Our main goal in this section is to prove Theorem \ref{thm:truncmain}, which states that $\beta_t$ will also satisfy condition (b) of Theorem \ref{thm:necsufcond}, and therefore be smooth, if we assume an additional hypothesis.  Much of our work involves the construction of a (near-)standard $x_i$-sub-bunch of arrows for $\beta$ that contains the lifting of a standard $x_i$-sub-bunch of arrows for $\beta_t$.  There are two cases: $x_i$ = $x_j$ (the easier case), discussed in Subsection \ref{SS:subbchs1}, and $x_i$ $=$ $x_k$  $\neq$  $x_j$, discussed in Subsection \ref{SS:subbchs2}. 

\subsection{The additional hypothesis} \label{SS:newhyp}

     By the first assertion of Lemma \ref{lem:mgtrunc}, we know that if $m$ is any minimal generator of $I_{\beta}$, then
\[
    \forall x_k \neq x_j,\ x_j\text{-degree}(m) \geq h\ \Rightarrow\ x_k\text{-degree}(m) \leq w^{(t)}_k.
\]
To prove Theorem \ref{thm:truncmain}, we need to control $x_k$-degree$(m)$ when $x_j$-degree($m$) $<$ $h$, using the following 

\begin{hyp} \label{hyp:btsmth}
    For every minimal generator $m$ of $I_{\beta}$, we have that
\[
    \forall x_k \neq x_j,\ x_j\text{-degree}(m) < h\ \text{and } x_k\, |\, m\ \Rightarrow\ x_k\text{-degree}(m) \geq w^{(t)}_k.
\]
\end{hyp}

\subsection{$x_j$-sub-bunches of arrows for $\beta$ and $\beta_t$} \label{SS:subbchs1}

\begin{lem} \label{lem:liftS1}
Let $\beta$ be an arbitrary (not necessarily smooth) basis set such that $\beta_t$, the $x_j$-truncation at height $h$ {\rm (\ref{E:truncdef})}, is defined.  Let $\lis(j)$ be a standard $x_j$-sub-bunch of arrows for $\beta$, constructed as in the proof of Theorem {\rm \ref{thm:S}}.  Then the set  of arrows $c^{\dlist}_{\jlist}$ $\in$ $\lis(j)$ with heads $\xlist^{\jlist}$ divisible by $x_j^h$ (so that $\xlist^{\jlist}/x_j^h$ $\in$ $\beta_t$) is the lifting of a standard $x_j$-sub-bunch of arrows $\lis^{(t)}(j)$ for $\beta_t$.
\end{lem}

\emph{Proof:}
    Consider a minimal $x_j$-standard arrow $c^{\dlist}_{\jlist}$ for $\beta$ such that 
\[
    x_j\text{-degree}(\xlist^{\jlist}) \geq h\ \Rightarrow \ x_j\text{-degree}(\xlist^{\dlist}) > h,
\]      
and suppose that $c^{\dlist}_{\jlist}$ can be advanced.  We can then promote the entire $x_j$-shadow of $c^{\dlist}_{\jlist}$, as described in Subsection \ref{SS:shadprom}; the promotion image is the $x_j$-shadow of an arrow $c^{\blist}_{\jlist_2}$, where $\xlist^{\blist}$ is a minimal generator of $I_{\beta}$, 
\[
    x_j\text{-degree}(\xlist^{\blist}) > x_j\text{-degree}(\xlist^{\jlist_2}) = x_j\text{-degree}(\xlist^{\jlist}) \geq h\,
\]
and
\[
    x_j\text{-degree}(\xlist^{\blist}) < x_j\text{-degree}(\xlist^{\dlist}).
\]
Let $c^{\dlist_t}_{\jlist_t}$ (resp.\ $c^{\blist_t}_{\jlist_{2,t}}$) denote the arrow that $c^{\dlist}_{\jlist}$ (resp.\ $c^{\blist}_{\jlist_2}$) descends to (on the truncation $\beta_t$).  Lemma \ref{lem:mgtrunc} implies that $\xlist^{\dlist_t}$ and $\xlist^{\blist_t}$ are minimal generators of $I_{\beta_t}$; moreover, it is clear that the $x_j$-shadow of $c^{\blist_t}_{\jlist_{2,t}}$ is the promotion image of the shadow of $c^{\dlist_t}_{\jlist_t}$, and that the lifting of the shadow of $c^{\dlist_t}_{\jlist_t}$ (resp.\ $c^{\blist_t}_{\jlist_{2,t}}$) to $\beta$ consists of the arrows in the shadow of $c^{\dlist}_{\jlist}$ (resp.\ $c^{\blist}_{\jlist_2}$) whose heads have $x_j$-degree $\geq$ $h$.  Applying this observation to the iterated shadow promotion operations used to construct $\lis(j)$, as described in Subsection \ref{SS:isp}, one sees that the desired result follows readily.
\qedsymbol

\subsection{$x_k$-standard arrows of $x_j$-height $\geq$ $h$} \label{SS:arrlift}

  Recall from Subsection \ref{SS:isp} that $\lis_0(k,v)$ denotes the set of all (minimal $x_k$-standard) arrows for $\beta$ having tail $x_k^{w_k}$ and offset $v$.  
Suppose that $c^{\dlist}_{\jlist}$ $\in$ $\lis_0(k,v)$ has head $\xlist^{\jlist}$ $\in$ $\beta_t$.  It is then clear that the entire $x_k$-shadow of $c^{\dlist}_{\jlist}$ has this property; therefore, the subset of all arrows in $\lis_0(k,v)$ with heads in $\beta_t$ is equal to the $x_k$-shadow of $c^{\dlist}_{\jlist_v}$, the member of the subset 
whose head has maximal $x_k$-degree
\begin{equation} \label{E:rkv}
    r(k,v)\ =\ x_k\text{-degree}(\xlist^{\jlist_v}).
\end{equation}
If no  arrow in $\lis_0(k,v)$ has head in $\beta_t$, we define $r(k,v)$ $=$ $-1$, and $\xlist^{\jlist_v}$ is undefined. 
 Recalling that the $x_j$-height of an arrow is the $x_j$-degree of its tail, we have the following

\begin{lem} \label{lem:r(k,v)}
    Let $c^{\blist}_{\jlist}$ be an $x_k$-standard arrow for $\beta$ having $x_j$-height $\geq$ $h$ and offset $v$.  Then 
\[
    x_k\text{-degree}(\xlist^{\jlist}) \leq r(k,v).
\]
\end{lem}

\emph{Proof:}
  If not, then, since $c^{\blist}_{\jlist}$ is $x_k$-standard, we have that
\[
  \begin{array}{c}
    x_j\text{-degree}(\xlist^{\jlist}) \geq  x_j\text{-degree}(\xlist^{\blist}) \geq h, \text{and}\\
    x_i\text{-degree}(\xlist^{\jlist}) \geq x_i\text{-degree}(\xlist^{\blist}),\ i \neq j,\ i \neq k.
  \end{array}
\]
Let
\[
    p = x_k\text{-degree}(\xlist^{\blist}),\ m = \xlist^{\blist} / x_k^{p}, \text{ and } \xlist^{\jlist_1} = \xlist^{\jlist}/m \in \beta.
\]
Then we have that $\xlist^{\jlist_1}$ has offset $v$ from the corner monomial $x_k^{w_k}$ $=$ $\xlist^{\dlist}$, and
\[
    x_k\text{-degree}(\xlist^{\jlist_1})\ =\ x_k\text{-degree}(\xlist^{\jlist}).
\]  
We further have that $\xlist^{\jlist_1}$ $\in$ $\beta_t$ because
\[
     x_j^h\, |\, m\ \Rightarrow (x_j^h \cdot \xlist^{\jlist_1})\, |\, (m \cdot \xlist^{\jlist_1}), \text{ and } m \cdot \xlist^{\jlist_1} = \xlist^{\jlist} \in \beta,\ \text{so } x_j^h \cdot \xlist^{\jlist_1} \in \beta.
\]
It then follows from the definition (\ref{E:rkv}) of $r(k,v)$  that 
\[
    x_k\text{-degree}(\xlist^{\jlist})\ = \ x_k\text{-degree}(\xlist^{\jlist_1})\ \leq \ r(k,v),
\]
as desired. 
\qedsymbol

\subsection{Linear independence of lifts of $x_k$-sub-bunches} \label{SS:lilift}

\begin{lem} \label{lem:lifttrans}
    Let $\beta$ and $\beta_t$ be as above, and let $x_k$ $\neq$ $x_j$.  Suppose in addition that the specialization of Hypothesis \emph{\ref{hyp:btsmth}} to $x_k$ holds; that is, for every minimal generator $m$ of $I_{\beta}$, we have that    
\begin{equation} \label{E:restrHyp}
  x_j\text{-degree}(m) < h\ \text{and } x_k\, |\, m\ \Rightarrow\ x_k\text{-degree}(m) \geq w^{(t)}_k.
\end{equation}
Then, given two $x_k$-standard arrows (for $\beta$)  $c^{\dlist_1}_{\jlist_1}$ $\sim$ $c^{\dlist_1^{\prime}}_{\jlist_1^{\prime}}$ of $x_j$-height $\geq$ $h$, there is a translation path from $c^{\dlist_1}_{\jlist_1}$ to $c^{\dlist_1^{\prime}}_{\jlist_1^{\prime}}$ such that every arrow in the path has $x_j$-height $\geq$ $h$.  
\end{lem}

\emph{Proof:}
Consider a translation path from $c^{\dlist_1}_{\jlist_1}$ to $c^{\dlist_1^{\prime}}_{\jlist_1^{\prime}}$  that at some point involves an arrow of $x_j$-height $<$ $h$, and let $c^{\dlist_2}_{\jlist_2}$ be the arrow in the path that immediately precedes the very first such arrow; in particular, since the next step must be in the negative $x_j$-direction to reach $x_j$-height $h-1$, we have that 
\[
     \xlist^{\dlist_2} / x_j \in I_{\beta},\ \text{and } x_j\text{-degree}(\xlist^{\dlist_2}/x_j) = h-1.
\]
There exists a minimal generator $m$ of $I_{\beta}$ such that 
\[
      m\, |\, (\xlist^{\dlist_2}/x_j)\ \Rightarrow\ x_j\text{-degree}(m) < h.
\]
Furthermore, we have that $x_k$ $|$ $m$, since the reasoning leading to the inequality (\ref{eqn:ineq1}) implies that 
\[
    x_k\text{-degree}(m)\ > \ x_k\text{-degree}(\xlist^{\jlist_2})\ \geq \ 0,
\]
otherwise we have the contradiction $m$ $|$ $\xlist^{\jlist_2}$ $\Rightarrow$ $m$ $\in$ $\beta$.
 
The hypothesis (\ref{E:restrHyp}) now implies that 
\[
    w^{(t)}_k \leq x_k\text{-degree}(m) \leq x_k\text{-degree}(\xlist^{\dlist_2}/x_j) = x_k\text{-degree}(\xlist^{\dlist_2}).
\]
Indeed, similar reasoning shows that any arrow in the path having $x_j$-height $<$ $h$ must have $x_k$-height $\geq$ $w^{(t)}_k$.  Since the path terminates in the arrow $c^{\dlist_1^{\prime}}_{\jlist_1^{\prime}}$ of $x_j$-height $\geq$ $h$, the path must eventually reach an arrow $c^{\dlist_3}_{\jlist_3}$ following $c^{\dlist_2}_{\jlist_2}$ such that the $x_j$-height of $c^{\dlist_3}_{\jlist_3}$ equals $h$, all subsequent arrows in the path have $x_j$-height $\geq$ $h$, and the arrow preceding $c^{\dlist_3}_{\jlist_3}$ has $x_j$-height $=$ $h-1$, so that the step to   $c^{\dlist_3}_{\jlist_3}$ is in the increasing $x_j$-direction.  An argument similar to that for $\xlist^{\dlist_2}$ yields that $x_k$-degree$(\xlist^{\dlist_3})$ $\geq$ $w^{(t)}_k$.  It follows that both $\xlist^{\dlist_2}$ and $\xlist^{\dlist_3}$ are divisible by $\xlist^{\blist}$ = $x_j^h \cdot x_k^{w^{(t)}_k}$ $\in$ $I_{\beta}$; therefore, each of the arrows $c^{\dlist_2}_{\jlist_2}$ and $c^{\dlist_3}_{\jlist_3}$ can be translated by degree-reducing steps to an arrow $c^{\blist}_{\jlist_4}$.  The heads of the arrows cannot leave the first orthant during these translations because, the arrows being $x_k$-standard (and of $x_k$-height $\geq$ $w^{(t)}_k$), the exit would have to be across the hyperplane ($x_k$-degree $=$ $0$), but then the original arrow $c^{\dlist_1}_{\jlist_1}$ of $x_k$-height $\leq$ $w^{(t)}_k$ could not have existed.  We now replace the original path segment from $c^{\dlist_2}_{\jlist_2}$ to $c^{\dlist_3}_{\jlist_3}$ by the translation from $c^{\dlist_2}_{\jlist_2}$ to $c^{\blist}_{\jlist_4}$ and the reversal of the translation from $c^{\dlist_3}_{\jlist_3}$ to $c^{\blist}_{\jlist_4}$.  Since the latter translations only involve arrows of $x_j$-height $=$ $h$, we have produced a path from $c^{\dlist_1}_{\jlist_1}$ to $c^{\dlist_1^{\prime}}_{\jlist_1^{\prime}}$ that involves only arrows  of $x_j$-height $\geq$ $h$, as desired.
\qedsymbol

\begin{cor} \label{cor:liftnoadv}
    With the hypotheses of Lemma \emph{\ref{lem:lifttrans}}, suppose given two $x_k$-standard arrows $c^{\dlist_1}_{\jlist_1}$ and $c^{\dlist_2}_{\jlist_2}$ for $\beta_t$, and let $c^{\dlist^{\prime}_1}_{\jlist^{\prime}_1}$, $c^{\dlist^{\prime}_2}_{\jlist^{\prime}_2}$ be the associated liftings to $\beta$.  If $c^{\dlist^{\prime}_1}_{\jlist^{\prime}_1}$ $\sim$ $c^{\dlist^{\prime}_2}_{\jlist^{\prime}_2}$, then $c^{\dlist_1}_{\jlist_1}$ $\sim$ $c^{\dlist_2}_{\jlist_2}$.  Furthermore, if $c^{\dlist^{\prime}_1}_{\jlist^{\prime}_1}$ can be advanced, then $c^{\dlist_1}_{\jlist_1}$ can be advanced.
\end{cor}

\emph{Proof:}
Since $c^{\dlist^{\prime}_1}_{\jlist^{\prime}_1}$ and $c^{\dlist^{\prime}_2}_{\jlist^{\prime}_2}$ are $x_k$-standard and have $x_j$-height $\geq$ $h$, Lemma \ref{lem:lifttrans} implies that there is a translation path from $c^{\dlist^{\prime}_1}_{\jlist^{\prime}_1}$ to $c^{\dlist^{\prime}_2}_{\jlist^{\prime}_2}$ consisting entirely of arrows of $x_j$-height $\geq$ $h$.  But then this path descends to give the translation equivalence of $c^{\dlist_1}_{\jlist_1}$ and $c^{\dlist_2}_{\jlist_2}$, which proves the first assertion.  

Note that the second assertion is trivially true if  
\[
    x_k\text{-degree}(\xlist^{\dlist_1}) > w^{(t)}_k,
\]
for then we can advance $c^{\dlist_1}_{\jlist_1}$ by moving its tail in degree-decreasing steps toward the minimal generator $x_k^{w^{(t)}_k}$.  It therefore remains to show that we can advance $c^{\dlist_1}_{\jlist_1}$ provided that its lifting $c^{\dlist^{\prime}_1}_{\jlist^{\prime}_1}$ can be advanced and 
\[
    x_k\text{-degree}(\xlist^{\dlist^{\prime}_1}) =  x_k\text{-degree}(\xlist^{\dlist_1}) \leq w^{(t)}_k.
\]
Since we can advance $c^{\dlist^{\prime}_1}_{\jlist^{\prime}_1}$, we can translate it to an arrow $c^{\dlist^{\prime}_3}_{\jlist^{\prime}_3}$ such that 
\[
    \xlist^{\dlist^{\prime}_3}/x_k \in I_{\beta} \text{ and } x_k\text{-degree}(\xlist^{\dlist^{\prime}_3}) = x_k\text{-degree}(\xlist^{\dlist^{\prime}_1}).
\]
Then there is a minimal generator $m$  of $I_{\beta}$ such that 
\[
    m\, |\, (\xlist^{\dlist^{\prime}_3}/x_k)\ \text{ and }\ x_k\text{-degree}(m) > 0,
\]
where the inequality follows from (\ref{eqn:ineq1}).

If $x_j$-degree$(m)$ $<$ $h$, then hypothesis (\ref{E:restrHyp}) implies that 
\[
  \begin{array}{rcl}
   x_k\text{-degree}(m)  \geq w^{(t)}_k & \Rightarrow & x_k\text{-degree}(\xlist^{\dlist^{\prime}_3}) > w^{(t)}_k\\
     {} & \Rightarrow & x_k\text{-degree}(\xlist^{\dlist^{\prime}_1}) > w^{(t)}_k,
  \end{array}
\] 
which is a contradiction.
We therefore have that 
\[
  \begin{array}{rcl}
     x_j\text{-degree}(m) \geq h & \Rightarrow & x_j\text{-degree}(\xlist^{\dlist^{\prime}_3}) \geq h;
  \end{array}
\]
  whence, Lemma \ref{lem:lifttrans} ensures that there is a translation path from $c^{\dlist^{\prime}_1}_{\jlist^{\prime}_1}$ to $c^{\dlist^{\prime}_3}_{\jlist^{\prime}_3}$ consisting entirely of arrows of $x_j$-height $\geq$ $h$, but this path then descends to advance $c^{\dlist_1}_{\jlist_1}$, which proves the second assertion.
\qedsymbol

\begin{cor} \label{cor:liftlinind}
    Again with the hypotheses of Lemma \emph{\ref{lem:lifttrans}}, we have that the lifting, to $\beta$, of a standard $x_k$-sub-bunch $\lis^{(t)}(k)$ of arrows for $\beta_t$, has maximal rank \emph{(mod $M^2$)}; furthermore, the lifted arrows cannot be advanced.
\end{cor}

\emph{Proof:}
    It suffices, by Theorem \ref{thm:lindepconds}, to show that for any two distinct arrows $c^{\dlist_1}_{\jlist_1}$, $c^{\dlist_2}_{\jlist_2}$ $\in$ $\lis^{(t)}(k)$, with liftings $c^{\dlist^{\prime}_1}_{\jlist^{\prime}_1}$ and $c^{\dlist^{\prime}_2}_{\jlist^{\prime}_2}$, we have that $c^{\dlist^{\prime}_1}_{\jlist^{\prime}_1}$ $\nsim$ $0$
and $c^{\dlist^{\prime}_1}_{\jlist^{\prime}_1}$ $\nsim$ $c^{\dlist^{\prime}_2}_{\jlist^{\prime}_2}$.  Arguing by contradiction, suppose that  $c^{\dlist^{\prime}_1}_{\jlist^{\prime}_1}$ $\sim$ $c^{\dlist^{\prime}_2}_{\jlist^{\prime}_2}$.  Then Corollary \ref{cor:liftnoadv} implies that $c^{\dlist_1}_{\jlist_1}$ $\sim$ $c^{\dlist_2}_{\jlist_2}$, which contradicts the hypothesis that $\lis^{(t)}(k)$ is an $x_k$-standard sub-bunch.  To prove that $c^{\dlist^{\prime}_1}_{\jlist^{\prime}_1}$ $\nsim$ $0$, it suffices to prove, more generally, that $c^{\dlist^{\prime}_1}_{\jlist^{\prime}_1}$ cannot be advanced.  However, if this were false, then Corollary \ref{cor:liftnoadv} would yield that $c^{\dlist_1}_{\jlist_1}$ can be advanced, which would again contradict the hypothesis that $\lis^{(t)}(k)$ is an $x_k$-standard sub-bunch.
\qedsymbol

\subsection{$x_k$-sub-bunches of arrows for $\beta$ and $\beta_t$, $x_k$ $\neq$ $x_j$} \label{SS:subbchs2}

We are now ready to prove

\begin{lem} \label{lem:liftS2}
     Let $\beta_t$ be the $x_j$-truncation of $\beta$ at height $h$, let $x_k$ $\neq$ $x_j$, and suppose that the specialization \emph{(\ref{E:restrHyp})} of Hypothesis \emph{\ref{hyp:btsmth}} to $x_k$ holds.  Then there exists a near-standard $x_k$-sub-bunch $\lis^{\prime}(k)$ of arrows for $\beta$ that contains the lifting of a standard $x_k$-sub-bunch $\lis^{(t)}(k)$ of arrows for $\beta_t$.
\end{lem}

\emph{Proof:}
     We begin by constructing a standard $x_k$-sub-bunch $\lis(k)$ for $\beta$ as in the proof of Theorem \ref{thm:S}.  With $r(k,v)$ defined as in (\ref{E:rkv}), we let 
\[
    \lis^{\prime}_1(k)\ = \ \{c^{\dlist_1}_{\jlist_1} \in \lis(k) \mid x_k\text{-deg} (\xlist^{\jlist_1}) > r(k,v),\ \text{where } v = x_k\text{-offset}(c^{\dlist_1}_{\jlist_1}) \};   
\]
it is evident that $\lis^{\prime}_1$ is a maximal rank (mod $M^2$) set of $x_k$-standard arrows that cannot be advanced.  It follows from Lemma \ref{lem:r(k,v)} that the $x_j$-height of an arrow $c^{\dlist_1}_{\jlist_1}$ $\in$ $\lis^{\prime}_1(k)$ is $<$ $h$.
Furthermore, the cardinality of this set is
\[
    |\lis^{\prime}_1(k)|\ =\ |\beta| - |\beta_t|\ =\ n - n_t,
\]
because shadow promotion does not change the $x_k$-heights of the heads of the promoted arrows (recall Equation (\ref{E:headDeg})), and the arrows in the original sets $\lis_0(k,v)$ with heads of height $\leq$  $r(k,v)$ are precisely the arrows with heads in $\beta_t$.  

We next construct a standard $x_k$-sub-bunch $\lis^{(t)}(k)$ for $\beta_t$, and denote the lifting of this set to $\beta$ by $\lis^{\prime}_2(k)$; this set consists of 
\[
    |\lis^{\prime}_2(k)|\ =\ |\beta_t|\ =\ n_t 
\]
$x_k$-standard arrows $c^{\dlist_2}_{\jlist_2}$ having $x_j$-height $\geq$ $h$.  Corollary \ref{cor:liftlinind} implies that $\lis^{\prime}_2(k)$ has maximal rank (mod $M^2$), and that its arrows cannot be advanced.  By Lemma \ref{lem:mgtrunc}, we know that the tail $\xlist^{\dlist_2}$ is a minimal generator of $I_{\beta}$ provided that its $x_j$-degree is $>$ $h$; however, if the $x_j$-degree $=$ $h$, we only know that $\xlist^{\dlist_2}$ is an $x_j$-multiple of a minimal generator; therefore, some of the arrows in $\lis^{\prime}_2(k)$ may not be minimal.  

The desired set is
\[
    \lis^{\prime}(k)\ =\ \lis^{\prime}_1(k) \cup \lis^{\prime}_2(k).
\]
By comparing $x_j$-heights of arrows, we see that the two sets in the union do not overlap; therefore,
\[
    |\lis^{\prime}(k)| \ = \ |\lis^{\prime}_1(k)| + |\lis^{\prime}_2(k)| \ = \ (n - n_t) + n_t \ = \ n \ = \ |\beta|.
\]
The arrows in $\lis^{\prime}(k)$ are $x_k$-standard arrows that cannot be advanced, but need not be minimal.  To show that $\lis^{\prime}(k)$ is a near-standard $x_k$-sub-bunch, it remains to show that it has maximal rank (mod $M^2$).  Since $\lis^{\prime}_1(k)$ and $\lis^{\prime}_2(k)$ each have maximal rank (mod $M^2$) and consist of $x_k$-unadvanceable arrows, it suffices (by Theorem \ref{thm:lindepconds}) to prove that no arrow in $\lis^{\prime}_1(k)$ can be translated to an arrow in $\lis^{\prime}_2(k)$.  So let $c^{\dlist_1}_{\jlist_1}$ $\in$ $\lis^{\prime}_1(k)$ have $x_k$-offset $v$.  By definition we have that $x_k$-degree$(\xlist^{\jlist_1})$ $>$ $r(k,v)$.  If this arrow could be translated to an arrow $c^{\dlist_2}_{\jlist_2}$ of $x_j$-height $\geq$ $h$, then Lemma \ref{lem:r(k,v)} would yield that $x_k$-degree$(\xlist^{\jlist_2})$ $\leq$ $r(k,v)$, implying that $c^{\dlist_1}_{\jlist_1}$ can be advanced, which is a contradiction.  Therefore, $\lis^{\prime}(k)$ is a near-standard $x_k$-sub-bunch that contains the lifting of a standard $x_k$-sub-bunch for $\beta_t$, as desired.  \qedsymbol

\subsection{Main theorem on truncations} \label{SS:truncmain}

\begin{thm} \label{thm:truncmain}
    Let $\beta$ be a smooth basis set, and $\beta_t$ the $x_j$-truncation of $\beta$ at height $h$.  If  in addition Hypothesis \emph{\ref{hyp:btsmth}} holds, then $\beta_t$ is a smooth basis set.
\end{thm}

\emph{Proof:}
     To show that $\beta_t$ is a smooth basis set, it suffices to prove that it satisfies conditions (a) and (b) of Theorem \ref{thm:necsufcond}.  Theorem \ref{thm:truncnsa} ensures that $\beta_t$ satisfies condition (a) --- that every non-\sa\ arrow is translation-equivalent to 0 --- because $\beta$, smooth by hypothesis, has this property; therefore, it remains to show that $\beta_t$ satisfies condition (b).      

     Let $\lis(j)$ be a standard $x_j$-sub-bunch of arrows for $\beta$, and $\lis^{(t)}(j)$ the associated $x_j$-sub-bunch for $\beta_t$ whose lifting to $\beta$ lies in $\lis(j)$, as in Lemma \ref{lem:liftS1}. For each variable $x_k$ $\neq$ $x_j$, let $\lis^{\prime}(k)$ be a near-standard $x_k$-sub-bunch of arrows for $\beta$ constructed as in Lemma \ref{lem:liftS2}, and $\lis^{(t)}(k)$ the associated $x_k$-sub-bunch for $\beta_t$ whose lifting to $\beta$ lies in $\lis^{\prime}(k)$.  Taking unions, we obtain a near-standard bunch $\lis^{\prime}$ for $\beta$ that contains the lifting of a standard bunch $\lis^{(t)}$ for $\beta_t$.  To complete the proof, it suffices to show that if $c^{\dlist}_{\jlist}$ $\nsim$ $0$ is a \sa\ arrow for $\beta_t$, then $c^{\dlist}_{\jlist}$ $\sim$ $c^{\dlist_1}_{\jlist_1}$ $\in$ $\lis^{(t)}$.

     By assertion (c) of Lemma \ref{lem:shadadv}, we may assume that $c^{\dlist}_{\jlist}$ is a minimal $x_i$-standard arrow for $\beta_t$ that cannot be advanced.  Consider the lifting $c^{\dlist^{\prime}}_{\jlist^{\prime}}$ of $c^{\dlist}_{\jlist}$ to $\beta$.  We have that $c^{\dlist^{\prime}}_{\jlist^{\prime}}$ cannot be advanced: Indeed, if $x_i$ $=$ $x_j$, then one checks easily that a translation path advancing $c^{\dlist^{\prime}}_{\jlist^{\prime}}$ would descend to a translation path advancing $c^{\dlist}_{\jlist}$, a contradiction.  The same contradiction arises in case $x_i$ $=$ $x_k$ $\neq$ $x_j$ by Corollary \ref{cor:liftnoadv}. Since $\beta$ is assumed smooth, we know that 
\[ 
    c^{\dlist^{\prime}}_{\jlist^{\prime}} \sim c^{\dlist^{\prime}_1}_{\jlist^{\prime}_1} \in \lis^{\prime}.
\] 
Furthermore, since neither of the arrows $c^{\dlist^{\prime}}_{\jlist^{\prime}}$, $c^{\dlist^{\prime}_1}_{\jlist^{\prime}_1}$ can be advanced, we see that these arrows must have the same $x_i$-height.  If $x_i$ $=$ $x_j$, we obtain that $c^{\dlist^{\prime}_1}_{\jlist^{\prime}_1}$ lies in the lifting of $\lis^{(t)}(j)$, and the translation path from $c^{\dlist^{\prime}}_{\jlist^{\prime}}$ to $c^{\dlist^{\prime}_1}_{\jlist^{\prime}_1}$ descends to yield  
\[
    c^{\dlist}_{\jlist} \sim c^{\dlist_1}_{\jlist_1} \in \lis^{(t)}(j).
\]
If $x_i$ $=$ $x_k$ $\neq$ $x_j$, we obtain a similar conclusion as follows: recalling that $\lis^{\prime}_2(k)$ $\subseteq$ $\lis^{\prime}(k)$ denotes the lifting of the $x_k$-standard sub-bunch $\lis^{(t)}(k)$ (in the notation of Theorem \ref{lem:liftS2}), we \emph{claim} that  
\begin{equation} \label{E:inlis2}
     c^{\dlist^{\prime}_1}_{\jlist^{\prime}_1} \in \lis^{\prime}_2(k) \subseteq \lis^{\prime}(k).
\end{equation}

\emph{Proof of claim:}
     The lifted arrow $c^{\dlist^{\prime}}_{\jlist^{\prime}}$ has $x_j$-height $\geq$ $h$; therefore, by Lemma \ref{lem:r(k,v)}, we have that
\[
    x_k\text{-degree}(\xlist^{\jlist^{\prime}}) \leq r(k,v),
\]
where $v$ $=$ $x_k$-offset$(c^{\dlist^{\prime}}_{\jlist^{\prime}})$.  Since the translation-equivalent arrow $c^{\dlist^{\prime}_1}_{\jlist^{\prime}_1}$ $\in$ $\lis^{\prime}(k)$ has the same $x_k$-height and $x_k$-offset, it satisfies
\[
    x_k\text{-degree}(\xlist^{\jlist^{\prime}_1}) =  x_k\text{-degree}(\xlist^{\jlist^{\prime}}) \leq r(k,v).
\]
The claim follows immediately, because $\lis^{\prime}(k)$ $=$ $\lis^{\prime}_1(k)$ $\cup$ $\lis^{\prime}_2(k)$, and by definition the arrows $c^{\dlist_*}_{\jlist_*}$ $\in$ $\lis^{\prime}_1(k)$ of $x_k$-offset $v$ satisfy
\[
    x_k\text{-degree}(\xlist^{\jlist_*}) > r(k,v).
\]
\medskip

The claim (\ref{E:inlis2}) implies that the arrow  $c^{\dlist^{\prime}_1}_{\jlist^{\prime}_1}$ is the lifting of an arrow $c^{\dlist_1}_{\jlist_1}$ $\in$ $\lis^{(t)}(k)$.  Corollary \ref{cor:liftnoadv} now yields that $c^{\dlist}_{\jlist}$ $\sim$ $c^{\dlist_1}_{\jlist_1}$, and the proof is complete.  
\qedsymbol


\section{Addition of boxes to basis sets} \label{S:boxes}

     In this section we explore another way to construct a smooth basis set from a given smooth basis set, by ``adding a box'' in an appropriate way.  The undoing of this operation (that is, the removal of the added box) is accomplished by a truncation.

\subsection{Definition of box addition} \label{SS:boxadddefn}

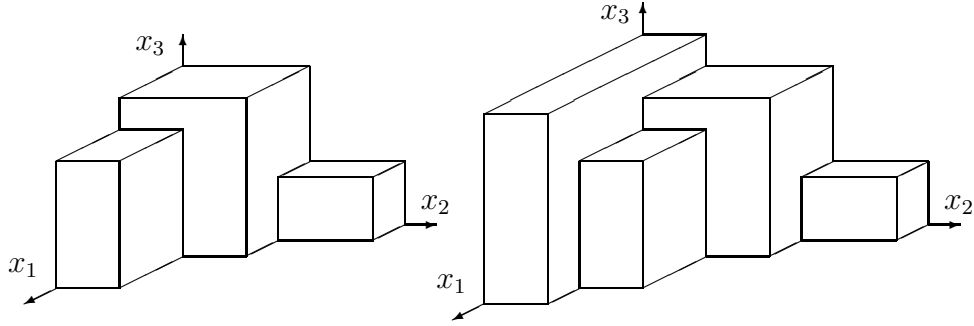
\begin{figure} 
\begin{picture}(374,108)
        \put(30,6){\line(1,0){24}}
        \put(30,6){\line(0,1){48}}
        \put(54,6){\line(0,1){48}}
        \put(30,54){\line(1,0){24}}

        \put(30,54){\line(2,1){24}}
        \put(54,54){\line(2,1){24}}
        \put(54,6){\line(2,1){24}}
        \put(78,18){\line(0,1){48}}

        \put(54,66){\line(1,0){24}}
        \put(54,66){\line(0,1){12}}
        \put(54,78){\line(1,0){48}}
        \put(54,78){\line(2,1){24}}
        \put(78,90){\line(1,0){48}}
        \put(102,78){\line(2,1){24}}
        \put(102,78){\line(0,-1){60}}
        \put(78,18){\line(1,0){24}}

        \put(102,18){\line(2,1){12}}
        \put(114,24){\line(1,0){36}}
        \put(114,24){\line(0,1){24}}
        \put(114,48){\line(2,1){12}}
        \put(114,48){\line(1,0){36}}

        \put(126,54){\line(0,1){36}}
        \put(126,54){\line(1,0){36}}

        \put(150,24){\line(2,1){12}}
        \put(150,24){\line(0,1){24}}

        \put(150,48){\line(0,-1){24}}
        \put(150,48){\line(2,1){12}}
        \put(162,54){\line(0,-1){24}}

        \put(216,0){\line(2,1){12}}
        \put(192,0){\line(1,0){24}}
        \put(192,0){\line(0,1){72}}
        \put(216,0){\line(0,1){72}}
        \put(192,72){\line(1,0){24}}
        \put(192,72){\line(2,1){60}}
        \put(216,72){\line(2,1){60}}
        \put(252,102){\line(1,0){24}}
        \put(276,102){\line(0,-1){12}}

        \put(228,6){\line(1,0){24}}
        \put(228,6){\line(0,1){48}}
        \put(252,6){\line(0,1){48}}
        \put(228,54){\line(1,0){24}}

        \put(228,54){\line(2,1){24}}
        \put(252,54){\line(2,1){24}}
        \put(252,6){\line(2,1){24}}
        \put(276,18){\line(0,1){48}}

        \put(252,66){\line(1,0){24}}
        \put(252,66){\line(0,1){12}}
        \put(252,78){\line(1,0){48}}
        \put(252,78){\line(2,1){24}}
        \put(276,90){\line(1,0){48}}
        \put(300,78){\line(2,1){24}}
        \put(300,78){\line(0,-1){60}}
        \put(276,18){\line(1,0){24}}

        \put(300,18){\line(2,1){12}}
        \put(312,24){\line(1,0){36}}
        \put(312,24){\line(0,1){24}}
        \put(312,48){\line(2,1){12}}
        \put(312,48){\line(1,0){36}}

        \put(324,54){\line(0,1){36}}
        \put(324,54){\line(1,0){36}}

        \put(348,24){\line(2,1){12}}
        \put(348,24){\line(0,1){24}}

        \put(348,48){\line(0,-1){24}}
        \put(348,48){\line(2,1){12}}
        \put(360,54){\line(0,-1){24}}

        \put(30,6){\vector(-2,-1){12}}
        \put(12,12){$x_1$}
        \put(162,30){\vector(1,0){12}}
        \put(168,36){$x_2$}
        \put(78,90){\vector(0,1){12}}
        \put(60,96){$x_3$}

        \put(192,0){\vector(-2,-1){12}}
        \put(174,6){$x_1$}
        \put(252,102){\vector(0,1){12}}
        \put(236,108){$x_3$}
        \put(360,30){\vector(1,0){12}}
        \put(366,36){$x_2$} 
\end{picture}
\caption{The basis set on the right is obtained by adding a box in the $x_2$-direction to the basis set on the left.} \label{Fig:boxaddn}
\end{figure}

Suppose that $\beta$ is an arbitrary basis set of monomials in the variables $x_1,\dots,x_r$; recall that $w_i$ denotes the $x_i$-width of $\beta$ for $1 \leq i \leq r$.  Choose one of the variables, say $x_j$, and an integer $h$ $\geq$ $1$, and form the set
\[
    x_j^h \cdot \beta = \{ x_j^h \cdot m \mid m \in \beta \},
\]
which can be viewed geometrically as the translation of $\beta$ a distance of $h$ steps in the positive $x_j$-direction.  Then, for each variable $x_k$ $\neq$ $x_j$, choose an integer $w_k^{\prime}$ $\geq$ $w_k$, and form the box (\ref{E:rbox})
\begin{equation} \label{E:thebox}
    \rbox = \rbox(w_1^{\prime}, w_2^{\prime},\dots, w_{j-1}^{\prime}, h, w_{j+1}^{\prime}, \dots, w_r^{\prime}).
\end{equation}
We then set
\begin{equation} \label{E:boxaddndef}
    \beta^{\prime} \ = \ (x_j^h \cdot \beta)\ \cup\ \rbox,
\end{equation}
and say that $\beta^{\prime}$  is obtained by \textbf{adding the box} $\rbox$ \textbf{in the} $x_j$\textbf{-direction to} $\beta$ (see Figure \ref{Fig:boxaddn}).

\begin{lem}  \label{lem:baibs} 
     Let $\beta^{\prime}$ be the set of monomials obtained by adding the box $\rbox$ \emph{(\ref{E:thebox})} to the basis set $\beta$ in the $x_j$-direction, as in \emph{(\ref{E:boxaddndef})}.  Then
\begin{enumerate}
    \item[(a)] $\beta^{\prime}$ is a basis set.
    \item[(b)] $\beta$ is the $x_j$-truncation of $\beta^{\prime}$ at height $h$.
\end{enumerate}
\end{lem}

\emph{Proof:}
      To prove that $\beta^{\prime}$ is a basis set, we must show that if $m_1$ and $m_2$  are monomials such that $m_1$ $\in$ $\beta^{\prime}$ and $m_2$ $|$ $m_1$, then $m_2$ $\in$ $\beta^{\prime}$.
If $m_1$ $\in$ $\rbox$, then it is clear that $m_2$ $\in$ $\rbox$ $\subseteq$ $\beta^{\prime}$.  If $m_1$ $\notin$ $\rbox$, then $m_1 / x_j^h$ $\in$ $\beta$.  Since $m_2$ $|$ $m_1$, it is clear that 
\[
    x_k \neq x_j \Rightarrow x_k\text{-degree}(m_2) \leq x_k\text{-degree}(m_1/x_j^h) < w_k \leq w_k^{\prime};
\]
therefore, 
\[
    x_j\text{-degree}(m_2) < h \Rightarrow m_2 \in \rbox \subseteq \beta^{\prime},
\]
 and 
\[
   \begin{array}{rcl}
    x_j\text{-deg}(m_2) \geq h & \Rightarrow & (m_2 / x_j^h)\, |\, (m_1 / x_j^h)\\
     {} &  \Rightarrow & m_2 /x_j^h \in \beta \\
     {} &  \Rightarrow  & m_2 \in x_j^h \cdot \beta \subseteq \beta^{\prime}.
   \end{array}
\]
This completes the proof that $\beta^{\prime}$ is a basis set, and the second statement follows easily from the definitions.     
\qedsymbol

\subsection{Minimal generators of $I_{\beta^{\prime}}$} \label{SS:boxaddmgs}

Let $\beta$, $\beta^{\prime}$, $\rbox$, \etc, be as above.  Since $\beta$ $=$ $\beta^{\prime}_t$ is the $x_j$-truncation of $\beta^{\prime}$ at height $h$, we have that
\[
    w_j^{\prime}\ =\ w_j + h
\]
is the $x_j$-width of $\beta^{\prime}$, by assertion (d) of Lemma \ref{lem:mgtrunc}.
We then have the following

\begin{lem} \label{lem:boxaddncms} 
    For each $i$, $1$ $\leq$ $i$ $\leq$ $r$, the $x_i$-width of $\beta^{\prime}$ is $w_i^{\prime}$; that is, $x_i^{w_i^{\prime}}$ is a minimal generator of $I_{\beta^{\prime}}$.  Furthermore, for all minimal generators $m$ of $I_{\beta^{\prime}}$ we have that 
\[
    x_j\text{-degree}(m) < h \ \Rightarrow\  m = x_k^{w_k^{\prime}} \text{ for some } x_k \neq x_j;
\]
in particular, Hypothesis \emph{\ref{hyp:btsmth}} holds for $\beta^{\prime}$ and its truncation $\beta$.
\end{lem}

\emph{Proof:}
For all variables $x_i$ $=$ $x_k$ $\neq$ $x_j$, we have that 
\[
    x_k^{w_k^{\prime}} \in I_{\beta^{\prime}} \subseteq I_{\rbox},
\]
and $x_k^{w_k^{\prime}}$ is a minimal generator of $I_{\rbox}$, by Lemma \ref{lem:boxmingen}; therefore, $x_k^{w_k^{\prime}}$ is a minimal generator of $I_{\beta^{\prime}}$.  For $x_i$ $=$ $x_j$, 
we have already observed that $w_j^{\prime}$ is the $x_j$-width of $\beta^{\prime}$, so the first assertion holds.  Suppose now that $m$ is a minimal generator of $I_{\beta^{\prime}}$ such that $x_j$-degree$(m)$ $<$ $h$.  A moment's reflection shows that in fact $m$ must be a minimal generator of the monomial ideal $I_{\rbox}$; the second assertion then follows at once from Lemma \ref{lem:boxmingen}. 
\qedsymbol
\bigskip

In addition, recall that Lemma \ref{lem:mgtrunc} further describes the relationship between the minimal generators of $I_{\beta^{\prime}}$ and the minimal generators of $I_{\beta}$, since $\beta$ is the $x_j$-truncation of $\beta^{\prime}$ at height $h$.

\subsection{Main theorem on box additions} \label{SS:boxaddmt}

     Box addition is a convenient tool for building smooth basis sets, as the following result suggests:

\begin{thm} \label{thm:boxaddnmain}
    Let $\beta$ be a basis set, and let $\beta^{\prime}$ be obtained by adding the box $\rbox$ \emph{(\ref{E:thebox})} in the $x_j$-direction.  Then: 
\[
    \beta \text{ is smooth } \Leftrightarrow \beta^{\prime} \text{ is smooth.}
\]
\end{thm}

\emph{Proof:} 
    $(\Leftarrow)$: Since $\beta^{\prime}$ is smooth and Hypothesis \ref{hyp:btsmth} holds, by Lemma \ref{lem:boxaddncms}, Theorem \ref{thm:truncmain} implies that the truncation $\beta^{\prime}_t$ $=$ $\beta$ is smooth.
\bigskip

    $(\Rightarrow)$: Given that $\beta$ is a smooth basis set, we must prove that $\beta^{\prime}$ is a smooth basis set.  To do this we will use Theorem \ref{thm:necsufcond}: we must show (for $\beta^{\prime}$) that every non-standard arrow  is translation-equivalent to $0$, and that there exists a near-standard bunch of arrows $\lis^{\prime}$ such that if $c^{\dlist}_{\jlist}$ $\nsim$ $0$ is an $x_i$-\sa\ arrow, then there is an $x_i$-\sa\ arrow $c^{\dlist^{\prime}}_{\jlist^{\prime}}$ $\in$  $\lis^{\prime}$ such that $c^{\dlist}_{\jlist}$ $\sim$ $c^{\dlist^{\prime}}_{\jlist^{\prime}}$.  

First, let $c^{\dlist}_{\jlist}$ be a non-\sa\ arrow, which we can assume is minimal by Lemma \ref{lem:msa}.  Since the arrow is non-\sa\, we know that the vector $\jlist$ $-$ $\dlist$ has negative components in at least two variable directions, say $x_{i_1}$ and $x_{i_2}$.  Note that the tail of the arrow must be a minimal generator of $I_{\beta^{\prime}}$ that has $x_j$-degree $\geq$ $h$, since the minimal generators of $\beta^{\prime}$ with $x_j$-degree $<$ $h$ are the corner monomials  $x_k^{w_k^{\prime}}$ (Lemma \ref{lem:boxaddncms}), and a non-standard arrow cannot have a corner monomial as its tail.

Suppose first that the head of the arrow $\xlist^{\jlist}$ $\in$ $\rbox$.  By translating the arrow in the increasing $x_{i_2}$-direction, one eventually reaches an arrow $c^{\dlist_1}_{\jlist_1}$ $\sim$ $c^{\dlist}_{\jlist}$ such that 
\[
    x_{i_2}\text{-degree}(\xlist^{\jlist_1}) = w_{i_2}^{\prime}-1 \Rightarrow x_{i_2}\text{-degree}(\xlist^{\dlist_1}) \geq w_{i_2}^{\prime};
\]
therefore, $c^{\dlist_1}_{\jlist_1}$ can be translated in the direction of decreasing $x_{i_1}$-degree until the head exits the first orthant; whence, $c^{\dlist}_{\jlist}$ $\sim$ $0$.  

If $\xlist^{\jlist}$ $\notin$ $\rbox$, then $\xlist^{\jlist}$ $\in$ $x_j^h \cdot \beta$, and the arrow $c^{\dlist}_{\jlist}$ descends to $c^{\dlist_0}_{\jlist_0}$.  Since $\beta$ is assumed smooth, we have that $c^{\dlist_0}_{\jlist_0}$ $\sim$ $0$.  If the associated translation path causes the head to cross the hyperplane ($x_k$-degree $=$ $0$) for any $x_k$ $\neq$ $x_j$, then we can lift the translation path to $\beta^{\prime}$ to obtain that $c^{\dlist}_{\jlist}$ $\sim$ $0$.  If the head crosses the hyperplane ($x_j$-degree $=$ $0$), then the lifted path shows that $c^{\dlist}_{\jlist}$ $\sim$ $c^{\dlist_2}_{\jlist_2}$ with $\xlist^{\jlist_2}$ $\in$ $\rbox$, and $c^{\dlist_2}_{\jlist_2}$ $\sim$ $0$ as before.

     Now let $\lis^{\prime}$ be the near-\sa\ bunch of arrows for $\beta^{\prime}$ that was constructed in the proof of Theorem \ref{thm:truncmain} (based on Lemmas \ref{lem:liftS1} and \ref{lem:liftS2}); recall that $\lis^{\prime}$ contains the lifting of a standard sub-bunch $\lis^{(t)}$ for $\beta$ $=$ $\beta^{\prime}_t$.  Suppose that $c^{\dlist}_{\jlist}$ $\nsim$ $0$ is an $x_i$-\sa\ arrow; we can assume that this arrow is minimal and cannot be advanced by assertion (c) of Lemma \ref{lem:shadadv}.  We must show that $c^{\dlist}_{\jlist}$ is translation-equivalent to an arrow in $\lis^{\prime}(i)$.

Suppose first that $x_i$ = $x_j$, so that $\lis^{\prime}(i)$ $=$ $\lis(j)$ is a standard $x_j$-sub-bunch.  If $\xlist^{\jlist}$ $\notin$ $\rbox$, then the arrow $c^{\dlist}_{\jlist}$ is the lifting of a minimal \sa\ arrow $c^{\dlist_0}_{\jlist_0}$ for $\beta$.  Since $\beta$ is smooth, we can find a minimal \sa\ arrow $c^{\dlist^{\prime}_0}_{\jlist^{\prime}_0}$ $\in$ $\lis^{(t)}(j)$ such that $c^{\dlist_0}_{\jlist_0}$ $\sim$ $c^{\dlist^{\prime}_0}_{\jlist^{\prime}_0}$.  Lifting the translation path, we find that 
$c^{\dlist}_{\jlist}$ $\sim$ $c^{\dlist^{\prime}}_{\jlist^{\prime}}$ $\in$ $\lis^{\prime}(j)$.  

If $\xlist^{\jlist}$ $\in$ $\rbox$, then let $c^{\dlist^{\prime}}_{\jlist^{\prime}}$ $\in$ $\lis^{\prime}(j)$  be the unique arrow  such that
\[
    x_j\text{-degree}(\xlist^{\jlist^{\prime}}) = x_j\text{-degree}(\xlist^{\jlist})\ \ \text{and \ $x_j$-offset}(c^{\dlist^{\prime}}_{\jlist^{\prime}}) = x_j\text{-offset}(c^{\dlist}_{\jlist}),
\]   
the existence of which is ensured by Corollary \ref{cor:degofffact}.  Suppose that the lengths of these two arrows differ; in other words, suppose that
\[
    x_j\text{-degree}(\xlist^{\dlist}) \neq x_j\text{-degree}(\xlist^{\dlist^{\prime}}).
\]
Recalling Lemma \ref{lem:boxlcm}, we let 
\[
    \xlist^{\jlist_1} = \lcm(\xlist^{\jlist}, \xlist^{\jlist^{\prime}}) \in \rbox \subseteq \beta^{\prime}.
\]
It is clear that we can translate both $c^{\dlist}_{\jlist}$ and $c^{\dlist^{\prime}}_{\jlist^{\prime}}$ by degree-increasing steps (excluding the $x_j$-direction) to arrows $c^{\dlist_1}_{\jlist_1}$ and $c^{\dlist^{\prime}_1}_{\jlist_1}$, respectively.  The tails $\xlist^{\dlist_1}$ and $\xlist^{\dlist_1^{\prime}}$ differ only in $x_j$-degree; it follows that the arrow corresponding to the tail of larger $x_j$-degree can be advanced, which is a contradiction, since neither $c^{\dlist}_{\jlist}$ nor $c^{\dlist^{\prime}}_{\jlist^{\prime}}$ can be advanced.  We therefore have that 
\[
    c^{\dlist_1}_{\jlist_1} = c^{\dlist^{\prime}_1}_{\jlist_1}\ \Rightarrow\ c^{\dlist}_{\jlist} \sim c^{\dlist^{\prime}}_{\jlist^{\prime}} \in \lis^{\prime}(j),
\]
as desired.

Finally, we have to consider the case in which $c^{\dlist}_{\jlist}$ is a minimal $x_k$-\sa\ arrow that cannot be advanced, where $x_k$ $\neq$ $x_j$.  Let $v$ denote the $x_k$-offset of $c^{\dlist}_{\jlist}$.  From Lemma \ref{lem:boxaddncms}, we know that the tail $\xlist^{\dlist}$ is either equal to the corner monomial $x_k^{w_k^{\prime}}$ or else $x_j$-degree$(\xlist^{\dlist})$ $\geq$ $h$.
  
In case $\xlist^{\dlist}$ $=$ $x_k^{w_k^{\prime}}$, if the head 
\[
    \xlist^{\jlist}\ \notin\ \beta = \beta^{\prime}_t,
\]
then
\[
    x_k\text{-degree}(\xlist^{\jlist}) > r(k,v)\ \Rightarrow\ c^{\dlist}_{\jlist} \in \lis^{\prime}_1(k) \subseteq \lis^{\prime}(k),
\]
in the notation of Lemma \ref{lem:liftS2}.  If the head  
\[
    \xlist^{\jlist}\ \in\ \beta = \beta^{\prime}_t,
\]
then, since $c^{\dlist}_{\jlist}$ cannot be advanced, but \emph{can} be translated $h$ steps in the direction of increasing $x_j$-degree to $c^{\dlist^{\prime}}_{\jlist^{\prime}}$, we must have that $w_k^{\prime}$ $=$ $w_k$.  This in turn implies that $c^{\dlist}_{\jlist}$, considered as an arrow for the truncation $\beta$, lies in the sub-bunch $\lis^{(t)}(k)$; therefore,
\[
    c^{\dlist}_{\jlist} \sim c^{\dlist^{\prime}}_{\jlist^{\prime}} \in \lis^{\prime}_2(k) \subseteq \lis^{\prime}(k).
\]

In case $x_j$-degree$(\xlist^{\dlist})$ $\geq$ $h$, we can descend the arrow to $c^{\dlist_0}_{\jlist_0}$; then, because $\beta$ is assumed smooth, we know that $c^{\dlist_0}_{\jlist_0}$ $\sim$ $c^{\dlist_0^{\prime}}_{\jlist_0^{\prime}}$ $\in$ $\lis^{(t)}(k)$, and the translation path lifts to yield
$c^{\dlist}_{\jlist}$ $\sim$ $c^{\dlist^{\prime}}_{\jlist^{\prime}}$ $\in$ $\lis^{\prime}(k)$.  This completes the proof of the theorem.
\qedsymbol

\subsection{Compound boxes} \label{SS:cpboxes}

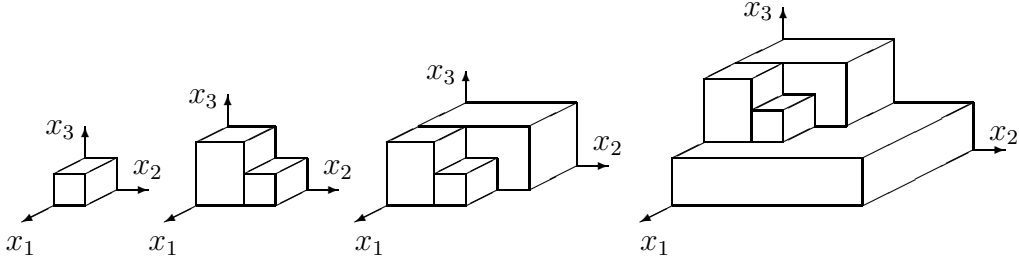
\begin{figure}
\begin{picture}(374,84)

        \put(6,29){\line(1,0){12}}
        \put(6,29){\line(2,1){12}}
        \put(6,29){\line(0,-1){12}}
        \put(6,17){\line(1,0){12}}
        \put(18,35){\line(1,0){12}}
        \put(30,35){\line(-2,-1){12}}
        \put(30,23){\line(-2,-1){12}}
        \put(18,17){\line(0,1){12}}
        \put(30,35){\line(0,-1){12}}
        \put(6,17){\vector(-2,-1){12}}
        \put(-12,0){$x_1$}
        \put(30,23){\vector(1,0){12}}
        \put(36,29){$x_2$}
        \put(18,35){\vector(0,1){12}}
        \put(3,44){$x_3$}

        \put(60,41){\line(1,0){18}} 
        \put(60,41){\line(2,1){12}}
        \put(72,47){\line(1,0){18}}
        \put(90,47){\line(0,-1){12}}
        \put(78,41){\line(2,1){12}} 
        \put(78,41){\line(0,-1){12}} 
        \put(60,41){\line(0,-1){24}} 
        \put(60,17){\line(1,0){30}}  
        \put(78,29){\line(1,0){12}}
        \put(78,29){\line(2,1){12}}
        \put(78,29){\line(0,-1){12}}
        \put(78,17){\line(1,0){12}}
        \put(90,35){\line(1,0){12}}
        \put(102,35){\line(-2,-1){12}}
        \put(102,23){\line(-2,-1){12}}
        \put(90,17){\line(0,1){12}}
        \put(102,35){\line(0,-1){12}}
        \put(60,17){\vector(-2,-1){12}}
        \put(42,0){$x_1$}
        \put(102,23){\vector(1,0){12}}
        \put(108,29){$x_2$}
        \put(72,47){\vector(0,1){12}}
        \put(57,56){$x_3$}

        \put(186,47){\line(2,1){18}}
        \put(162,47){\line(1,0){12}}
        \put(186,47){\line(0,-1){24}}
        \put(186,23){\line(2,1){18}}
        \put(186,23){\line(-1,0){12}}
        \put(144,47){\line(2,1){18}}
        \put(162,56){\line(1,0){42}}
        \put(204,56){\line(0,-1){24}}
          \put(132,41){\line(1,0){18}} 
          \put(132,41){\line(2,1){12}}
          \put(144,47){\line(1,0){42}} 
          \put(162,47){\line(0,-1){12}}
          \put(150,41){\line(2,1){12}} 
          \put(150,41){\line(0,-1){12}} 
          \put(132,41){\line(0,-1){24}} 
          \put(132,17){\line(1,0){30}}  
          \put(150,29){\line(1,0){12}}
          \put(150,29){\line(2,1){12}}
          \put(150,29){\line(0,-1){12}}
          \put(150,17){\line(1,0){12}}
          \put(162,35){\line(1,0){12}}
          \put(174,35){\line(-2,-1){12}}
          \put(174,23){\line(-2,-1){12}}
          \put(162,17){\line(0,1){12}}
          \put(174,35){\line(0,-1){12}}
        \put(132,17){\vector(-2,-1){12}}
        \put(120,0){$x_1$}
        \put(204,32){\vector(1,0){12}}
        \put(210,38){$x_2$}
        \put(162,56){\vector(0,1){12}}
        \put(147,65){$x_3$}

           \put(252,41){\line(-2,-1){12}}
           \put(240,35){\line(0,-1){18}} 
           \put(240,17){\line(1,0){72}}
           \put(240,35){\line(1,0){72}}
           \put(312,35){\line(0,-1){18}}
           \put(312,35){\line(2,1){42}}
           \put(312,17){\line(2,1){42}}
           \put(354,56){\line(0,-1){18}}
           \put(354,56){\line(-1,0){30}}
           \put(306,71){\line(2,1){18}}
           \put(282,71){\line(1,0){12}}
           \put(306,71){\line(0,-1){24}}
           \put(306,47){\line(2,1){18}}
           \put(306,47){\line(-1,0){12}}
           \put(264,71){\line(2,1){18}}
           \put(282,80){\line(1,0){42}}
           \put(324,80){\line(0,-1){24}}
           \put(252,65){\line(1,0){18}} 
           \put(252,65){\line(2,1){12}}
           \put(264,71){\line(1,0){42}} 
           \put(282,71){\line(0,-1){12}}
           \put(270,65){\line(2,1){12}} 
           \put(270,65){\line(0,-1){12}} 
           \put(252,65){\line(0,-1){24}} 
           \put(252,41){\line(1,0){30}}  
           \put(270,53){\line(1,0){12}}
           \put(270,53){\line(2,1){12}}
           \put(270,53){\line(0,-1){12}}
           \put(270,41){\line(1,0){12}}
           \put(282,59){\line(1,0){12}}
           \put(294,59){\line(-2,-1){12}}
           \put(294,47){\line(-2,-1){12}}
           \put(282,41){\line(0,1){12}}
           \put(294,59){\line(0,-1){12}}
        \put(240,17){\vector(-2,-1){12}}
        \put(228,0){$x_1$}
        \put(354,38){\vector(1,0){12}}
        \put(360,42){$x_2$}
        \put(282,80){\vector(0,1){12}}
        \put(267,89){$x_3$}
\end{picture}
\caption{A sequence of compound boxes, beginning on the left with a box.} \label{fig:cpbox}
\end{figure}

     By a \textbf{compound box}, we mean a basis set $\beta$ that is constructed by starting with a box, and then performing a finite sequence of box additions in various variable directions; Figure \ref{fig:cpbox} illustrates the idea.  Note that a box $\rbox$ is a compound box, since it can be generated by starting with itself and performing a sequence of box additions of length $0$.

Since the starting point (a box) is smooth, by Proposition \ref{prop:boxsmth}, and adding a box to a smooth basis set yields a smooth basis set, by Theorem \ref{thm:boxaddnmain}, we obtain the following corollary by induction:

\begin{cor} \label{cor:cpboxsmth}
    If $\beta$ is a compound box, then $\beta$ is a smooth basis set. \qedsymbol
\end{cor}      

\subsection{Example: Basis sets in two variables} \label{SS:2varcase}

\begin{figure}
\begin{picture}(374,72)
        \put(148,0){\line(1,0){72}}
        \put(148,0){\line(0,1){60}}
        \put(148,60){\line(1,0){24}}
        \put(172,60){\line(0,-1){12}}
        \multiput(172,47)(0,-4){12}{\line(0,-1){2}}
        \put(172,48){\line(1,0){18}}
        \put(190,48){\line(0,-1){24}}
        \put(190,24){\line(1,0){30}}
        \multiput(190,23)(0,-4){6}{\line(0,-1){2}}
        \put(220,0){\line(0,1){24}}
        \put(220,0){\vector(1,0){12}}
        \put(226,6){$x_1$}
        \put(148,60){\vector(0,1){12}}
        \put(130,66){$x_2$}

\end{picture}
\caption{Every basis set in two variables is a compound box.} \label{fig:2varfig}
\end{figure}
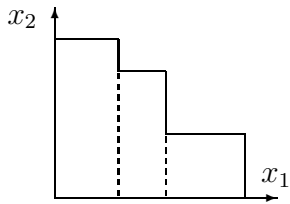

It is easy to verify (see Figure \ref{fig:2varfig}) that  \emph{every} basis set $\beta$ in two variables is a compound box; whence, Corollary \ref{cor:cpboxsmth} yields the following

\begin{cor} \label{cor:2varcase}
    Every basis set $\beta$ in two variables is smooth. \qedsymbol
\end{cor}

Haiman's lovely proof of this result (part of the proof of \cite[Proposition 2.4]{Haiman:CN-HS}) introduced the idea of arrow translation, and inspired the present paper.  Note that the $\gf$-basis of $M/M^2$ that he obtains is slightly different from ours; his basis arrows are typically not \emph{minimal} \sa\ arrows.  From the smoothness of the points $t_{\beta}$, Haiman deduces that $\Hn^n$ is everywhere nonsingular and irreducible (facts first proved by Fogarty \cite{Fogarty:AFAS}).  Corollary \ref{cor:cpboxsmth} can be viewed as a generalization of the two-variable smoothness phenomenon to higher dimensions.

\subsection{Example: The lexicographic point} \label{SS:lexpoint}

  The ``lexicographic point'' of $\Hilb^{p(z)}_{\mathbb{P}^r_{\gf}}$ is the point corresponding to the unique saturated lexicographic ideal $L$ such that $\gf[X_0, \dots, X_r]/L$ has Hilbert polynomial $p(z)$.  A.\ Reeves and M.\ Stillman  prove in general that the lexicographic point is a smooth point \cite{Reeves-Stillman-LexPt}.  In the case of a constant Hilbert polynomial $p(z)$ $=$ $n$, one checks that
\[
    L = (X_0,\ X_1,\ \dots, X_{r-2},\ X_{r-1}^{n}).
\]
Dehomogenizing with respect to the variable $X_r$, we obtain the ideal $I_{\beta}$ $\subseteq$ $\gf[x_0, \dots, x_{r-1}]$, where
\[
    \beta = \{1,\ x_{r-1},\ x_{r-1}^2,\dots,\ x_{r-1}^{n-1} \}.
\]
This is clearly a smooth basis set, since it is a special case of a box.

\subsection{Example: $\beta$ $=$ $\{1,\, x_1,\, x_2,\, x_1 x_2,\, x_3 \}$} \label{SS:spcpbox}
In Example \ref{SS:nonsmootheg} we considered the basis set
\[
    \{1,\ x_1,\ x_2,\ x_3  \},
\]
which is non-smooth.  If we add the monomial $x_1 x_2$, we obtain a compound box $\beta$ (see Figure \ref{fig:cpboxeg}); we will verify ``by hand'' that $\beta$ is smooth.
\begin{figure} 
\begin{picture}(374,84)
        \put(156,13){$x_1$}
        \put(122,0){$\mathbf{x_1^2}$}
        \put(190,26){$1$}
        \put(224,26){$x_2$}
        \put(258,26){$\mathbf{x_2^2}$}
        \put(188,12){${x_1 x_2}$}
        \put(190,50){$x_3$}
        \put(190,74){$\mathbf{x_3^2}$}
        \put(152,37){$\mathbf{x_1 x_3}$}
        \put(214,50){$\mathbf{x_2 x_3}$}
        \put(178,39){\vector(4,-1){40}}  
        \put(195,70){\vector(-1,-2){26}} 
        \put(164,33){\vector(0,-1){14}}  
        \put(227,47){\vector(0,-1){14}}  
        \put(190,22){\vector(-2,-1){60}} 
        \put(190,22){\vector(1,0){90}}   
        \put(190,22){\vector(0,1){70}}   
\end{picture}
\caption{The compound box $\beta$ $=$ $\{1, x_1, x_2, x_1 x_2, x_3 \}$, with the minimal generators of $I_{\beta}$ shown in boldface.  Note, for example, that the arrows $c^{(1,0,1)}_{(0,1,0)}$ and $c^{(0,0,2)}_{(1,0,0)}$ are translation-equivalent to $0$, and  $c^{(1,0,1)}_{(1,0,0)}$ $\sim$ $c^{(0,1,1)}_{(0,1,0)}$.} \label{fig:cpboxeg}
\end{figure}

The minimal generators of $I_{\beta}$ are the monomials
\[
    x_1^2,\ x_2^2,\ x_1 x_3,\ x_2 x_3,\ x_3^2;  
\]
whence, one has $5$ $\cdot$ $5$ $=$ $25$ minimal arrows, and these span $M/M^2$ by Lemma \ref{lem:msa}.  However, inspecting the minimal arrows, we find that:
\begin{itemize}
  \item Four are non-\sa\ arrows, all of which are translation-equivalent to $0$: $c^{(1,0,1)}_{(0,0,0)}$, $c^{(1,0,1)}_{(0,1,0)}$, $c^{(0,1,1)}_{(0,0,0)}$, and $c^{(0,1,1)}_{(1,0,0)}$.
  \item Five are \sa\ arrows that are translation-equivalent to $0$:  $c^{(2,0,0)}_{(0,0,0)}$, $c^{(0,2,0)}_{(0,0,0)}$, $c^{(0,0,2)}_{(0,0,0)}$, $c^{(0,0,2)}_{(1,0,0)}$, and $c^{(0,0,2)}_{(0,1,0)}$.
  \item Two are \sa\ arrows that are not translation-equivalent to $0$, but \emph{are} translation-equivalent to each other:   $c^{(1,0,1)}_{(1,0,0)}$ $\sim$ $c^{(0,1,1)}_{(0,1,0)}$.
\end{itemize}
This means that there are at most 15 non-trivial translation-equivalence classes of arrows available to span $M/M^2$, but $r$ $\cdot$ $n$ $=$ $3$ $\cdot$ $5$ = 15 is a lower bound on the $\gf$-dimension of the cotangent space, by Proposition \ref{prop:freedim}.  It follows that $M/M^2$ has $\gf$-dimension 15; whence, $\beta$ is smooth.


\section{Smooth basis sets in three variables are compound boxes} \label{S:threevars}

The main goal of this section is to prove Theorem \ref{thm:3varmain}, which states that a basis set $\beta$ in three variables is smooth if and only if $\beta$ is a compound box.  In the next section we will show by example that smooth basis sets in four or more variables need not be compound boxes.

\subsection{The main lemma} \label{sec:3varmainlem} 

Let $\beta$ be a basis set in the variables $x_1$, $x_2$, $x_3$, and $I_{\beta}$ the associated monomial ideal.  For 
\begin{equation} \label{E:3tups}
    (i,j) \in \{(1, 2), (1, 3), (2, 3)\}, 
\end{equation}
we write $G(i,j)$ for the set of minimal generators of $I_{\beta}$ that involve only variables in the set $\{x_i, x_j \}$.  For example, we always have that
\[
     x_i^{w_i} \in G(i, j)\ \text{ and }\ x_j^{w_j} \in G(i,j),
\]
so the number of elements 
\[
    |G(i,j)| \geq 2.
\]

\begin{lem} \label{lem:3varmainlem}

    If $|G(i,j)|$ $>$ $2$ for all of the ordered pairs $(i,j)$ in \emph{(\ref{E:3tups})}, then there exists a non-standard arrow for $\beta$ that is not translation-equivalent to $0$; consequently, $\beta$ is not smooth by Theorem \emph{\ref{thm:necsufcond}}.
    
\end{lem}

\subsection{Proof of Lemma \ref{lem:3varmainlem}} \label{SS:mainlempf}
    By the hypothesis, we may choose 
\begin{equation} \label{E:selmonoms}
  \begin{array}{rcccl}
    m_{(1,3)} & = & x_1^a x_3^b & \in & G(1,3),\\ 
    m_{(2,3)} & = & x_2^c x_3^d & \in & G(2,3),\\ 
    m_{(1,2)} & = & x_1^e x_2^f & \in & G(1,2)
  \end{array}
\end{equation}
such that all the exponents are positive and $m_{(1,3)}$ (resp.\ $m_{(2,3)}$) has maximal $x_1$-degree (resp.\ maximal $x_2$-degree) subject to the stated constraints; note that this implies that $m_{(1,3)}$ (resp.\ $m_{(2,3)}$) has minimal $x_3$-degree subject to the stated constraints.  We proceed to construct the desired non-standard arrow; there are two cases (see Figure \ref{fig:3vars}):

\begin{figure}
\begin{picture}(374,108)

        \put(30,0){\line(1,0){24}}
        \put(30,0){\line(0,1){48}}
        \put(54,0){\line(0,1){48}}
        \put(30,48){\line(1,0){24}}

        \put(30,48){\line(2,1){24}}
        \put(54,48){\line(2,1){24}}
        \put(54,0){\line(2,1){24}}
        \put(78,12){\line(0,1){48}}

        \put(54,60){\line(1,0){24}}
        \put(54,60){\line(0,1){12}}
        \put(54,72){\line(1,0){48}}
        \put(54,72){\line(2,1){24}}
        \put(78,84){\line(1,0){48}}
        \put(102,72){\line(2,1){24}}
        \put(102,72){\line(0,-1){60}}
        \put(78,12){\line(1,0){24}}

        \put(102,12){\line(2,1){12}}
        \put(114,18){\line(1,0){36}}
        \put(114,18){\line(0,1){24}}
        \put(114,42){\line(2,1){12}}
        \put(114,42){\line(1,0){36}}

        \put(126,48){\line(0,1){36}}
        \put(126,48){\line(1,0){36}}

        \put(150,18){\line(2,1){12}}
        \put(150,18){\line(0,1){24}}

        \put(150,42){\line(0,-1){24}}
        \put(150,42){\line(2,1){12}}
        \put(162,48){\line(0,-1){24}}

        \put(30,0){\vector(-2,-1){12}}
        \put(12,6){$x_1$}
        \put(162,24){\vector(1,0){12}}
        \put(168,30){$x_2$}
        \put(78,84){\vector(0,1){12}}
        \put(60,90){$x_3$}

        \put(54,63){\vector(4,-1){94}}
        \put(12,60){$m_{(1,3)}$}
        \put(42,62){\vector(1,0){10}}
        \put(138,60){$m_{(2,3)}$}
        \put(137,59){\vector(-1,-1){10}}
        \put(132,0){$m_{(1,2)}$}
        \put(126,6){\vector(-1,1){10}}

        \put(234,12){\line(0,1){36}}
        \put(234,12){\vector(-2,-1){12}}
        \put(208,12){$x_1$}
        \put(234,48){\line(1,0){36}}
        \put(234,48){\line(2,1){36}}
        \put(270,66){\line(0,1){24}}
        \put(282,96){\line(-2,-1){12}}
        \put(282,96){\line(1,0){54}} 
        \put(282,96){\vector(0,1){12}}
        \put(264,102){$x_3$}
        \put(270,90){\line(1,0){54}}
        \put(336,96){\line(0,-1){48}}
        \put(336,48){\line(-2,-1){24}}
        \put(336,96){\line(-2,-1){12}}
        \put(336,48){\line(1,0){12}}
        \put(348,48){\line(0,-1){12}}
        \put(348,36){\vector(1,0){12}}
        \put(354,42){$x_2$}
        \put(348,48){\line(-2,-1){24}}
        \put(324,36){\line(-1,0){12}}
        \put(312,36){\line(0,1){24}}
        \put(270,48){\line(2,1){24}}
        \put(294,60){\line(1,0){18}}
        \put(312,60){\line(2,1){12}}
        \put(324,66){\line(0,1){24}}
        \put(234,12){\line(1,0){36}}
        \put(270,66){\line(1,0){54}}
        \put(270,12){\line(0,1){36}}
        \put(270,12){\line(2,1){24}}
        \put(294,24){\line(0,1){36}}
        \put(294,24){\line(1,0){30}}
        \put(324,24){\line(0,1){12}}
        \put(324,24){\line(2,1){24}}

        \put(258,67){\vector(1,0){10}}
        \put(228,65){$m_{(1,3)}$}
        \put(348,60){$m_{(2,3)}$}
        \put(347,59){\vector(-1,-1){10}}
        \put(312,6){$m_{(1,2)}$}
        \put(306,12){\vector(-1,1){10}}
        \put(270,66){\vector(3,-1){68}}

\end{picture}
\caption{The two cases considered in the proof of Lemma \ref{lem:3varmainlem}.  In the left-hand figure, the monomial $m_{(1,3)}$ is dominant; the associated rigid non-standard arrow is also illustrated in Figure \ref{fig:nonstrigid}.  In the right-hand figure, none of the monomials $m_{(i,j)}$ are dominant.  The associated non-standard arrow is not translation-equivalent to $0$.} \label{fig:3vars}
\end{figure}
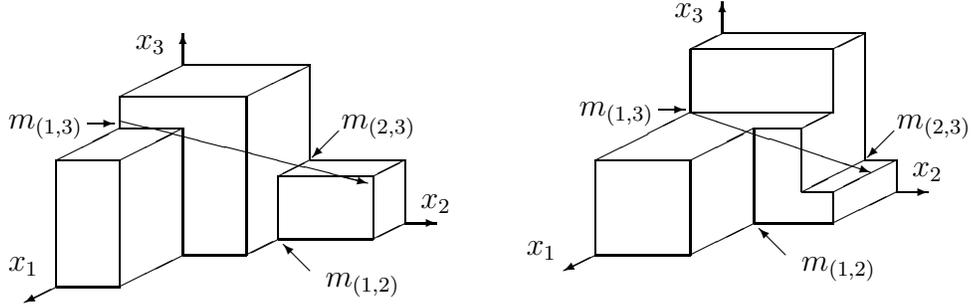

\subsubsection{Case 1:\ one of the monomials in \emph{(\ref{E:selmonoms})} is dominant} \label{SSS:pfcase1}

We say that $m_{(i,j)}$ is \textbf{dominant} provided that the degree of each of its constituent variables is greater than or equal to the degree of the same variable in the other monomial in which it appears.  For example, $m_{(1,3)}$ is dominant provided that (as shown in the left-hand portion of Figure \ref{fig:3vars})
\[
    a \geq e\ \text{ and }\ b \geq d.
\]
In this case, let $g$ = $\max(c,f)$, and note that $x_2^g$ is a basis monomial.  Starting at $x_2^g$, we can move to a maximal basis monomial $m^*$ $=$ $x_1^p x_2^q x_3^r$ by a sequence of degree-increasing steps.  Since $x_1^e x_2^f$ is a minimal generator of the ideal, we know that $p$ $<$ $e$.  Similarly, since $x_2^c x_3^d$ is a minimal generator, we have that $r$ $<$ $d$.  Then the rigid arrow $A$ with tail $m_{(1,3)}$ and head $m^*$ is also non-standard, since its vector 
\[
    (p,q,r) - (a,0,b)\ =\ (p-a, q, r-b)
\]
has the first and third coordinates negative:
\[
  \begin{array}{rcccccl}
       p-a & < & e-a & \leq & a - a & = & 0,\\
       r-b & < & d-b & \leq & b - b & = & 0.
  \end{array}
\]
Since $A$ is rigid, it is not translation-equivalent to $0$, as desired.

\subsubsection{Case 2:\ None of the monomials in \emph{(\ref{E:selmonoms})} is dominant} \label{SSS:pfcase2}

Writing out what this condition says, we find that:
\[
  \begin{array}{c}
    (\neg(a \geq e \wedge b \geq d))\ \wedge\ (\neg(c \geq f \wedge d \geq b))\ \wedge\ (\neg(e \geq a \wedge f \geq c))\\
\equiv\\
(a < e \vee b < d)\ \wedge\ (c < f \vee d < b)\ \wedge\ (e < a \vee f < c)\\
\equiv\\
(a < e \wedge f < c \wedge d < b)\  \vee\  (b < d \wedge c < f \wedge e < a).
  \end{array}
\]
Suppose the first alternative in the last line holds; that is, suppose (as shown in the right-hand portion of Figure \ref{fig:3vars}) that
\[
    a < e,\ f < c,\  \text{and } d < b.
\]
Note first of all that the monomial
\[
    m \ = \ x_2^{(w_2-1)} x_3^{(d-1)} \ \in \ \beta.
\]
If not, there would exist a minimal generator $m^{\prime}$ of $I_{\beta}$ such that 
\[
    m^{\prime}\, |\, m  \ \Rightarrow\ m^{\prime} \in G(2,3).
\]
It is clear that $m^{\prime}$ cannot equal either of the corner monomials $x_2^{w_2}$, $x_3^{w_3}$, so $m^{\prime}$ must involve both $x_2$ and $x_3$ nontrivially.  However, the $x_3$-degree of $m^{\prime}$ is $\leq$ $d-1$, which contradicts our choice of $m_{(2,3)}$ as having minimal $x_3$-degree $d$ among the members of $G(2,3)$ that involve two variables nontrivially.

Let $s$ $>$ $0$ denote the minimal exponent such that
\[
     x_1^s \cdot m \notin \beta;
\]
since
\[
     (x_1^e x_2^f = m_{(1,2)})\, |\, (x_1^e \cdot m)\ \Rightarrow \ (x_1^e \cdot m) \notin \beta,
\]
we have that 
\begin{equation} \label{E:sineq}
    s \leq e < w_1,
\end{equation}
where the second inequality holds because the minimal generator $x_1^{w_1}$ cannot divide the minimal generator $m_{(1,2)}$.  
Let 
\begin{equation} \label{E:alphadef}
   \alpha \    =\  \left\{ \begin{array}{l}
                              s-1, \text{ if } s \leq a,\\
                              a-1, \text{ if } s > a,
                           \end{array} \right. 
\end{equation}
and form the arrow 
\[
    A\ = \  c^{(a,0,b)}_{(\alpha,w_2-1, d-1)}
\]
with tail $m_{(1,3)}$ and head 
\[
  m_h\ =\ x_1^{\alpha} x_2^{w_2-1} x_3^{d-1}\ \in\ \beta\ \ (\text{ since } \alpha < s).
\]
Notice that $A$ is non-standard, since 
its vector
\[
     (\alpha, w_2 - 1, d-1) - (a, 0, b)\ = \ (\alpha - a, w_2-1, (d-1) - b)
\]
has negative first and third components (and non-negative second component $w_2 - 1$).

We assert that this arrow is not translation-equivalent to $0$.  If it were, the head would have to exit the first octant across either the hyperplane ($x_1$-degree $=$ $0$) or the hyperplane ($x_3$-degree $=$ $0$). In the former case, we would have to translate $A$ to an arrow $A^{\prime}$ of the same $x_1$-height ($=$ $a$), but with tail divisible by a minimal generator $\tilde{m}$ of $x_1$-degree $<$ $a$.  Since $\tilde{m}$ cannot divide $m_{(1,2)}$ and $m_{(1,3)}$, we must have that
\begin{equation} \label{E:mprime}
    \text{either } \ x_2\text{-degree}(\tilde{m}) > f \ \text{ or }\ x_3\text{-degree}(\tilde{m}) > b.
\end{equation}
However, we have that $A$ is ``rigid'' with respect to motion in the $x_2$-direction; that is, neither $A$ nor any of its translates $A^{\prime}$ can be moved in either the increasing or decreasing $x_2$-direction.  Indeed, an easy induction on the length of the path from $A$ to $A^{\prime}$ shows that 
\[
    x_2\text{-degree}(\text{head}(A^{\prime})) = w_2-1\ \text{ and }\ x_2\text{-degree}(\text{tail}(A^{\prime})) = 0,
\]
and clearly such an arrow cannot be translated in either the increasing or decreasing $x_2$-degree directions.  Furthermore, neither $A$ nor any of its translates $A^{\prime}$ has $x_3$-height greater than the initial value $b$, for we have already seen that the $x_2$-degree of the head of $A^{\prime}$ is invariantly $w_2-1$, and if the $x_3$-height of $A^{\prime}$ were to exceed $b$, then the $x_3$-degree of the head would be $\geq$ $d$, implying that the head would be divisible by $m_{(2,3)}$.  Therefore, we see that it is impossible to translate $A$ to $A^{\prime}$ with tail divisible by $\tilde{m}$ as in (\ref{E:mprime}), so we cannot move the head of $A$ across the hyperplane ($x_1$-degree $=$ $0$).

The only remaining possibility is to translate $A$ so that the head exits the first octant across the hyperplane ($x_3$-degree $=$ $0$).  This requires us to translate $A$ to an arrow $A^{\prime}$ of the same $x_3$-height ($=$ $b$), but with tail divisible by a minimal generator $\hat{m}$ of $x_3$-degree $<$ $b$.  Since the $x_2$-height of $A^{\prime}$ is invariantly $0$, we see that 
\[
  x_2\text{-degree}(\hat{m}) = 0\ \Rightarrow\  \hat{m} \in G(1,3)\ \Rightarrow\ \hat{m} = x_1^{w_1},
\]
where the last implication follows from our choice of $m_{(1,3)}$  as the monomial of minimal $x_3$-degree in $G(1,3)$ among those involving both $x_1$ and $x_3$ nontrivially.  In view of the constraints on the motion of $A$, we see that we would have to be able to translate $A$ a distance of $w_1-a$ units in the positive $x_1$-direction, but this motion would move the head to 
\[
    x_1^{\alpha + (w1-a)} x_2^{w_2-1} x_3^{d-1},
\]
which lies outside of $\beta$, because (recalling (\ref{E:sineq}) and (\ref{E:alphadef}))
\[
   \begin{array}{rcl}
       s \leq a & \Rightarrow & \alpha + (w_1-a) = (s-1)+(w_1-a) \geq s,\ \text{ and }\\
       s > a    & \Rightarrow & \alpha + (w_1-a) = (a-1) + (w_1-a) = w_1-1 \geq s. 
  \end{array}
\]
The required translation is therefore impossible; whence, $A$ is not translation-equivalent to $0$, and  the proof of Lemma \ref{lem:3varmainlem} is complete. \qedsymbol

\subsection{The main theorem} \label{SS:sbstv}

\begin{thm} \label{thm:3varmain}
     Let $\beta$ be a basis set in the variables $x_1$, $x_2$, $x_3$.  Then
\[
     \beta \text{ is smooth } \ \Leftrightarrow \ \beta \text{ is a compound box.}
\]
\end{thm}

\emph{Proof:}
    ($\Leftarrow$): Immediate from Corollary \ref{cor:cpboxsmth}.

    ($\Rightarrow$): We proceed by induction on $n$ $=$ $|\beta|$.

\emph{Base case:} $n$ $=$ $1$. The only basis set of cardinality $1$ is the box
\[ 
    \{1\}\ = \ \rbox(1,1,1),
\]
which is smooth by Proposition \ref{prop:boxsmth}, and trivially a compound box.

\emph{Inductive step:} Suppose that $|\beta|$ $=$ $n$, and that any smooth basis set of cardinality $<$ $n$ is a compound box.  Since $\beta$ is by hypothesis smooth, Lemma \ref{lem:3varmainlem} implies that for at least one of the pairs  $(i, j)$ in (\ref{E:3tups}), we have that 
\[
    G(i,j)\ =\ \{ x_i^{w_i}, x_j^{w_j} \};
\] 
that is, no minimal generator of $I_{\beta}$ exists that involves both $x_i$ and $x_j$ nontrivially, but does not involve the third variable $x_k$.  Let $h$ be the minimal $x_k$-degree among the minimal generators of $I_{\beta}$ that \emph{do} involve $x_k$.  If $h$ $=$ $w_k$, then
\[
    \beta = \rbox(w_1, w_2, w_3)
\] 
is a (compound) box, as desired.  If $h$ $<$ $w_k$, we let $\beta_t$ be the $x_k$-truncation of $\beta$ at height $h$ (\ref{E:truncdef}).  Then one sees easily that $\beta$ is obtained from $\beta_t$ by adding a box in the $x_k$-direction, as described in Section \ref{S:boxes}: more precisely, the added box has dimensions $w_i$ in the $x_i$-direction, $w_j$ in the $x_j$-direction, and $h$ in the $w_k$-direction.  It follows from Theorem \ref{thm:boxaddnmain} that $\beta_t$ is smooth, since $\beta$ is smooth by hypothesis.  The induction hypothesis now implies that $\beta_t$ is a compound box $\Rightarrow$ $\beta$ is a compound box, and we are done. \qedsymbol
\bigskip

     More can be gleaned from the preceding proof: suppose that $\beta$ is a basis set in the variables $x_1$, $x_2$, $x_3$ for which condition (a) of Theorem \ref{thm:necsufcond} holds; that is, every non-standard arrow for $\beta$ is translation-equivalent to $0$.  Then Lemma \ref{lem:3varmainlem} implies that for at least one of the pairs  $(i, j)$ in (\ref{E:3tups}), we have that 
\[
    G(i,j)\ =\ \{ x_i^{w_i}, x_j^{w_j} \};
\] 
arguing as in the preceding proof, we then see that $\beta$ is either a box or is the result of adding a box to a truncation $\beta_t$.  In the former case, $\beta$ is smooth, and in the latter case, $\beta_t$ also satisfies condition (a) of Theorem \ref{thm:necsufcond}, by Theorem \ref{thm:truncnsa}; therefore, induction yields that $\beta_t$ is smooth $\Rightarrow$ $\beta$ is smooth, by Theorem \ref{thm:boxaddnmain}.  Whence:
 
\begin{cor} \label{cor:dim3nsa}
    Let $\beta$ be a basis set in the variables $x_1$, $x_2$, $x_3$.  Then
\[
    \beta \text{ is smooth} \ \Leftrightarrow\ \left\{ \begin{array}{l}
                                                      \text{condition \emph{(a)} of Theorem \emph{\ref{thm:necsufcond}} holds; that is, every}\\ \text{non-standard arrow for $\beta$ is translation-equivalent to $0$.}\ \ \text{\qedsymbol}
                                               \end{array} \right.  
\]
\end{cor}

As of this writing, I do not know if this result extends to higher dimensions.


\section{The union of two boxes} \label{S:twoboxes}

To end this paper, we study one more family of smooth basis sets that does not consist entirely of compound boxes (in four or more variables): basis sets that are unions of two boxes.

\subsection{Notation} \label{SS:2bunot} 
     We will use the following notation throughout this section.  Let 
\[
  \begin{array}{rcl}
    \rbox_1 & = & \rbox(w_{1,1},\ w_{2,1},\ \dots,\ w_{r,1}),\\
    \rbox_2 & = & \rbox(w_{1,2},\ w_{2,2},\ \dots,\ w_{r,2})
  \end{array}
\]
be two boxes  (\ref{E:rbox}) in the variables $x_1$, \dots, $x_r$, and let 
\begin{equation} \label{E:2boxunion}
    \beta\ =\ \rbox_1 \cup \rbox_2.
\end{equation}
One checks easily that $\beta$ is a basis set.  Lacking inspiration, we call $\beta$ a \textbf{two-box union}.

\subsection{Minimal generators of $I_{\beta}$} \label{SS:U2mingen}

As usual, we write $w_i$ for the $x_i$-width of $\beta$; that is, $x_i^{w_i}$ is the corner (minimal) monomial divisible by $x_i$.  The following result is clear:

\begin{lem} \label{lem:2bcormon}
    The $x_i$-width $w_i$ of the two-box union $\beta$ \emph{(\ref{E:2boxunion})} is given by
\[
    w_i \ = \ \max(w_{i,1}, w_{i,2}),\ \ 1 \leq i \leq r.\ \ \text{\qedsymbol}
\]
\end{lem}

We now write the set of variables as a union
\[
    \{ x_1,\, x_2,\, \dots,\, x_r\}\ = \ V_1 \cup V_2 \cup V_3,
\]
where
\begin{equation} \label{E:Vs}
    \begin{array}{rcl}
        V_1 & = & \{x_j \mid w_{j,1} > w_{j,2}  \},\\
        V_2 & = & \{x_k \mid w_{k,1} < w_{k,2}  \},\\
        V_3 & = & \{x_{\ell} \mid w_{\ell,1} = w_{\ell,2}  \};
    \end{array}
\end{equation}
henceforth, we will use the subscripts $j$, $k$, and $\ell$ to denote membership in $V_1$, $V_2$, and $V_3$, respectively.  One sees easily that 
\[
  \begin{array}{l}
    V_1 = \emptyset \ \Rightarrow\ \rbox_1 \subseteq \rbox_2\ \Rightarrow\ \beta = \rbox_2,\\
    V_2 = \emptyset \ \Rightarrow\ \rbox_2 \subseteq \rbox_1\ \Rightarrow\ \beta = \rbox_1,
  \end{array}
\]
so the most interesting case is when both $V_1$ and $V_2$ are non-empty.

We have the following

\begin{lem} \label{lem:u2boxmingen}
    Let $\beta$ be a two-box union \emph{(\ref{E:2boxunion})}.  Then a minimal generator $m$ of the monomial ideal $I_{\beta}$ is either a corner monomial $x_i^{w_i}$ or a two-variable monomial of the form $x_j^{w_{j,2}} x_k^{w_{k,1}}$, with $x_j$ $\in$ $V_1$, $x_k$ $\in$ $V_2$ \emph{(\ref{E:Vs})}.
\end{lem}

\emph{Proof:}
Let 
\[
    m\ =\ x_1^{s_1} x_2^{s_2} \dots x_r^{s_r}
\]
be a minimal generator of $I_{\beta}$, and suppose that three or more of the exponents (say $s_1$, $s_2$, and $s_3$) are positive.  Then each of the monomials
\[
    x_1^{s_1-1} x_2^{s_2} x_3^{s_3} \dots x_r^{s_r},\ \ x_1^{s_1} x_2^{s_2-1} x_3^{s_3} \dots x_r^{s_r},\ \ x_1^{s_1} x_2^{s_2} x_3^{s_3-1} \dots x_r^{s_r} 
\]
belong to $\beta$, which implies that two of these monomials must belong to the same box ($\rbox_1$ or $\rbox_2$).  But then the least common multiple of the two monomials must also belong to this box (Lemma \ref{lem:boxlcm}); that is, $m$ $\in$ $\beta$, which is a contradiction.  We conclude that $\leq$ $2$ of the exponents $s_i$ can be positive.  If only one of the exponents $s_i$ is positive, then $m$ is the corner monomial $x_i^{(s_i\, =\, w_i)}$.  If two of the exponents (say $s_1$ and $s_2$) are positive, so that 
\[
    m \ = \ x_1^{s_1} x_2^{s_2},
\]
we again have that the monomials
\[
    x_1^{s_1-1} x_2^{s_2},\ \ x_1^{s_1} x_2^{s_2-1}
\]
belong to $\beta$; if both belonged to the same box, then we would arrive once more at the contradiction $m$ $\in$ $\beta$.  Therefore, the latter two monomials belong to different boxes, say
\[
     x_1^{s_1-1} x_2^{s_2} \in \rbox_2 \setminus \rbox_1,\ \ x_1^{s_1} x_2^{s_2-1} \in \rbox_1 \setminus \rbox_2;
\]
     consequently,
\[
  \begin{array}{l}
    x_1^{s_1-1} x_2^{s_2} \in \rbox_2 \ \Leftrightarrow \ s_1 - 1 < w_{1,2} \text{ and } s_2 < w_{2,2},\\
    x_1^{s_1} x_2^{s_2-1} \in \rbox_1 \ \Leftrightarrow \ s_1  < w_{1,1} \text{ and } s_2-1 < w_{2,1}, 
  \end{array}
\]
and
\[
  \begin{array}{l}
       x_1^{s_1-1} x_2^{s_2} \notin \rbox_1\ \Leftrightarrow\ s_1-1 \geq w_{1,1}\ \text{ or }\ s_2 \geq w_{2,1},\\
       x_1^{s_1} x_2^{s_2-1} \notin \rbox_2\ \Leftrightarrow\ s_1 \geq w_{1,2} \ \text{ or }\ s_2-1 \geq w_{2,2}.
  \end{array}
\]
Note that
\[  
     s_2  < w_{2,2}\ \Rightarrow\ s_2 - 1 < w_{2,2}\ \Rightarrow s_1 \geq w_{1,2};
\] 
whence,
\[
    (s_1 \geq w_{1,2} \text{ and } s_1-1 < w_{1,2})\ \Rightarrow s_1 = w_{1,2},
\]
and
\[
    w_{1,2} = s_1 < w_{1,1} \ \Rightarrow \ x_1 \in V_1,
\]
as desired.  A similar argument yields  
\[
    s_2 = w_{2,1} \ \text{ and } \ x_2 \in V_2,
\]
and the  proof is complete.
\qedsymbol

\subsection{Two-box unions are smooth basis sets}

Retaining the notation of Subsections \ref{SS:2bunot} and \ref{SS:U2mingen}, we begin with the following 

\begin{lem} \label{lem:2boxstarr} 
    Let $\beta$ be a two-box union, and let 
\[
    \xlist^{\dlist} = x_j^{w_{j,2}} x_k^{w_{k,1}},\ \ x_j \in V_1,\ x_k \in V_2,
\]
 be a two-variable minimal generator of $I_{\beta}$, as in Lemma \emph{\ref{lem:u2boxmingen}}.  If $c^{\dlist}_{\jlist}$ is an $x_j$- (resp.\ $x_k$-) \sa\ arrow for $\beta$, then the head of the arrow $\xlist^{\jlist}$ $\in$ $\rbox_2$ (resp.\ $\rbox_1$).
\end{lem}

\emph{Proof:}
    If $c^{\dlist}_{\jlist}$ is an $x_j$-standard arrow, then we have that
\[
    x_k\text{-degree}(\xlist^{\jlist}) \geq x_k\text{-degree}(\xlist^{\dlist}) = w_{k,1} \ \Rightarrow\  \xlist^{\jlist} \notin \rbox_1.
\]
Since $\xlist^{\jlist}$ $\in$ $\beta$, we must have that $\xlist^{\jlist}$ $\in$ $\rbox_2$, as asserted.  A similar argument applies in the case that $c^{\dlist}_{\jlist}$ is $x_k$-standard.
\qedsymbol

\begin{thm} \label{thm:U2boxsmth}
    A two-box union $\beta$ \emph{(\ref{E:2boxunion})} is a smooth basis set.
\end{thm}
    
\emph{Proof:} 
    We will show that the conditions (a) and (b) of Theorem \ref{thm:necsufcond} hold for $\beta$.

First suppose that $c^{\dlist}_{\jlist}$ is a non-standard arrow; we must show that $c^{\dlist}_{\jlist}$ is translation-equivalent to $0$.  We may assume, by Lemma \ref{lem:msa}, that $c^{\dlist}_{\jlist}$ is a minimal arrow; that is, its tail $\xlist^{\dlist}$ is a minimal generator of $I_{\beta}$.  Since every arrow with tail a corner monomial $x_i^{w_i}$ is \sa, Lemma \ref{lem:u2boxmingen} implies that
\[
    \xlist^{\dlist} = x_j^{w_{j,2}} x_k^{w_{k,1}},\ \ x_j \in V_1,\ x_k \in V_2.
\]
Without loss of generality, suppose that the head $\xlist^{\jlist}$ $\in$ $\rbox_1$.  Then we can translate the arrow in the direction of increasing $x_j$-degree to reach $c^{\dlist_1}_{\jlist_1}$, where 
\[
  \begin{array}{c}
    x_j\text{-degree}(\xlist^{\jlist_1}) = w_j-1 = w_{j,1}-1 \ \text{ and }\\ x_{i}\text{-degree}(\xlist^{\jlist_1}) = x_{i}\text{-degree}(\xlist^{\jlist}),\ \text{for } i \neq j.
  \end{array}
\]
Then 
\[
    x_j\text{-degree}(\xlist^{\dlist_1}) > x_j\text{-degree}(\xlist^{\jlist_1})\ \Rightarrow \  x_j\text{-degree}(\xlist^{\dlist_1}) \geq  w_j,
\]
so we may translate $c^{\dlist_1}_{\jlist_1}$ by degree-decreasing steps so that its tail approaches the corner monomial $x_j^{w_j}$. The head must eventually exit the first orthant, since a non-standard arrow cannot have a corner monomial as its tail.  Therefore, $c^{\dlist}_{\jlist}$ $\sim$ $0$, and condition (a) holds.

Now let $\lis$ be a standard bunch of arrows for $\beta$, and $c^{\dlist}_{\jlist}$ $\nsim$ $0$ be a \sa\ arrow for $\beta$.  We must show that there exists an arrow  $c^{\dlist^{\prime}}_{\jlist^{\prime}}$ $\in$ $\lis$ such that $c^{\dlist}_{\jlist}$ $\sim$ $c^{\dlist^{\prime}}_{\jlist^{\prime}}$.  We may assume that $c^{\dlist}_{\jlist}$ is a minimal $x_i$-\sa\ arrow that cannot be advanced, by assertion (c) of Lemma \ref{lem:shadadv}.  If the tail $\xlist^{\dlist}$ is the corner monomial $x_i^{w_i}$, we are done, since then
\[
    c^{\dlist}_{\jlist} \in \lis(i) \subseteq \lis
\]
by Corollary \ref{cor:crnrArrLem}.  Otherwise, by Lemma \ref{lem:u2boxmingen}, we have that
\[
    \xlist^{\dlist} = x_j^{w_{j,2}} x_k^{w_{k,1}},\ \ x_j \in V_1,\ x_k \in V_2,
\]
and we may assume without loss of generality that $x_i$ = $x_j$.  Let $v$ denote the $x_j$-offset of $c^{\dlist}_{\jlist}$, and let $c^{\dlist^{\prime}}_{\jlist^{\prime}}$ be the unique arrow 
in $\lis(j)$ such that 
\[
  x_j\text{-degree}(\xlist^{\jlist^{\prime}}) = x_j\text{-degree}(\xlist^{\jlist})\ \text{ and\ \  $x_j$-offset}(c^{\dlist^{\prime}}_{\jlist^{\prime}}) = v;
\]
the existence of $c^{\dlist^{\prime}}_{\jlist^{\prime}}$ is guaranteed by Corollary \ref{cor:degofffact}.  The tail $\xlist^{\dlist^{\prime}}$ is a minimal generator of $I_{\beta}$ that is divisible by $x_j$.  We proceed to show that $c^{\dlist}_{\jlist}$ $\sim$ $c^{\dlist^{\prime}}_{\jlist^{\prime}}$.

If $\xlist^{\dlist^{\prime}}$ $=$ $x_j^{w_j}$, then one sees easily that
\[
     \xlist^{\jlist^{\prime}}\ =\ \xlist^{\jlist}/x_k^{w_{k,1}};
\]
it is then apparent that $c^{\dlist^{\prime}}_{\jlist^{\prime}}$ can be translated $w_{k,1}$ steps in the direction of increasing $x_k$-degree to reach an arrow $c^{\dlist^{\prime}_1}_{\jlist}$, whose tail $\xlist^{\dlist^{\prime}_1}$ is divisible by $\xlist^{\dlist}$; therefore, $c^{\dlist^{\prime}_1}_{\jlist}$ and $c^{\dlist^{\prime}}_{\jlist^{\prime}}$ can be advanced, which is a contradiction.  Hence, Lemma \ref{lem:2bcormon} yields that 
\[
     \xlist^{\dlist^{\prime}} = x_j^{w_{j,2}} x_{k^{\prime}}^{w_{k^{\prime},1}},\ \ x_{k^{\prime}} \in V_2.
\]

We now know that $c^{\dlist}_{\jlist}$ and $c^{\dlist^{\prime}}_{\jlist^{\prime}}$ have the same length, $x_j$-offset, and $x_j$-height $=$ $w_{j,2}$. By Lemma \ref{lem:2boxstarr}, we know that  
\[
    \xlist^{\jlist},\ \xlist^{\jlist^{\prime}} \in \rbox_2\ \Rightarrow\ \lcm(\xlist^{\jlist},\ \xlist^{\jlist^{\prime}}) = \xlist^{\jlist^*} \in \rbox_2 \subseteq \beta,
\]
where the implication uses Lemma \ref{lem:boxlcm}.  
It is now clear that $c^{\dlist}_{\jlist}$ and $c^{\dlist^{\prime}}_{\jlist^{\prime}}$ can each be translated by degree-increasing steps to the same arrow $c^{\dlist^*}_{\jlist^*}$; whence, $c^{\dlist}_{\jlist}$ $\sim$ $c^{\dlist^{\prime}}_{\jlist^{\prime}}$, as desired.  Therefore, condition (b) holds, and the proof is complete.  
\qedsymbol

\begin{rem}
     Note that \emph{three}-box unions need not be smooth; for example, the left-hand basis set illustrated in Figure \ref{fig:3vars} is a three-box union, but is not smooth, since a non-standard rigid arrow exists.
\end{rem}

\subsection{Example: $\beta$ $=$ $\{1,\, x_1,\, x_2,\, x_1 x_2,\, x_3,\, x_4,\, x_3 x_4 \}$} \label{SS:2boxeg}
This basis set $\beta$ is the two-box union 
\[
    \beta \ = \ \rbox(2,2,1,1) \cup \rbox(1,1,2,2)\ \subseteq\ \gf[x_1, x_2, x_3, x_4]. 
\]
Therefore, by Theorem \ref{thm:U2boxsmth}, $\beta$ is a smooth basis set.  However, it is clear that $\beta$ is not a compound box, so the characterization of smooth basis sets in three variables given by Theorem \ref{thm:3varmain} does not extend to higher dimensions.


\bibliographystyle{amsplain}
\providecommand{\bysame}{\leavevmode\hbox to3em{\hrulefill}\thinspace}

\end{document}